\documentclass[10pt]{article}
\usepackage{latexsym}
\usepackage{amsmath}
\usepackage{cite}
\usepackage{amsthm,color,bm}
\usepackage{amssymb}
\usepackage{mathrsfs}
\usepackage{indentfirst}
\usepackage{times,cases}
\usepackage{lineno}
\usepackage{dsfont}
\usepackage{extarrows}

\usepackage{lscape}
\newtheorem{theorem}{\noindent Theorem}[section]
\newtheorem{proposition}[theorem]{\noindent Proposition}
\newtheorem{definition}[theorem]{\noindent Definition}
\newtheorem{lemma}[theorem]{\noindent Lemma}
\newtheorem{remark}[theorem]{\noindent Remark}
\newtheorem{corollary}[theorem]{\noindent Corollary}
\numberwithin{figure}{section}
\numberwithin{equation}{section}

\renewcommand{\theequation}{\thesection.\arabic{equation}}

\newcommand{\cA}{{\mathcal A}}
\newcommand{\cB}{{\mathcal B}}
\newcommand{\cC}{{\mathcal C}}
\newcommand{\cD}{{\mathcal D}}

\newcommand{\cF}{{\mathcal F}}

\newcommand{\cI}{{\mathcal I}}
\newcommand{\cK}{{\mathcal K}}

\newcommand{\cM}{{\mathcal M}}

\newcommand{\cP}{{\mathcal P}}
\newcommand{\cR}{{\mathcal R}}
\newcommand{\cS}{{\mathcal S}}
\newcommand{\cT}{{\mathcal T}}

\newcommand{\cX}{{\mathcal X}}

\newcommand{\fU}{\mathfrak{U}}

\newcommand{\sA}{{\mathscr A}}

\newcommand{\sD}{{\mathscr D}}

\newcommand{\sS}{{\mathscr S}}

\def\T{\mathbb{T}}

\def\R{\mathbb{R}}
\def\bE{\mathbb{E}}
\def\bP{\mathbb{P}}

\def\N{\mathbb{N}}
\def\1{\mathds{1}}
\newcommand{\di}{\mathrm{d}}
\newcommand{\me}{\mathrm{e}}

\def\disp{\displaystyle}

\def\bc{\begin{center}}
\def\ec{\end{center}}
\def\be{\begin{equation}}
\def\ee{\end{equation}}
\def\bea{\begin{eqnarray}}
\def\eea{\end{eqnarray}}
\def\ba{\begin{array}}
\def\ea{\end{array}}
\def\benu{\begin{enumerate}}
\def\eenu{\end{enumerate}}
\def\bt{\begin{theorem}}
\def\et{\end{theorem}}
\def\bl{\begin{lemma}}
\def\el{\end{lemma}}
\def\bco{\begin{corollary}}
\def\eco{\end{corollary}}
\def\bn{\begin{numcases}}
\def\en{\end{numcases}}
\def\br{\begin{remark}}
\def\er{\end{remark}}
\def\bd{\begin{definition}}
\def\ed{\end{definition}}
\def\bp{\begin{proposition}}
\def\ep{\end{proposition}}
\def\bo{\begin{proof}}
\def\eo{\end{proof}}
\def\bx{\begin{example}}
\def\ex{\end{example}}
\def\bal{\begin{align}}
\def\eal{\end{align}}

\def\pa{\partial}
\def\al{\alpha}
\def\De{\Delta} \def\de{\delta}
\def\na{\nabla}
\def\lam{\lambda} 

\def\ve{\varepsilon}
\def\sig{\sigma}
\def\vsig{\varsigma}
\def\vp{\varphi}
\def\w{\omega}\def\W{\Omega}
\def\Gam{\Gamma}

\def\~{\widetilde}

\def\ol{\overline}

\def\la{\leftarrow}
\def\ra{\rightarrow}

\def\Lra{\Leftrightarrow}

\def\8{\infty}
\def\X{\times}

\def\mb{\mbox}
\def\di{{\rm d}}
\def\me{{\rm e}}

\def\es{\emptyset}

\def\sm{\setminus}

\def\Hs{\hspace{0.8cm}}
\def\hs{\hspace{0.4cm}}
\def\Vs{\vskip10pt}
\def\vs{\vskip5pt}

\def\({\left(}
\def\){\right)}

\begin{document}
\pagewiselinenumbers


\begin{center}
    {\large \bf Global martingale and pathwise solutions and infinite regularity\\
    of invariant measures for a stochastic modified Swift-Hohenberg equation
    }
\vspace{0.5cm}\\
{Jintao Wang$^{*}$,\hs Xiaoqian Zhang,\hs Chunqiu Li}\\\vspace{0.3cm}

{\small Department of Mathematics, Wenzhou University, Wenzhou 325035, China}
\end{center}


\renewcommand{\theequation}{\arabic{section}.\arabic{equation}}
\numberwithin{equation}{section}


\begin{abstract}
We consider a 2D stochastic modified Swift-Hohenberg equations
with multiplicative noise and periodic boundary.
First, we establish the existence of local and global martingale and pathwise solutions
in the regular Sobolev space $H^{2m}$ for each $m\geqslant1$.
Associated with the unique global pathwise solution, we obtain a Markovian transition semigroup.
Then, we show the existence of invariant measures and ergodic invariant measures
for this Markovian semigroup on $H^{2m}$.
At last, we improve the regularity of the obtained invariant measures to $H^{2(m+1)}$.
With appropriate conditions on the diffusion coefficient, we can deduce the infinite regularity of
the invariant measures, which was conjectured by Glatt-Holtz \textit{et al.} in their situation \cite{GKV14}.
\vs

\vs


\end{abstract}

\vspace{-1 cm}

\footnote[0]{\hspace*{-7.4mm}
$^{*}$ Corresponding author.\\
E-mail address: wangjt@wzu.edu.cn (J.T. Wang); zhangxiaoqian@stu.wzu.edu.cn (X.Q. Zhang);
lichunqiu@wzu.edu.cn (C.Q. Li).}

\tableofcontents
\section{Introduction}
\label{s1}

In this paper, we investigate several stochastic problems,
including the local and global existence of martingale and pathwise solutions,
existence of ergodic invariant measures and infinite regularity of
the invariant measures,
for the following two-dimensional stochastic modified Swift-Hohenberg equation
(MSHE for short) with multiplicative noise and periodic boundary,
\bn{}
\di u+[\De^2u+2\De u+au+b|\na u|^2+u^3]\di t=\phi(u)\di W,\;\;~ t>0,\label{1.1}\\
u(x,y,t)=u(x+2\pi,y,t)=u(x,y+2\pi,t),\;\; ~t\geqslant0,\label{1.2}\\
u(x,y,0)=u_0(x,y),\;\;\label{1.3}
\en
where $u(x,y,t)$ is the unknown amplitude function, $(x,y)\in\R^2$,
$\De=\pa_{xx}+\pa_{yy}$, $a,\,b\in\R$
and $u_0$ is $2\pi$-periodic with respect to $x$ and $y$ respectively.
The stochastic term $\phi(u)\di W$ in \eqref{1.1} will be specified
in Subsection \ref{subs2.2}.

The Swift-Hohenberg equation was named after Swift and Hohenberg (see \cite{SH77}),
and arose from convective phenomena in the research of geophysical fluid flows
in atmosphere, oceans and earth's mantle.
It is closely contacted with nonlinear Navier-Stokes equations coupled
with the temperature equation.
The Swift-Hohenberg equation has played a significant role in different branches of physics,
ranging from hydrodynamics to nonlinear optics (see e.g. \cite{HS92,LMR75,LMN94,PM80}).
For this modified case, the cubic term $u^3$ in \eqref{1.1} is used as an approximation of
a nonlocal integral term (\cite{WLYJ21,WWL22}).
The modified term $b|\na u|^2$ in \eqref{1.1} comes from the various pattern formation phenomena
involving some kind of phase turbulence or phase transition (\cite{S77}),
which eliminates the symmetry $u\ra-u$.

There have been lots of research on the subject of MSHEs.
Roughly speaking, these works mainly include three aspects: attractors
(\cite{GGL17,LWZ20,WD17,XM15}) and the regularity (\cite{SZM10}),
bifurcations of solutions (\cite{C15,CHHL17,XG10})
and optimal control (\cite{DG12,S18}).
What is more, for the nonautonomous MSHE,
Wang, Yang and Duan presented in \cite{WYD20} a lower number of recurrent solutions
by topological methods (see more in \cite{LLW21,WL21,WLD16,WLD21});
Wang, Zhang and Zhao studied the existence of invariant measures and statistical solutions
in \cite{WZZ21}.
Due to the delicate impact on the deterministic cases in real world,
randomness has captured more and more concern over the development of evolution systems
recently.
It becomes centrally significant to take stochastic effects into account for mathematical models
of complex phenomena in engineering and science.
However, the study of stochastic MSHEs is still inadequate in the literature up to now,
especially that of existence and regularity of the invariant measures.

In physics, invariant measures have performed a significant role in the research of
turbulence (\cite{FMRT01}).
In the field of dynamical systems, invariant measures are to the measure space
for the phase space what invariant sets (even attractors) are to the phase space.
Essentially, invariant measures also illustrate the long-time dynamical behaviors for
a given dynamical system, but in view of measure.
Particularly, for stochastic systems, the solution process is random.
If we fix each random variable, the Wiener process $W$ would become a deterministic function.
Then we can also study the attractor and its properties
for the corresponding deterministic system.
This treatment seems like straying from the essence of randomness.
According to this observation, we are encouraged to consider the existence and regularity
of invariant measures for the stochastic dynamical systems.

In recent decades, the research of (ergodic) invariant measures has witnessed
the development of (stochastic) dynamical systems.
References \cite{CDJNZ,CDJZ,HJLY1,HJLY2,HJLY3,HJLY4,HJLY5,HJLY6} studied the limit
property or stochastic stability of some invariant measures for stochastic processes
as the noise vanishes or perturbs.
There have been more references concerning the existence of invariant measures,
such as \cite{BKR96,DaPZ14,EH01,EKZ17,GKV14,MSY16,WLYJ21,WZ23,WZC20,WSD05} and their references.
For the existence of invariant measures for stochastic partial differential equation,
it often needs to assume that the drift terms to be dissipative.
If the drift terms contains nonlinear ones, the nonlinear terms are often hoped to be
globally Lipschitz or dissipative.
In the book \cite{DaPZ14}, Da Prato and Zabczyk provided basic theory on
the existence of invariant measures for stochastic partial differential equations
with globally Lipschitz or dissipative nonlinearity.
Relatively, regularity of invariant measures has attracted little attention until now.
Bogachev, Krylov and R\"ockner provided in \cite{BKR96} a method to
improve the regularity of an invariant measure.
By adopting this idea, Glatt-Holtz et al. in \cite{GHV14} presented
another result, showing higher regularity of the invariant measures
for the three dimensional stochastic primitive equations.
The authors of \cite{GHV14} even conjecture that ``\textit{the invariant measures for
the 3D primitive equation are in fact supported on $\cC^{\8}$ function}"
with appropriate continuity and compatibility assumptions on the force.
\vs

Now, as presented in \eqref{1.1}, the modified nonlinear terms lack
the global Lipschitz continuity or dissipation,
and the lack brings new challenges to our discussion.
Moreover, what we study in this paper is the stochastic MSHE
with a general multiplicative noise, which is allowed to be nonlinear.
In this situation, we will not consider the mild solutions
(studied in \cite{LWZ20,WLYJ21,WSD05}),
but discuss the weak and strong solutions in the stochastic sense, saying,
the martingale and pathwise solutions, correspondingly, and establish their global existence.
With the unique global pathwise solution in hand,
we further explore the existence of ergodic invariant measures for the Markovian semigroup
associated with the solution in spaces of high regularity.
We also turn out that each invariant measure obtained is supported on $\cC^{\8}\cap H$
(defined in Subsection \ref{subs2.1}),
which is said to be the infinite regularity of the invariant measure,
under some appropriate assumption on the diffusion coefficient $\phi$.
This means that we have actually proved the conjecture of Glatt-Holtz et al. in our situation.

Since the equation \eqref{1.1} contains the modified nonlinearity, the existence of martingale
and pathwise solutions is far from holding naturally without a detailed argument.
We initiate our topics by introducing cut-off functions to treat the nonlinear terms,
and then apply the Galerkin approximating methods to the cut-off system.
By this setting, we are prompted to use the unique global existence of solutions (\cite{E13}) for
finite-dimensional stochastic differential equation with globally Lipschitz nonlinearity
to start our work.

To show the local existence of martingale and pathwise solutions in $H^{2m}$ (the Sobolev space
with periodic boundary given in Subsection \ref{subs2.1}) for each $m\geqslant1$,
we adopt the classical Galerkin method and Yamada-Wannabe theorem (see \cite{DGHT11,DZ20,GHV14})
for stochastic partial differential equations.
We actually conclude that the martingale solution (see Theorem \ref{th4.2})
or pathwise solution (see Theorem \ref{th5.3}) $u(t)$ satisfies
$$u(\cdot\wedge\tau)\in L^2(\W;\cC([0,\8);H^{2m})\cap L_{\rm loc}^2(0,\8;H^{2(m+1)}))$$
for some stopping time $\tau\geqslant 0$ or $\tau>0$,
where $\W$ is a probability space that is described in Subsection 2.2.
Here it is only required that $\phi$ is locally bounded (for martingale case)
or locally Lipschitz (for pathwise case).
As to the uniqueness of solutions, the usage of Skorohod embedding theorem
can not ensure the local uniqueness in the argument for existence of martingale solutions,
albeit with Lipschitz condition of $\phi$,
which, however, is sufficient to imply the local uniqueness in the pathwise case,
as the stochastic basis is chosen fixed in advance.
The uniqueness enables us to extend every local pathwise solution to a maximal one.
In particular, the maximal pathwise solution undergoes a blowup at the maximal living time.

For the global existence of martingale solutions, Dhariwal et al. in \cite{DJZ19} considered
a stochastic population cross-diffusion system and gave
the existence of global martingale solutions, in the vension that the solution depends on
the time interval $[0,T]$, where $T$ can be chosen arbitrarily.
In this paper, we consider a new version that the global martingale is defined on the entire
positive axis.
In order to obtain the existence of such global martingale solutions,
we need to assume $|b|<4$ to constrain the modified term and $\phi$ to be globally Lipschitz.
Although infeasible to ensure uniqueness for martingale solutions, the Lipschitz
condition makes the solutions coincide locally for the same finite-dimensional
approximating equation \eqref{3.2} with the same initial datum but different cut-off functions.
Based on this coincidence, we can then make good use of the diagonal method to pick
a subsequence of the approximating solutions, which converges all over $t\in[0,\8)$.
Thus, we successfully obtain a global martingale solution on $[0,+\8)$
for \eqref{1.1} -- \eqref{1.3}.
This process can be well applied to many other stochastic partial differential equations
to construct a global martingale solution.

For the global existence of pathwise solutions, we need the same assumptions as the case for
global martingale solutions, i.e., $|b|<4$ and the global Lipschitz condition for $\phi$.
Here we prove a general and simple conclusion (Lemma \ref{le5.5}),
similar to the Dominated Convergence Theorem,
but the limit function takes values at the infinity almost surely.
The global unique existence follows immediately from this conclusion and a uniform estimate
for the maximal pathwise solution given above by means of reduction to absurdity.

The unique existence of the global pathwise solution draws forth a Markovian transition
semigroup $\cP_t$ on the Borel $\sig$-algebra $\cB(H^{2m})$ of $H^{2m}$,
which is proved to be Feller by a process referring to \cite{WLYJ21}
with an additional condition \eqref{2.4} on $\phi$.
The existence of the invariant measures is a direct consequence of tightness for
the measure family $\{\nu_T\}_{T>0}$ with
$$\nu_T(\cdot)=\frac1T\int_0^T\cP_t(u,\cdot)\di t,$$
by Krylov-Bogoliubov existence theorem and Prokhorov's theorem (see \cite{DaPZ14}).
To guarantee the existence of ergodic measures, we utilize the Krein-Milman Theorem
(\cite{Conway90}) to find an extreme point of the set of invariant measures,
and this extreme point is exactly an ergodic measure.
Again we necessarily show the tightness of the set of all invariant measures,
which is assured in nature by the uniform integrability of $\|u\|_{2m}^2$,
the square of $H^{2m}$-norm of $u$, over $H^{2m}$ with respect to each invariant measure.
This uniform integrability is not an easy deduction.
Thanks to Song, Zhang and Ma's work \cite{SZM10}, in which they proved that the MSHE
possesses a global attractor in the Sobolev space $H^k$ for each $k\geqslant0$ by using
iteration procedure, we can as well use iteration procedure or induction method
to raise the order $k\in[0,m]$ one by one such that $\|u\|_{2k}^2$ is uniformly
integrable over $H^{2m}$ with respect to each fixed invariant measure.

For the rise of regularity of invariant measures, we consider to use the a posteriori method
by the idea propounded in \cite{BKR96}, and first show that $\|u\|^{2}_{2(m+1)}$ is uniformly
integrable over $H^{2m}$ with respect to each invariant measure under the same condition that
ensures the existence of invariant measures on $H^{2m}$.
This adds one to the regularity of invariant measures.
In this way, by assuming an appropriate stronger condition \eqref{2.5}
than \eqref{2.3} and \eqref{2.4},
we can eventually obtain the infinite regularity of the invariant measures.
\vs

The rest part of this paper is arranged as follows.
Section \ref{s2} provides the basic concepts, notions, settings and theorems that will be
used frequently in the following sections.
In Section \ref{s3} we construct the cut-off system and Galerkin scheme of the original
MSHE and give a uniform estimate for the solutions of Galerkin approximating equations.
In Section \ref{s4} we are devoted to proving the existence of local and global martingale
solutions.
And Section \ref{s5} is for the demonstration of existence of local and
global pathwise solutions.
In Section \ref{s6}, we prove the existence of invariant measures and then
ergodic invariant measures.
In Section \ref{s7}, the regularity of the invariant measures is discussed and improved
to the infinity under some appropriate conditions.
Section \ref{s8} is about the summary and remarks for this work.

\section{Preliminaries}
\label{s2}

We first present some basic notations and properties that will be used frequently in this paper.
\subsection{On spaces and operators}\label{subs2.1}

For metric spaces $X$ and $Y$, we conventionally denote by $\cC(X,Y)$ ($\cC_{\rm b}(X,Y)$)
the collection of continuous (and bounded) functionals from $X$ to $Y$.
When $Y=\R$, we simply use $\cC(X)$ ($\cC_{\rm b}(X)$)
to represent $\cC(X,\R)$ ($\cC_{\rm b}(X,\R)$).

Let $I$ denote the interval $[0,2\pi]$.
The periodic boundary condition \eqref{1.2} prompts us to study
the modified Swift-Hohenberg problem \eqref{1.1} -- \eqref{1.3} on the 2-dimensional torus $\T^2$,
saying, the quotient space $I\X I/\sim$,
where $\sim$ is the equivalence relation such that
$$(x_1,y_1)\sim(x_2,y_2)\hs\mb{if and only if}\hs (x_1,y_1)=(x_2,y_2),$$
$$\mb{ or }x_1=x_2\mb{ and }y_1,y_2\in\{0,2\pi\}\mb{ or }x_1,x_2\in\{0,2\pi\}\mb{ and }y_1=y_2.$$
This indicates that $(0,y)\sim(2\pi,y)$, $(x,0)\sim(x,2\pi)$ for each $x,y\in[0,2\pi]$ and
specifically $(0,0)\sim(0,2\pi)\sim(2\pi,0)\sim(2\pi,2\pi)$.
Hence each function $u$ on $\T^2$ is defined on $I\X I$
with $u(0,y)=u(2\pi,y)$ and $u(x,0)=u(x,2\pi)$.
By this setting, the space $L^2(\T^2)$ is a Hilbert space with the scalar product
$$(u,v):=\int_0^{2\pi}\int_0^{2\pi}u(x,y)v(x,y)\di x\di y,$$
and the norm denoted by $\|\cdot\|$.

Observe that (see also \cite{XG10}) the eigenvalues of $-\De:\sD(-\De)\subset L^2(\T^2)\ra L^2(\T^2)$ are
nonnegative integers $n$ that satisfies $n=k^2+l^2$ ($k,l\in\N$),
with the corresponding eigenfunctions
$$\cos(kx\pm ly),\hs\sin(kx\pm ly),\hs\cos(lx\pm ky)\hs\mb{and}\hs\sin(lx\pm ky),$$
for all possible $k,l\in\N$.
We give an order to these eigenvalues and denote them by new notations
$\{\lam_i\}_{i\in\N}$ as follows
$$
0=\lam_0<\lam_1\leqslant\lam_2\leqslant\lam_3\leqslant\cdots\leqslant\lam_i
\leqslant\cdots\ra\8,
$$
with the unit eigenfunction of $\lam_i$ denoted by $w_i$,
where we have taken the multiplicities into account.
In this way $\{w_i\}_{i\in\N}$ constitutes an orthonormal basis of $L^2(\T^2)$.
The Sobolev space $H^\al(\T^2)$ ($\al\in\R$) is the Banach space such that
$$
u=\sum_{i\geqslant0}a_iw_i\in H^{\al}(\T^2)\hs\Lra\hs
\sum_{i\geqslant0}(1+\lam_i)^\al|a_i|^2<+\8,
$$
which gives $u\in H^{\al}(\T^2)$ a norm
$$\|u\|_{\al}:=\(\sum_{i\geqslant0}(1+\lam_i)^\al|a_i|^2\)^{\frac12}.$$
Particularly, if $m$ is an integer, we also have
\be\label{2.1}\|(-\De+1)^{m}u\|=\(\sum_{i\geqslant0}(1+\lam_i)^{2m}|a_i|^2\)^{\frac12}
=\|u\|_{2m}.\ee
This also makes $H^\al(\T^2)$ a Hilbert space for each $\al\in\R$ such that
$$
(u,v)_{\al}:=(u,v)_{H^\al}=\sum_{i\geqslant0}(1+\lam_i)^\al a_ib_i\
\hs\mb{for }u=\sum_{i\geqslant0}a_iw_i,\,v=\sum_{i\geqslant0}b_iw_i.
$$
It is well known that the embedding of $H^{\al_2}$ into $H^{\al_1}$ is dense and compact
for all $\al_2>\al_1$.
Note that $L^2(\T^2)=H^0(\T^2)$.
Denote also the norm of $L^p:=L^p(\T^2)$ by $|\cdot|_p$ for all $2<p\leqslant\8$.
For notational simplicity, we denote
$$H:=L^2(\T^2),\hs H^\al:=H^\al(\T^2),\hs L^p:=L^p(\T^2),\hs W^{\al,p}:=W^{\al,p}(\T^2)$$
for $\al\in\R$ and $p\in[2,\8]$, with norms of the latter two spaces denoted
by $|\cdot|_p$ and $\|\cdot\|_{\al,p}$, respectively.
Here we need to keep in mind that $L^2(\T^2)=L^2(I\X I)$,
$W^{\al,p}(\T^2)=W^{\al,p}(I\X I)$ and $W^{\al,p}$ means the Sobolev space
with the usual Sobolev norm, which is equivalent to the norm of $H^\al$ when $p=2$.
We also use $\cC^{\8}$ to denote the set of all functions defined on $\T^2$ with derivatives
of arbitrary order.

Let
$$A=-\De+1:\sD(A)\subset H\ra H.$$
We know the eigenvalues of $A$ are $\{\lam_i+1\}_{i\in\N}$ with the eigenfunctions
$\{w_i\}_{i\geqslant0}$, correspondingly.
Recalling the basic knowledge of the fractional power of
sectorial operator (see \cite{Hen81,LW18}),
we can define the fractional power $A^{\al}$ of $A$ as
$$
A^{\al}u=\sum_{i\geqslant0}(\lam_i+1)^{\al}(u,w_i)w_i,\hs
\mb{for all }u\in H^{2\al}\mb{ and }\al\in\R.
$$
We can deduce from \eqref{2.1} that $\sD(A^{\al})=H^{2\al}$ and $\|A^{\al}u\|=\|u\|_{2\al}$
for all $u\in H^{2\al}$.
\vs

In process of completing our work, it is necessary to recall some spaces of fractional (in time)
derivative and some compact embedding results (see \cite{DGHT11,FG95}).
Let $X$ be a separable Hilbert space with the norm denoted by $\|\cdot\|_X$.
For fixed $p>1$ and $\al\in(0,1)$, define
$$
W^{\al,p}(0,T;X):=\left\{u\in L^p(0,T;X):\,
\int_0^T\int_0^T\frac{\|u(s)-u(\sig)\|_X^p}{|s-\sig|^{1+\al p}}\di s\di\sig<\8\right\}
$$
with the norm
$$
\|u\|_{W^{\al,p}(0,T;X)}^p:=\int_0^T\|u(s)\|_X^p\di s
+\int_0^T\int_0^T\frac{\|u(s)-u(\sig)\|_X^p}{|s-\sig|^{1+\al p}}\di s\di\sig.
$$

For the case when $\al=1$, we take
$$W^{1,p}(0,T;X):=\{u\in L^p(0,T;X):\frac{\di u}{\di t}\in L^p(0,T;X)\}$$
to be the classical Sobolev space with its usual norm
$$
\|u\|_{W^{1,p}(0,T;X)}^p:=\int_0^T\(\|u(s)\|_X^p+\left\|\frac{\di u(s)}{\di s}\right\|_X^p\)\di s.
$$
Note that for $\al\in(0,1)$,
\be\label{2.2}
W^{1,p}(0,T;X)\subset W^{\al,p}(0,T;X)\mb{ and }\|u\|_{W^{\al,p}(0,T;X)}
\leqslant\|u\|_{W^{1,p}(0,T;X)}.
\ee
With these settings, we have the compact embeddings below (\cite{DGHT11,FG95}).

\bl\label{le2.1}
\benu\item[(i)] Suppose that $X_2\supset X_0\supset X_1$ are Banach spaces
with $X_2$ and $X_1$ reflexive, and the embedding of $X_1$ into $X_0$ is compact.
Then for each $1<p<\8$ and $0<\al<1$, the embedding
$$L^p(0,T;X_1)\cap W^{\al,p}(0,T;X_2)\subset\subset L^p(0,T;X_0)$$
is compact (the notation $\subset\subset$ is used to denote the compact embedding).
\item[(ii)] Suppose that $Y_0\supset Y$ are Banach spaces with $Y$ compactly embedded in $Y_0$.
Let $\al\in(0,1]$ and $p\in(1,\8)$ be such that $\al p>1$ then
$$W^{\al,p}(0,T;Y)\subset\subset\cC([0,T];Y_0)$$
and the embedding is compact.
\eenu
\el

\subsection{On stochastic framework}\label{subs2.2}

As to determine the stochastic term $\phi(u)\di W$ in \eqref{1.1},
we recall some basic knowledge of stochastic analysis (see \cite{DaPZ14} for more details).

Fix a stochastic basis $\cS:=(\W,\cF,\{\cF_t\}_{t\geqslant0},\bP,\{W_i\}_{i\geqslant0})$,
that is, a filtered probability space $\W$ with a sequence $\{W_i\}_{i\geqslant0}$ of independent
standard one-dimensional Brownian motions relative to $\cF_t$.
For the sake of avoiding unnecessary complications,
we may assume that $\cF_t$ is complete and right continuous (\cite{DaPZ14}).
Fix a separable Hilbert space $\fU$ with an associated orthonormal basis $\{e_i\}_{i\geqslant1}$.
We may formally define $W$ as
$$W=\sum_{i\geqslant0}W_ie_i,$$
which is known as a ``cylindrical Brownian motion" over $\fU$.

Given another Hilbert space $X$ with its associated inner product denoted by $(\cdot,\cdot)_X$.
Let $L_2(\fU,X)$ be the separable Hilbert space of all
Hilbert-Schmidt operators (\cite[Appendix C]{DaPZ14}) from $\fU$ to $X$ with the inner product
$$
(\psi_1,\psi_2)_2=\sum_{i\geqslant0}(\psi_1e_i,\psi_2e_i)_X,\hs
\mb{for all }\psi_1,\psi_2\in L_2(\fU,X).
$$

To impose some conditions on the stochastic term, we introduce some new notations.
Given each pair of Banach spaces $X$ and $Y$ with $X\subset L^\8$,
we denote by ${\rm Bnd}_{u,{\rm loc}}(X,Y)$, the collection of
all continuous mappings $\vp:X\X[0,\8)\ra Y$ so that
$$\|\vp(x,t)\|_Y\leqslant \kappa(|x|_{\8})(1+\|x\|_X),\hs x\in X,\,t\geqslant 0,$$
where $\kappa:\R^+\ra(0,\8)$ is an increasing and locally bounded function, independently of $t$.
If, in addition,
$$\|\vp(x,t)-\vp(y,t)\|_Y\leqslant \kappa(|x|_\8+|y|_\8)\|x-y\|_X,$$
we say $\vp$ is in ${\rm Lip}_{u,{\rm loc}}(X,Y)$.
If the coefficient $\kappa$ involved above is a constant, i.e., independent of $t$, $x$ and $y$,
we denote ${\rm Bnd}_u(X,Y)$ and ${\rm Lip}_u(X,Y)$ instead, correspondingly.
Obviously ${\rm Lip}_{u,{\rm loc}}(X,Y)\subset{\rm Bnd}_{u,{\rm loc}}(X,Y)$
and ${\rm Lip}_u(X,Y)\subset {\rm Bnd}_u(X,Y)$.

Thus the stochastic term in \eqref{1.1} can be written as the formal expansion
$$\phi(u)\di W=\sum_{i\geqslant0}\phi_i(u)\di W_i,\hs\mb{with}\hs\phi_i(u)=\phi(u)e_i.$$
Let $m\geqslant1$.
For the global existence of martingale and pathwise solutions in Sections \ref{s4} and \ref{s5},
we assume that the mapping $\phi$ satisfies
\be\label{2.3}
\phi\in {\rm Lip}_{u}(H,L_2(\fU,H))\cap{\rm Lip}_{u}(H^{2m},L_2(\fU,H^{2m})).
\ee
For the existence of ergodic invariant measures in Section \ref{s6},
we assume further that
\be\label{2.4}
\phi\in {\rm Lip}_{u}(H^{2(m+1)},L_2(\fU,H^{2(m+1)})).
\ee
For the infinite regularity of the invariant measures discussed in Section \ref{s7},
we need to require that
\be\label{2.5}
\phi\in \bigcap_{m=0}^{\8}{\rm Lip}_{u}(H^{2m},L_2(\fU,H^{2m})).
\ee
such as $\phi(u)\in L_2(\fU,L^2)$, $u\in L^2$,
(we use $(\cdot,\cdot)_{\fU}$ the inner product of $\fU$)
so that
$$\phi(u)=k\sum_{i=1}^{\8}(u,w_i)w_i\otimes e_i\hs\mb{and}\hs
\phi(u)\circ W=k\sum_{i=1}^{\8}(u,w_i)(e_i,W)_{\fU}w_i,$$
for $W\in\fU$, where $\otimes$ is the tensor product between Hilbert spaces.

We also assume additionally
\be\label{2.6}
\phi\in{\rm Bnd}_{u,{\rm loc}}(H^{2m},L_2(\fU,H^{2m})).
\ee
or
\be\label{2.7}
\phi\in{\rm Lip}_{u,{\rm loc}}(H^{2m},L_2(\fU,H^{2m})),
\ee
for local consequences in Sections \ref{s4} and \ref{s5}.

Given an arbitrary separable Hilbert space $X$,
let
$$\vp\in L^2(\W,L_{\rm loc}^2([0,\8),L_2(\fU,X)))$$
be an $X$-valued predictable process.
The Burkholder-Davis-Gundy inequality (BDG inequality for short) holds in the following form,
$$
\bE\sup_{t\in[0,T]}\left\|\int_0^t\vp\di W\right\|_X^p\leqslant
c\bE\(\int_0^T\|\vp\|_{L_2(\fU,X)}^2\di t\)^{\frac{p}2},\hs\mb{for all }p\geqslant1,
$$
where $c$ only depends upon $p$.
Applying a variation of BDG inequality, for all $p\geqslant2$ and $\al\in[0,1/2)$,
we also have (\cite{DGHT11,FG95})
\be\label{2.8}
\bE\left\|\int_0^t\vp\di W\right\|_{W^{\al,p}(0,T;X)}^p\leqslant
c\bE\int_0^T\|\vp\|_{L_2(\fU,X)}^p\di t,
\ee
for all $X$-valued predictable $\vp\in L^2(\W,L_{\rm loc}^p([0,\8),L_2(\fU,X)))$.

We consider the initial condition $u_0$ to be random in general
if there is no special instructions.
Let $m\geqslant1$ and $p\geqslant2$.
When we study the case of martingale solutions (the definition of martingale and pathwise solutions mentioned below can be found in Subsection 2.3), the stochastic basis is unknown for the problem,
and hence we are only able to specify $u_0$ as an initial probability measure $\mu_0$
on $H^{2m}$, i.e., $\mu_0(\cdot)=\bP(u_0\in\cdot)$, such that
\be\int_{H^{2m}}\|u\|_{2m}^q\di\mu_0(u)<+\8,\mb{ with }q\geqslant p.\label{2.9}\ee
For the global existence of martingale solutions, we further assume that
\be\int_{H^{2m}}\(\|u\|_{2m}^{q}+\|u\|^{q'}\)\mu_0(u)<+\8,\mb{ with }q\geqslant p
\mb{ and }q'\geqslant(2m+3)p.\label{2.10}\ee
For the case of pathwise solutions where the stochastic basis $\cS$ is fixed, we assume that,
relative to this basis, $u_0$ is an $H^{2m}$-valued random variable and $\cF_0$ measurable.
For some auxiliary results during the procedure, we assume further that
\be\label{2.11}
\bE\|u_0\|_{2m}^{q}<+\8,\mb{ with }q\geqslant p.
\ee
Moreover, for the existence of global solutions and invariant measures,
we sometimes also need to assume that $u_0$ satisfies
\be\label{2.12}\bE\(\|u_0\|_{2m}^{q}+\|u_0\|^{q'}\)<+\8\mb{ with }q\geqslant p\mb{ and }
q'\geqslant(2m+3)p,\ee
where $q'$ will be specified when necessary in the sequel.

\subsection{On martingale and pathwise solutions}

In this subsection, we give the definitions of martingale and pathwise solutions
for the stochastic modified Swift-Hohenberg problem \eqref{1.1} -- \eqref{1.3}
with multiplicative noise.
By the settings above, the original problem \eqref{1.1} -- \eqref{1.3} can be rewritten as
\be\label{2.13}\di u+[A^2u+f(u)]\di t=\phi(u)\di W,\;\;~ t>0;\hs u(0)=u_0,\ee
where $f$ is given by
$$f(u):=(a+3)u-4Au+b|\na u|^2+u^3.$$
We present the definitions of local and global solutions of \eqref{2.13} in
both martingale and pathwise senses below.

\bd\label{de2.1}
Suppose $\mu_0$ is a probability measure on $H^{2m}$ ($m\geqslant1$) satisfying \eqref{2.9}.
Assume that $\phi$ satisfies \eqref{2.6} or \eqref{2.3}.
\benu\item[(i)] A triple $(\cS,u,\tau)$ is a \textbf{local martingale solution}
if
$$\cS=(\W,\cF,\{\cF_t\}_{t\geqslant0},\bP,W)$$
is a stochastic basis,
$\tau$ is a stopping time relative to $\cF_t$ and
$u(\cdot\wedge\tau):\W\X[0,\8)\ra H^{2m}$ is
an $\cF_t$ adapted process such that
\be\label{2.14}u(\cdot\wedge\tau)\in
L^2(\W;\cC([0,\8);H^{2m})\cap L_{\rm loc}^2(0,\8;H^{2(m+1)})),\ee
the law of $u(0)$ is $\mu_0$ and $u$ satisfies for almost every $(t,\w)\in[0,\8)\X\W$,
\be\label{2.15}
u(t\wedge\tau)+\int_0^{t\wedge\tau}[A^2u+f(u)]\di s=u(0)+\int_0^{t\wedge\tau}\phi(u)\di W
\ee
with the equality understood in $H^{2(m-1)}$.
\item[(ii)] We say that the martingale solution $(\cS,u,\tau)$ is \textbf{global} if $\tau=\8$,
for almost surely $\w\in\W$.
\eenu\ed

We recall a convergence lemma for stochastic integrals (see \cite[Lemma 2.1]{DGHT11}).
\bl\label{le2.2}
Let $(\W,\cF,\bP)$ be a fixed probability space, $X$ a separable Hilbert space.
Consider a sequence of stochastic bases $\sS_n=(\W,\cF,\{\cF_t^n\}_{t\geqslant0},\bP,W^n)$,
that is a sequence so that each $W^n$ is cylindrical Brownian motion (over $\fU$)
with respect to $\cF_t^n$.
Assume that $\{G^n\}_{n\geqslant1}$ are a collection of $X$-valued $\cF_t^n$ predictable processes
such that $G^n\in L^2(0,T;L_2(\fU,X))$ almost surely.
Finally consider $\sS=(\W,\cF,\{\cF_t\}_{t\geqslant0},\bP,W)$ and $G\in L^2(0,T;L_2(\fU,X))$,
which is $\cF_t$ predictable.
If
$$W^n\ra W\mb{ in probability in }\cC([0,T];\fU),$$
$$ G^n\ra G\mb{ in probability in }L^2(0,T;L_2(\fU,X)),$$
then
$$\int_0^t G^n\di W^n\ra\int_0^t G\di W\mb{ in probability in }L^2(0,T;X).$$
\el

\bd\label{de2.2}
Let $\cS=(\W,\cF,\{\cF_t\}_{t\geqslant0},\bP,W)$ be a stochastic basis.
Suppose that $u_0$ is an $H^{2m}$-valued ($m\geqslant1$) random variable (relative to $\cS$)
and $\cF_0$-measurable and $\phi$ satisfies \eqref{2.7} or \eqref{2.3}.
\benu\item[(i)]
A pair $(u,\tau)$ is a \textbf{local pathwise solution} of \eqref{2.13}
if $\tau$ is a strictly positive stopping time and $u:\W\X\R^+\ra H^{2m}$
is an $\cF_t$-adapt process satisfying \eqref{2.14} and \eqref{2.15}
for each $t\geqslant0$ and $\bP$-almost surely, in $H^{2m}$.
\item[(ii)]
Pathwise solutions of \eqref{2.13} are called (pathwise) \textbf{unique},
if given two pathwise solutions $(u,\tau)$ and $(u',\tau')$
such that $u(0)=u'(0)$ on a subset $\W_0$ of $\W$, then
$$\bP((u(t)-u'(t))\chi_{\W_0}=0,\mb{ for all }t\in[0,\tau\wedge\tau'])=1,$$
where $\chi_A$ is the characteristic function of $A$, i.e., $\chi_A(x)=1$, for $x\in A$
and $\chi_A(x)=0$, for $x\notin A$.
\item[(iii)]
Suppose that $\{\tau_n\}_{n\geqslant1}$ is a strictly increasing sequence of stopping times
converging to a (possibly infinite) stopping time $\tau_0$ and
assume that $u$ is a predictable process in $H$.
We say that the triple $(u,\tau_0,\{\tau_n\}_{n\geqslant1})$ is a \textbf{maximal pathwise solution}
if
$(U,\tau_n)$ is a local pathwise solution for each $n$ and
$$\sup_{s\in[0,\tau_0]}\|u(s)\|_{2}^2+\int_0^{\tau_0}\|A^2u(s)\|^2\di s=\8$$
almost surely on the set $\{\tau_0<\8\}$. If, moreover, as $n$ tends to the infinity,
\be\label{2.16}\sup_{s\in[0,\tau_n]}\|u(s)\|_{2}^2+\int_0^{\tau_n}\|A^2u(s)\|^2\di s\ra+\8,\ee
for almost all $\w$ with $\{\tau_0<\8\}$, then we say that the sequence announces any finite time
\textbf{blowup}.
\item[(iv)] If $(u,\tau_0)$ is a maximal pathwise solution and $\tau_0=\8$ almost surely,
then we say that the solution is \textbf{global}.
\eenu
\ed

In order to obtain pathwise solutions from martingale ones, we adopt a convergence conclusion
in probability given in \cite{GK96}.
Let $\{Y_n\}_{n\geqslant 0}$ be a sequence of $X$-valued random variables on a probability space
$(\W,\cF,\bP)$.
Let $\{\mu_{n,m}\}_{n,m\geqslant0}$ be the collection of joint laws of $\{Y_n\}_{n\geqslant0}$,
that is
$$\mu_{n,m}(E):=\bP((Y_n,Y_m)\in E),\hs E\in\cB(X\X X).$$
The convergence conclusion reads as follows.
\bp\label{pro2.1}
A sequence of $X$-valued random variables $\{Y_n\}_{n\geqslant0}$ converges in probability
if and only if for every subsequence of joint probability laws, $\{\mu_{n_k,m_k}\}_{k\geqslant0}$,
there exists a further subsequence which converges weakly to a probability measure $\mu$ such that
\be\label{2.17}\mu(\{(x,y)\in X\X X:x=y\})=1.\ee
\ep
\subsection{On invariant measures}

Let $X$ be an arbitrary Banach space with its Borel $\sig$-algebra $\cB(X)$.
We use $\cB_{\rm b}(X)$ to denote the space of bounded Borel measurable functions on $X$.
In the following we consider the $X$-valued stochastic process $u(t;x)$
with the initial datum $x\in X$.
For a set $\Gam\in\cB(X)$, we define the {\em transition functions}
$$\cP_t(x,\Gam)=\bP(u(t;x)\in\Gam)\hs\mb{for all }t\geqslant0.$$
The {\em Markovian transition semigroup} $\cP_t$ on $\cB_{\rm b}(X)$ is defined as
$$\cP_t\vp(x)=\bE\vp(u(t;x))=\int_{X}\vp(y)\cP_t(x,\di y),\hs t\geqslant0,\;
\vp\in\cB_{\rm b}(X),\;x\in X.$$
Note also that $\cP_t(x,\Gam)=\cP_t\chi_{\Gam}(x)$.
The dual semigroup $\cP^*_t$ of $\cP_t$ is defined on and into the set of
Borel probability measures $\nu$ on $X$ by
$$\cP^*_t\nu(\Gam)=\int_{X}\cP_t(x,\Gam)\nu(\di x),\hs\mb{for each }\Gam\in\cB(X).$$

The Markovian transition semigroup $\cP_t$ ($t\geqslant0$) is said to be {\em Feller}
if for arbitrary $\vp\in \cC_{\rm b}(X)$
and $t>0$, the mapping $x\mapsto\cP_t\vp(x)$ is continuous.
An {\em invariant measure} for the stochastic process $u(t;x)$ is a probability measure $\nu$ on $X$,
which is a fixed point for $\cP^*_t$, that is to say,
$$\int_X\cP_t(x,\Gam)\nu(\di x)=\nu(\Gam),\hs\mb{for each }\Gam\in\cB(X)\mb{ and }t\in\R^+.$$
Let $\nu$ be an invariant measure for $\cP_t$.
We say that $\nu$ is {\em ergodic} if
$$
\lim_{T\ra\8}\frac1T\int_0^T\cP_t\vp\di t=\int_X\vp(x)\nu(\di x)\hs\mb{for all }\vp\in L^2(X,\mu).
$$
\Vs

In what follows, we denote $c$ as an arbitrary positive constant,
which only depends on the parameters of
the original problem and the assumptions, i.e., $a$, $b$, $\T^2$, $N$
(only when the cut-off function $\de_N$ is involved below) and $\kappa$
(only when \eqref{2.3} is satisfied),
and may be different from line to line and even in the same line.
When the constant $c$ depends on some extra parameters, such as $m$, $p$, $k$ and $l$,
we denote it with subscripts $c_m$, $c_p$, $c_{m,p}$, $c_k$ or $c_l$ instead for emphasis.

\section{Estimates on the Galerkin Scheme for the Cut-off System}\label{s3}

In this section, we introduce the \textbf{Galerkin method} to open the discussions
in the following two sections.

We first apply the \textbf{cut-off function} to the original problem,
so as to make use of the global existence of solutions for finite-dimensional differential equations
with globally Lipschitz drift terms and diffusion coefficients when adopting the Galerkin scheme.
Specifically, we consider the cut-off system for \eqref{2.13},
\be\label{3.1}
\di u+[A^2u+\de_N(\|u\|_{2})f(u)]\di t=\de_N(\|u\|_{2})\phi(u)\di W,\hskip2mm t>0;\hs u(0)=u_0,
\ee
where $\de_N:\R^+\ra[0,1]$, $N\in\N^+$, is a smooth cut-off function such that
$$\de_N(r)=1\mb{ for }r\in[0,N]\hs\mb{and}\hs\de_N(r)=0\mb{ for }r\geqslant N+1.$$

For the Galerkin scheme we will use below,
we introduce the subspace $H_n$ of $H$, with dimension $n\geqslant0$, such that
$$H_n={\rm span}\{w_0,w_1,\cdots,w_n\}$$
and denote $P_n$ to be the projection operator from $H$ onto $H_n$.
Let
$$u^n(t)=\sum_{i=0}^na_i(t)w_i\hs\mb{and}\hs u^n_0=P_nu_0.$$
Note that all norms of $H_n$ are equivalent.
We consider the finite-dimensional stochastic differential equation
\be\label{3.2}
\di u^n+[A^2u^n+\de_N(\|u^n\|_{2})f_n(u^n)]\di t=\de_N(\|u^n\|_{2})P_n\phi(u^n)\di W,
\ee
with $u^n(0)=u^n_0$, where $f_n:H_n\ra H_n$ is defined as
$$f_n(u^n)=(a+3)u^n-4Au^n+bP_n|\na u^n|^2+P_n(u^n)^3.$$

There is no hard to check that the drift term and diffusion coefficient of \eqref{3.2}
are globally Lipschitz continuous in the finite-dimensional space $H_n$
by the involvement of the cut-off function $\de_N(\|u^n\|_{2})$ and \eqref{2.6}.
Referring to \cite[Section 5.2]{E13},
we know that there exists a unique global solution $u^n$ to \eqref{3.2} in $H_n$
for each initial datum $u_0$ satisfying \eqref{2.9}.

Following this we need the uniform estimates to deduce the tightness of the law of $u^n$.
We recall the Gagliardo-Nirenberg inequality (see \cite{SY02,WYD20}) for the discussion below.

\bl\label{le3.1}
Let $U$ be an open, bounded domain of the Lipschitz class in $\R^n$.
Assume that $\hat p\geqslant1$, $1\leqslant\hat q,\hat r\leqslant\8$, $0\leqslant k\leqslant m$,
$0<\theta\leqslant1$ and that
$$k-\frac{n}{\hat p}\leqslant\theta\(m-\frac{n}{\hat q}\)-(1-\theta)\frac{n}{\hat r}.$$
Then there is a positive constant $C$ such that
$$
\|u\|_{W^{k,\hat p}(U)}\leqslant C\|u\|_{W^{m,\hat q}(U)}^{\theta}\|u\|_{L^{\hat r}(U)}^{1-\theta},
\hs\mb{for all }u\in W^{m,\hat q}(U).
$$
\el

\bl\label{le3.2}
Let $\cS=(\W,\cF,\{\cF_t\}_{t\geqslant0},\bP,W)$ be a stochastic basis $p\geqslant2$
and $T\geqslant0$.
Suppose that $\phi$ satisfies \eqref{2.6}.
Given every $n\geqslant0$ and each $u_0$ satisfying \eqref{2.11},
there exists a positive constant $C_1$ independent of $n$,
such that the solution $u^n$ of \eqref{3.2} has the following uniform estimates,
\be
\bE\(\sup_{s\in[0,T]}\|A^mu^n(s)\|^p+\int_{0}^T\|A^mu^n\|^{p-2}
\|A^{m+1}u^n\|^2\di s\)\leqslant C_1,
\label{3.3}\ee
\be\label{3.4}
\bE\(\int_0^T\|A^{m+1}u^n(s)\|^2\di s\)^{\frac{p}2}\leqslant C_1,
\ee
\be
\bE\left\|A^{m-1}u^n(t)-\int_0^t\de_N(\|u^n\|_{2})P_nA^{m-1}
\phi(u^n(s))\di s\right\|^2_{W^{1,2}(0,T;H)}\leqslant C_1,
\label{3.5}\ee
and given $\al\in(0,1/2)$ additionally, it holds that
\be
\bE\left\|\int_0^t\de_N(\|u^n\|_{2})P_nA^m\phi(u^n)
\di W\right\|_{W^{\al,2}(0,T;H)}^p\leqslant C_1.
\label{3.6}\ee
\el

\bo We first show \eqref{3.3}.
Apply $A^m$ to \eqref{3.2} and It\^o's Formula to $\|A^mu^n\|^p$.
We obtain that for all $s\in[0,T]$,
\begin{align*}
&
  \frac1p\di\|A^mu^n\|^p+\|A^mu^n\|^{p-2}
  \|A^{m+1}u^n\|^2\di s\\
=&
  -\de_N(\|u^n\|_{2})\|A^mu^n\|^{p-2}
  \(A^mu^n,A^mf_n(u^n)\)\di s\\
&
  +\frac12\de_N^2(\|u^n\|_{2})\|A^mu^n\|^{p-4}
  \left[(p-2)(A^mu^n,P_nA^m\phi(u^n))_{L_2(\fU,\R)}^2\right.\\
&
  \left.+\|A^mu^n\|^2\|P_nA^m\phi(u^n)\|_{L_2(\fU,H)}^2\right]\di s\\
&
  +\de_N(\|u^n\|_{2})\|A^mu^n\|^{p-2}
  \(A^mu^n,P_nA^m\phi(u^n)\di W\)\\
:=&
 (I_1^p+I_2^p)\di s+I_3^p\di W.
\end{align*}
For arbitrary stopping times $\tau'$ and $\tau''$ with
$0\leqslant\tau'\leqslant s\leqslant\tau''\leqslant T$, we have
\begin{align}&
  \frac1p\|A^mu^n(s)\|^p+\int_{\tau'}^s\|A^mu^n\|^{p-2}
  \|A^{m+1}u^n\|^2\di s'\notag\\
\leqslant&
  \frac1p\|A^mu^n(\tau')\|^p+\int_{\tau'}^s\(|I_1^p|+|I_2^p|\)\di s'
  +\left|\int_{\tau'}^sI_3^p\di W\right|.
\label{3.7}\end{align}
Take the expectation of the supremum over $s\in[\tau',\tau'']$ for \eqref{3.7}, we have
\begin{align}&
  \bE\(\frac1p\sup_{s\in[\tau',\tau'']}\|A^mu^n(s)\|^p+\int_{\tau'}^{\tau''}
  \|A^mu^n\|^{p-2}\|A^{m+1}u^n\|^2\di s'\)\notag\\
\leqslant&
  \frac1p\bE\|A^mu^n(\tau')\|^p+\bE\int_{\tau'}^{\tau''}\(|I_1^p|+|I_2^p|\)\di s'
  +\bE\sup_{s\in[\tau',\tau'']}\left|\int_{\tau'}^sI_3^p\di W\right|.
\label{3.8}\end{align}

The key difficulty here lies in the estimation of $\(A^mu^n,A^mf_n(u^n)\)$
in $I_1^p$.
We are now devoted to addressing it.
First, we see that
\be\label{3.9}|(A^mu^n,A^mAu^n)|\leqslant\frac{2p-1}{6p}\|A^{m+1}u^n\|^2
+c_p\|A^mu^n\|^2.\ee
Observe by \eqref{2.1} that
\begin{align}
|\(A^mu^n,A^m|\na u^n|^2\)|
=&\left|\(A^{m+1}u^n,A^{m-1}|\na u^n|^2\)\right|\notag\\
\leqslant&\|A^{m+1}u^n\|\|A^{m-1}|\na u^n|^2\|\label{3.10}
\end{align}
for every $m\in\N^+$.
We can expand $A^{m-1}|\na u^n|^2$
to be the sum of finitely many (set to be $l_1$) summands (with each coefficient being $1$)
of the following form
$$
\frac{\pa^{m_1}u^n}{\pa x^{m_{11}}\pa y^{m_{12}}}
\frac{\pa^{m_2}u^n}{\pa x^{m_{21}}\pa y^{m_{22}}}\Hs\mb{with}\hs
m_{i1}+m_{i2}=m_i\geqslant1,\hs i=1,\,2,
$$
\be\label{3.11}2\leqslant\tilde{m}_j:=m_1+m_2\leqslant2m,\hs
j=1,\cdots,l_1.
\ee
It is easy to see that $m_i\leqslant2m-1$ for $i=1,2$.
By selecting $\theta=\frac{2m_i+1}{4(m+1)}$ in Lemma \ref{le3.1}, we can deduce that
$$
\left|\frac{\pa^{m_i} u^n}{\pa x^{m_{i1}}\pa y^{m_{i2}}}\right|_4
\leqslant\|u^n\|_{m_i,4}
\leqslant c_m\|u^n\|_{2(m+1)}^{\frac{2m_i+1}{4(m+1)}}
\|u^n\|^{\frac{4m-2m_i+3}{4(m+1)}}
$$
and hence with \eqref{3.11} in mind,
\begin{align*}
\left\|\frac{\pa^{m_1}u^n}{\pa x^{m_{11}}\pa y^{m_{12}}}
\frac{\pa^{m_2}u^n}{\pa x^{m_{21}}\pa y^{m_{22}}}\right\|
\leqslant&
  \|u^n\|_{m_1,4}\|u^n\|_{m_2,4}\\
\leqslant&c_m\|u^n\|_{2(m+1)}^{\frac{\tilde{m}_j+1}{2(m+1)}}
  \|u^n\|^{\frac{4m-\tilde{m}_j+3}{2(m+1)}}.
\end{align*}
Now we get back to \eqref{3.10} and immediately infer that
\begin{align}
&\|A^{m+1}u^n\|\|A^{m-1}|\na u^n|^2\|\notag\\
\leqslant&
  c_m\sum_{l=1}^{l_1}\|A^{m+1}u^n\|^{\frac{2m+\tilde{m}_j+3}{2(m+1)}}
  \|u^n\|^{\frac{4m-\tilde{m}_j+3}{2(m+1)}}\notag\\
\leqslant&
  \sum_{l=1}^{l_1}\(\frac{2p-1}{6l_1p}\|A^{m+1}u^n\|^2
  +c_{m,p}\|u^n\|^{\frac{2(4m-\tilde{m}_j+3)}{2m-\tilde{m}_j+1}}\)\notag\\
\leqslant&
  \frac{2p-1}{6p}\|A^{m+1}u^n\|^2+c_{m,p}\(\|u^n\|^{2(2m+3)}+1\),\label{3.12}
\end{align}
where we have used $\ve$-Young inequality and noticed that
$\frac{4m-\tilde{m}_j+3}{2m-\tilde{m}_j+1}\leqslant2m+3$.
Next we treat the second nonlinear term $P_n(u^n)^3$ of $f_n(u^n)$.
Similarly, we know
\begin{align}|\(A^mu^n,A^m(u^n)^3\)|
=&\left|\(A^{m+1}u^n,A^{m-1}(u^n)^3\)\right|\notag\\
\leqslant&\|A^{m+1}u^n\|\|A^{m-1}(u^n)^3\|.\label{3.13}\end{align}
And then $A^{m-1}(u^n)^3$ can be expanded as the sum of $l_2$ summands
(with each coefficient being $1$) of the form
$$
\frac{\pa^{m_1}u^n}{\pa x^{m_{11}}\pa y^{m_{12}}}
\frac{\pa^{m_2}u^n}{\pa x^{m_{21}}\pa y^{m_{22}}}
\frac{\pa^{m_3}u^n}{\pa x^{m_{31}}\pa y^{m_{32}}}\Hs\mb{with}\hs
m_{i1}+m_{i2}=m_i
$$
and $\tilde{m}_k:=m_1+m_2+m_3\leqslant2m-2$, $i=1,\,2,\,3$, $k=1,\,\cdots,l_2$.
By choosing $\theta=\frac{3m_i+2}{6(m+1)}$ in Lemma \ref{le3.1}, one deduces that
$$
\left|\frac{\pa^{m_i}u^n}{\pa x^{m_{i1}}\pa y^{m_{i2}}}\right|_6
\leqslant\|u^n\|_{m_i,6}
\leqslant c_m\|u^n\|_{2(m+1)}^{\frac{3m_i+2}{6(m+1)}}
\|u^n\|^{\frac{6m-3m_i+4}{6(m+1)}}
$$
and
$$
\left\|\frac{\pa^{m_1}u^n}{\pa x^{m_{11}}\pa y^{m_{12}}}
\frac{\pa^{m_2}u^n}{\pa x^{m_{21}}\pa y^{m_{22}}}
\frac{\pa^{m_3}u^n}{\pa x^{m_{31}}\pa y^{m_{32}}}\right\|
\leqslant c_m\|u^n\|_{2(m+1)}^{\frac{\tilde{m}_k+2}{2(m+1)}}
\|u^n\|^{\frac{6m-\tilde{m}_k+4}{2(m+1)}}.
$$
And hence
\begin{align}
&\|A^{m+1}u^n\|\|A^{m-1}(u^n)^3\|\notag\\
\leqslant&
  \sum_{l=1}^{l_2}\|A^{m+1}u^n\|^{\frac{2m+\tilde{m}_k+4}{2(m+1)}}
  \|u^n\|^{\frac{6m-\tilde{m}_k+4}{2(m+1)}}\notag\\
\leqslant&
  \sum_{l=1}^{l_2}\(\frac{2p-1}{6l_2p}\|A^{m+1}u^n\|^2
  +c_{m,p}\|u^n\|^{\frac{2(6m-\tilde m_k+4)}{2m-\tilde{m}_k}}\)\notag\\
\leqslant&
  \frac{2p-1}{6p}\|A^{m+1}u^n\|^2+c_{m,p}\(\|u^n\|^{2(2m+3)}+1\),
\label{3.14}\end{align}
where the relation $\frac{6m-\tilde m_k+4}{2m-\tilde{m}_k}\leqslant2m+3$ is used.
Combining the inequalities from \eqref{3.9} to \eqref{3.14},
we obtain
\begin{align}|I_1^p|\leqslant&
  \frac{2p-1}{2p}\|A^mu^n\|^{p-2}\|A^{m+1}u^n\|^2\notag\\
&  +c_{m,p}\de_N(\|u^n\|_{2})\(\|A^mu^n\|^{p}
  +\|A^mu^n\|^{p-2}\(\|u^n\|^{2(2m+3)}+1\)\)\notag\\
\leqslant&
  \frac{2p-1}{2p}\|A^mu^n\|^{p-2}\|A^{m+1}u^n\|^2\notag\\
& +c_{m,p}\de_N(\|u^n\|_{2})\(\|A^mu^n\|^{p}+\|u^n\|^{(2m+3)p}+1\).
\label{3.15}\end{align}

For $I_2^p$, by \eqref{2.6}, we have
\begin{align}
\|A^m\phi(u^n)\|_{L_2(\fU,H)}\leqslant&
  \kappa(|u^n|_\8)\(1+\|A^mu^n\|\)\hs\mb{and hence}\notag\\
\label{3.16}
\|(A^mu^n,P_n A^m\phi(u^n))\|^2_{L_2(\fU,\R)}\leqslant&
  2\kappa^2(|u^n|_\8)\|A^mu^n\|^2\(1+\|A^mu^n\|^2\),\\
\notag
\|P_nA^m\phi(u^n)\|_{L_2(\fU,H)}^2\leqslant&
  2\kappa^2(|u^n|_\8)\(1+\|A^mu^n\|^2\).
\end{align}
Therefore, we obtain
\begin{align}|I_2^p|\leqslant&
  c\de_N(\|u^n\|_{2})\|A^mu^n\|^{p-2}\kappa^2(|u^n|_\8)
  (1+\|A^mu^n\|^2)\notag\\
\leqslant&
  c_p\de_N(\|u^n\|_{2})\kappa^2(|u^n|_\8)
  (1+\|A^mu^n\|^p).
\label{3.17}\end{align}
For the stochastic term, using H\"older's inequality,
the BDG inequality and  \eqref{3.16},
we have
\begin{align*}&\bE\sup_{s\in[\tau',\tau'']}\left|\int_{\tau'}^sI_3^p\di W\right|\\
\leqslant&
  c_p\bE\(\int_{\tau'}^{\tau''}\de_N^2(\|u^n\|_{2})\|A^mu^n\|^{2(p-1)}
  \kappa^2(|u^n|_\8)(1+\|A^mu^n\|^2)\di s'\)^{\frac12}\\
\leqslant&
  c_p\bE\(\sup_{s\in[\tau',{\tau''}]}\|A^mu^n(s)\|^p\)^{\frac12}\\
&  \(\int_{\tau'}^{\tau''}\de_N(\|u^n\|_{2})\kappa^2(|u^n|_\8)\|A^mu^n\|^{p-2}
  (1+\|A^mu^n\|^2)\di s'\)^{\frac12},
\end{align*}
and hence
\begin{align}\bE\sup_{s\in[\tau',\tau'']}\left|\int_{\tau'}^sI_3^p\di W\right|\leqslant&
  c_p\bE\int_{\tau'}^{\tau''}\de_N(\|u^n\|_{2})\kappa^2(|u^n|_\8)
  \(1+\|A^mu^n\|^p\)\di s'\notag\\
&  +\frac1{2p}\bE\sup_{s\in[\tau',{\tau''}]}\|A^mu^n(s)\|^p.\label{3.18}
\end{align}
It follows from \eqref{3.8}, \eqref{3.15}, \eqref{3.17}, \eqref{3.18} that
\begin{align}
&
 \bE\(\sup_{s\in[\tau',\tau'']}\|A^mu^n(s)\|^p+\int_{\tau'}^{\tau''}\|A^mu^n\|^{p-2}
 \|A^{m+1}u^n\|^2\di s'\) \label{3.19}\\
\leqslant&
c_{m,p}\bE\int_{\tau'}^{\tau''}\de_N(\|u^n\|_{2})
 \left[\(1+\kappa^2(|u^n|_\8)\)\(1+\|A^mu^n\|^p\)+\|u^n\|^{(2m+3)p}\right]
 \di s'\notag\\
&
 +2\bE\|A^mu^n(\tau')\|^p.\notag
\end{align}

In order that the expectation in \eqref{3.19} makes sense, we define
$$\ol\rho_R:=\inf\{t\geqslant0:\|A^mu^n(t)\|> R\}.$$
Since $H^2\subset L^\8\subset H$, then $\kappa$ and $\|u^n\|$ is bounded
over the set $\{u^n:\de_N(\|u^n\|_2)>0\}$,
which excludes the set $\{u^n:|u^n|_\8>\epsilon(N+1)\}$.
Here we use $\epsilon$ to denote the embedding constant that $|u|_\8\leqslant \epsilon\|u\|_2$ for each $u\in H^2$.
We thus note that the upper bounds of the terms
$$\de_N(\|u^n\|_{2})\(1+\kappa^2(|u^n|_\8)\)\mb{ and }
\de_N(\|u^n\|_{2})\kappa^2(|u^n|_{\8})\|u^n\|^{(2m+3)p}$$
only depend upon $N$ (via $\de_N$ and $\kappa$), but independent of $R$ and $n$.
Now we can apply the stochastic Gronwall's inequality (see \cite{GHZ09,WYD20}) to
\eqref{3.19} and obtain
\be\label{3.20}
\bE\(\sup_{t\in[0,T\wedge\ol\rho_R]}\|A^mu^n(t)\|^p+\int_{0}^{T\wedge\ol\rho_R}
 \|A^mu^n\|^{p-2}\|A^{m+1}u^n\|^2\di s\)\leqslant C'_1,
\ee
where $C'_1=C'_1(a,b,p,T,N,u_0,\kappa,\T^2)$ is a positive constant independent of $R$ and $n$.
Now let $R\ra\8$. One sees $T\wedge\ol\rho_R\ra T$ in \eqref{3.20} and hence \eqref{3.3}
by the Dominated Convergence Theorem.
\Vs

Next, we show \eqref{3.4}.
By replacing $\tau'$ and $s$ by $0$ and $s$ respectively and setting $p=2$ in \eqref{3.7},
we can obtain
$$\int_0^s\|A^{m+1}u^n\|^2\di s'\leqslant
  \frac12\|A^mu^n_0\|^2+\int_0^s\(|I_1^2|+|I_2^2|\)\di s'
  +\left|\int_{0}^sI_3^2\di W\right|.
$$
By \eqref{3.15} with a slight adjustment to the parameters and \eqref{3.17}, we can also get
\begin{align}&\int_0^s\|A^{m+1}u^n\|^2\di s'\notag\\
\leqslant&
  c_m\int_0^s\de_N(\|u^n\|_{2})
  \left[\(1+\kappa^2(|u^n|_\8)\)\(1+\|A^mu^n\|^2\)+\|u^n\|^{2(2m+3)}\right]\di s'\notag\\
&
 +\|A^mu^n_0\|^2+\left|\int_{0}^sI_3^2\di W\right|.\label{3.21}
\end{align}
Similar to \eqref{3.18} and using BDG inequality, we know
\begin{align}&\bE\sup_{s\in[0,T]}\left|\int_{0}^sI_3^2\di W\right|^{\frac{p}{2}}\notag\\
\leqslant&c_p\bE\(\sup_{s\in[0,T]}\|A^mu^n(s)\|^2\)^{\frac{p}2}\notag\\
&  +c_p\bE\(\int_0^T\de_N(\|u^n\|_{2})\kappa^2(|u^n|_\8)
  \(1+\|A^mu^n\|^2\)\di s'\)^{\frac{p}2}\notag\\
\leqslant&
  c_p\bE\sup_{s\in[0,T]}\|A^mu^n(s)\|^p\notag\\
& +c_p\bE\int_0^T\de_N(\|u^n\|_{2})
  \kappa^p(|u^n|_\8)\(1+\|A^mu^n\|^p\)\di s'.
\label{3.22}\end{align}
Then considering the $p/2$ power of \eqref{3.21} and taking the expected value of the supremum,
we have by \eqref{3.22}
\begin{align}&\bE\(\int_0^s\|A^{m+1}u^n\|^2\di s'\)^{\frac{p}2}\notag\\
\leqslant&
  c_{m,p}\bE\|A^mu^n_0\|^p+c_p\bE\sup_{s\in[0,T]}\|A^mu^n(s)\|^p\notag\\
&  +c_{m,p}\bE\int_0^T\de_N(\|u^n\|_{2})\|u^n\|^{(2m+3)p}\di s'\notag\\
&
  +c_{m,p}\bE\int_0^T\de_N(\|u^n\|_{2})\(1+\kappa^p(|u^n|_\8)\)\(1+\|A^mu^n\|^p\)\di s'.
\label{3.23}\end{align}
The inequality \eqref{3.4} follows immediately from \eqref{3.3}, \eqref{3.23} and the properties of
$\de_N$ and $\kappa$.
\Vs

Now we show \eqref{3.5}. Applying  \eqref{3.2}, we have the following equation
with stochastic integral in $H$,
\begin{align}
&A^{m-1}u^n(t)-\int_0^t\de_N(\|u^n\|_{2})P_nA^{m-1}\phi(u^n)\di W\notag\\
=&
  A^{m-1}u^n_0-\int_0^tA^{m+1}u^n\di s
  -\int_0^t\de_N(\|u^n\|_{2})P_nA^{m-1}f_n(u^n)\di s
  \notag\\
:=&
  J_1+J_2+J_3.
\label{3.24}\end{align}
Note that
\be\bE\|J_1\|^2\leqslant\bE\|A^{m-1}u_0\|^2.
\ee
By \eqref{3.3} with $p=2$, we obtain
\begin{align}
\bE\|J_2\|_{W^{1,2}(0,T;H)}^2=&\bE\int_0^T(\|J_2(s)\|^2+\|A^{m+1}u^n\|^2)\di s\notag\\
\leqslant&(T^2+1)\bE\int_0^{T}\|A^{m+1}u^n\|^2\di s\leqslant C'_2.
\end{align}
According to the discussion from \eqref{3.10} to \eqref{3.14}, we can similarly have
\begin{align}
  &\bE\|J_3\|^2_{W^{1,2}(0,T;H)}\notag\\
\leqslant&
  (T^2+1)\bE\int_0^T\de_N(\|u^n\|_{2})\|A^{m-1}f_n(u^n(t))\|^2\di s\notag\\
\leqslant&
  c(T^2+1)\bE\int_0^T\de_N(\|u^n\|_{2})\(\|A^{m+1}u^n\|^{2}+\|u^n\|^{2}
  +\|u^n\|^{\frac{2(2m+3)}{m+2}}+1\)\di s\notag\\
\leqslant& C'_3.
\label{3.27}\end{align}
Then \eqref{3.5} follows immediately from \eqref{3.24} to \eqref{3.27}.
\Vs

To prove \eqref{3.6}, since $\al\in(0,1/2)$, we adopt \eqref{2.8}, \eqref{2.6}
and \eqref{3.3} and obtain
\begin{align*}&
  \bE\left\|\int_0^t\de_N(\|u^n\|_{2})P_nA^m\phi(u^n(s))
  \di W\right\|_{W^{\al,2}(0,T;H)}^p\\
\leqslant&
  c\bE\int_0^T\de_N(\|u^n\|_{2})\kappa^p(|u^n|_\8)(1+\|A^mu^n\|^p)\di s
  \leqslant C'_4,
\end{align*}
where we used the properties of $\de_N$ and $\kappa$ that have been applied above.
This is \eqref{3.6}.
The proof is complete.
\eo

\section{Local and Global Existence of Martingale Solutions}\label{s4}

\subsection{Local existence of martingale solutions}\label{ss4.1}

Let $\mu_0$ be a given initial distribution on $H^{2m}$.
We fix a stochastic basis $\cS=(\W,\cF,\{\cF_t\}_{t\geqslant0},\bP,W)$
upon which is defined an $\cF_0$ measurable random element $u_0$ with distribution $\mu_0$.
Consider the sequence of Galerkin approximations $\{u^n\}_{n\geqslant1}$ solving \eqref{3.2}
relative to this basis and initial condition.
We consider the phase spaces
\be\label{4.1}
\cX_u=L^2(0,T;H^{2m})\cap \cC([0,T];H^{2(m-2)}),\hs
\cX_W=\cC([0,T];\fU)
\ee
and $\cX=\cX_u\X\cX_W$, with $\cX_u$ the space
where the solution $u^n$ lives and $\cX_W$ the set where
the driving Brownian motions are defined.
We consider the probability measures
\be\label{4.2}\mu_u^n(\cdot)=\bP(u^n\in\cdot)\in{\rm Pr}(\cX_u)\hs\mb{and}\ee
\be\label{4.2A}\mu_W^n(\cdot)=\mu_W(\cdot)=\bP(W\in\cdot)\in{\rm Pr}(\cX_W).\ee
This defines a sequence of probability measures $\mu^n:=\mu^n_u\X\mu_W^n$ on the phase space $\cX$.
We now show that this sequence is tight as follows by using Lemma \ref{le3.2}.
\bl\label{le4.1}Suppose that $\phi$ satisfies \eqref{2.6} and $\mu_0$ satisfies \eqref{2.9}.
Then the sequence $\{\mu^n\}_{n\geqslant1}$ is tight and weakly compact over the phase space $\cX$.
\el

\bo
Applying Lemma \ref{le2.1} (i) by setting $X_2=H^{2(m-1)}$, $X_0=H^{2m}$, $X_1=H^{2(m+1)}$,
$p=2$ and $\al\in(0,1/2)$, we have that
$$L^2(0,T;H^{2(m+1)})\cap W^{\al,2}(0,T;H^{2(m-1)})\subset\subset L^2(0,T;H^{2m}).$$
Given $R>1$, define a set
\begin{align*}
B_R^1=\{&u\in L^2(0,T;H^{2(m+1)})\cap W^{\al,2}(0,T;H^{2(m-1)}):\\
&\|u^n\|^2_{L^2(0,T;H^{2(m+1)})}
+\|u^n\|^2_{W^{\al,2}(0,T;H^{2(m-1)})}\leqslant R^2\}.
\end{align*}
Then $B_R^1$ is compact in $L^2(0,T;H^{2m})$.
The Chebyshev inequality applies and we have
\begin{align}&\mu_u^n((B_R^1)^C)\notag\\=&
  \bP\(\|u^n\|^2_{L^2(0,T;H^{2(m+1)})}+\|u^n\|^2_{W^{\al,2}(0,T;H^{2(m-1)})}\geqslant R^2\)\notag\\
\leqslant&
  \bP\(\|u^n\|^2_{L^2(0,T;H^{2(m+1)})}\geqslant\frac{R^2}{2}\)
  +\bP\(\|u^n\|^2_{W^{\al,2}(0,T;H^{2(m-1)})}\geqslant\frac{ R^2}2\)\notag\\
\leqslant&
  \frac{2}{R^2}\bE\(\int_0^T\|A^{m+1}u^n(s)\|^2\di s
  +\|u^n\|^2_{W^{\al,2}(0,T;H^{2(m-1)})}\)\leqslant\frac{C'_5}{R^2},
\label{4.3}\end{align}
independent of $n$, where we have used \eqref{3.3} with $p=2$, \eqref{3.6},
\eqref{3.5} and \eqref{2.2}.

Take $\al\in(\frac1q,\frac12)$ with $\al q>1$.
By Lemma \ref{le2.1} (ii) with $Y_0=H^{2(m-2)}$ and $Y=H^{2(m-1)}$, we have the compact embeddings
$$
W^{1,2}(0,T;H^{2(m-1)})\subset\subset\cC([0,T];H^{2(m-2)}),$$$$
W^{\al,q}(0,T;H^{2(m-1)})\subset\subset \cC([0,T];H^{2(m-2)}).
$$
Given $R>1$, let $B_R^{21}$ and $B_R^{22}$ be the closed balls of radius $R$ in the spaces
$W^{1,2}(0,T;H^{2(m-1)})$ and $W^{\al,q}(0,T;H^{2(m-1)})$ respectively.
It follows that for $R>1$, $B^2_R:=B_R^{21}+B_R^{22}$ is compact in $\cC([0,T];H^{2(m-2)})$.
Due to the inclusion
$$
\{u^n\in B_R^2\}\supset\left\{u^n(t)-\int_0^tP_n\phi(u^n)\di W\in B_R^{21}\right\}\cap
\left\{\int_0^tP_n\phi(u^n)\di W\in B_R^{22}\right\},
$$
we obtain by Chebyshev inequality and \eqref{3.6}, \eqref{3.5} that
\begin{align}\mu_u^n((B_R^2)^C)\leqslant&
  \bP\(\left\|u^n(t)-\int_0^tP_n\phi(u^n)\di W\right\|_{W^{1,2}(0,T;H^{2(m-1)})}\geqslant R^2\)\notag\\
&
  +\bP\(\left\|\int_0^tP_n\phi(u^n)\di W\right\|_{W^{\al,q}(0,T;H^{2(m-1)})}\geqslant R^q\)
  \leqslant\frac{C'_6}{R^2},
\label{4.4}\end{align}
also independent of $n$.

It is trivial that $B_R^1\cap B_R^2$ is compact in $L^2(0,T;H^{2m})\cap\cC([0,T];H^{2(m-2)})$
for every $R>0$.
It follows from \eqref{4.3} and \eqref{4.4} that
$$\mu^n_u((B_R^1\cap B_R^2)^C)\leqslant\mu_u^n((B_R^1)^C)+\mu^n_u((B_R^2)^C)\leqslant\frac{C'_5+C'_6}{R^2}.$$
Then for every $\ve>0$, we are allowed to pick a set $\cA_\ve$ (by increasing $R$ for $B_R^1\cap B_R^2$)
such that
\be\mu_u^n(\cA_\ve)\geqslant 1-\frac\ve2,\Hs\mb{for all }n\geqslant0.\label{4.5}\ee

Now we consider the sequence $\{\mu_W^n\}_{n\geqslant0}$, which actually identically equals to $\mu_W$,
and is hence weakly compact.
By Prokhorov's Theorem (\cite{DaPZ14}), $\{\mu_W^n\}_{n\geqslant0}$ is surely tight.
This helps find a compact set $\cB_\ve$ in $\cC([0,T];\fU_0)$ such that for all $n\geqslant0$,
\be\label{4.6}\mu_W^n(\cB_\ve)\geqslant1-\frac\ve2.\ee
Combining \eqref{4.5} and \eqref{4.6}, we know that for every $\ve>0$,
the compact set $\cA_\ve\X\cB_\ve$ in $\cX$ satisfies that for all $n\geqslant0$,
$$\mu^n(\cA_\ve\X\cB_\ve)\geqslant1-\ve,$$
and therefore $\{\mu^n\}_{n\geqslant0}$ is tight in $\cX$ and finally weakly compact.
The proof is finished.
\eo

Now we can deduce the following theorem to guarantee the existence of martingale solutions
of \eqref{3.1}.

\bt\label{th4.2} Suppose that $\phi$ satisfies \eqref{2.6} and
$\mu_0$ is a probability measure on $H^{2m}$ satisfying \eqref{2.9}.
Then there exists a subsequence $n_k$ and a probability space $(\tilde\W,\tilde\cF,\tilde\bP)$,
on which there lies a sequence of $\cX$-valued random variables $(\tilde u^{n_k},\tilde{W}^{n_k})$
such that
\benu\item[(i)] $(\tilde{u}^{n_k},\tilde{W}^{n_k})$ converges almost surely, in the topology of $\cX$,
    to an element $(\tilde u,\tilde W)\in\cX$,
    the law of $(\tilde u^{n_k},\tilde W^{n_k})$ is $\mu^{n_k}$ for each $k$ and
    $\mu^{n_k}$ weakly converges to the law $\mu$ of $(\tilde u,\tilde W)$.
\item[(ii)] $\tilde{W}^{n_k}$ is a cylindrical Wiener process,
    relative to the filtration $\tilde\cF_t^{n_k}$,
    given by the completion of
    $$\sig((\tilde u^{n_k}(s),\tilde W^{n_k}(s));s\leqslant t).$$
\item[(iii)] every pair $(\tilde u^{n_k},\tilde{W}^{n_k})$ satisfies \eqref{3.2} with only $u^n$ and $W$
    replaced by $\tilde{u}^{n_k}$ and $\tilde W^{n_k}$ therein.
\eenu
Let $\tilde{\cS}=(\tilde{\W},\tilde{\cF},\{\tilde\cF_t\}_{t\geqslant0},\tilde{\bP},\tilde{W})$,
with $\tilde{\cF}_t$ the completion of $\sig((\tilde{u}(s),\tilde{W}(s));s\leqslant t)$.
Then $(\tilde{\cS},\tilde{u})$ is a global martingale solution of \eqref{3.1}.
As a result, by defining a stopping time
\be\label{4.7}\rho_N:=\inf\{t\geqslant0:\|\tilde{u}(t)\|_{2}> N\},\hs\mb{for }N>0,\ee
the triple $(\tilde{\cS},\tilde{u},\rho_N)$ is a local martingale solution of \eqref{2.13}.
\et
\br Here the martingale solution of \eqref{3.1} is a slight modification of that of \eqref{2.13}, i.e.,
Definition \ref{de2.1}.
Similarly, the pathwise solution of \eqref{3.1} adopted in the next section is
also a slight modification of Definition \ref{de2.2}.
\er
\noindent\textit{Proof of Theorem \ref{th4.2}.}\hs
The conclusion (i) is easily obtained by the Skorohod embedding theorem
(see \cite[Theorem 2.4]{DaPZ14}).
The conclusions (ii) and (iii) can be ensured by a similar procedure to that in
\cite[Section 4.3.4]{B95}.

Next we show that $(\tilde{\cS},\tilde{u})$ is a global martingale solution of \eqref{3.1}
through three steps.
\vs

\textbf{Step 1.} \textit{Improvement of regularity on the phase space}

Observe from the conclusion (i) that
\begin{align}\tilde{u}^{n_k}\ra\tilde u&
  \hs\mb{in }\cX_u\mb{ for }\bP\mb{ a.s. }\w\in\tilde\W,
  \label{4.8}\\
\tilde{W}^{n_k}\ra\tilde{W}&
  \hs\mb{in }\cX_W\mb{ for }\bP\mb{ a.s. }\w\in\tilde\W.\label{4.9}
\end{align}
According to the conclusion (iii) above, we know $\tilde{u}^{n_k}$ has the same estimates as $u^{n_k}$
stated in Lemma \ref{le3.2}.
Using Banach-Alaoglu theorem, \eqref{3.3} and \eqref{3.6} for $\tilde{u}^{n_k}$, we can obtain
$\ol{u}\in L^2(\tilde\W;L^2(0,T;H^{2(m+1)}))$ and $\ol{\ol u}\in L^p(\tilde\W;L^\8(0,T;H^{2m}))$
such that
\begin{align}\label{4.10}\tilde{u}^{n_k}\ra&
  \ol{u}\hs\mb{in }L^2(\tilde\W;L^2(0,T;H^{2(m+1)}))\;\mb{weakly},\hs\mb{and}\\
\label{4.11}\tilde{u}^{n_k}\ra&
  \ol{\ol u}\hs\mb{in }L^p(\tilde\W;L^\8(0,T;H^{2m}))\;\mb{weakly star,}
\end{align}
where the sequence $\{\tilde{u}^{n_k}\}$ is perhaps chosen to be a subsequence.
On account of \eqref{2.9} and with Lemma \ref{le3.2} applied, we have
$$
\sup_{k\in\N}\bE\sup_{t\in[0,T]}\|A^{m-2}\tilde{u}^{n_k}\|^p
\leqslant c\sup_{k\in\N}\bE\sup_{t\in[0,T]}\|A^{m}\tilde{u}^{n_k}\|^p<\8.
$$
Then use the Vitali convergence theorem (see \cite{F48}) to the convergence in the conclusion (i)
and base the estimation of \eqref{3.3}, we have
\be\label{4.12}\tilde{u}^{n_k}\ra\tilde{u}\hs\mb{in }L^p(\tilde\W;L^{\8}(0,T;H^{2(m-2)})).\ee
Now given each measurable subset $\cR$ of $[0,T]\X\tilde{\W}$ and $v\in H^{2(m+1)}$, by \eqref{4.10},
\eqref{4.11} and \eqref{4.12}, we have
$$\bE\int_0^T\chi_{\cR}\langle\tilde{u},v\rangle_{2m}\di s
=\bE\int_0^T\chi_{\cR}\langle\ol{\ol u},v\rangle_{2m}\di s
=\bE\int_0^T\chi_{\cR}\langle\ol{u},v\rangle_{2m}\di s,
$$
where $\langle\cdot,\cdot\rangle_{2m}$ means the dual product between $H^{-2m}$ and $H^{2m}$
and we have used the dense embeddings.
This indicates that $\tilde{u}=\ol{u}=\ol{\ol u}$ and therefore,
\be\label{4.13}\tilde{u}\in L^p(\tilde{\W};L^\8(0,T;H^{2m}))\cap L^2(\tilde{\W};L^2(0,T;H^{2(m+1)})).\ee
\vs

\textbf{Step 2.} \textit{Convergence to the integral form of \eqref{3.1}}

In this step, we show that for the limit $(\tilde{u},\tilde{W})$, the following integral equation
\be\label{4.14}\tilde{u}(t)+\int_0^{t}[A^2\tilde{u}+\de_N(\|\tilde{u}\|_{2})f(\tilde{u})]\di s=
u_0+\int_0^{t}\de_N(\|\tilde{u}\|_{2})\phi(\tilde{u})\di W\ee
for almost every $(t,\w)\in[0,T]\X\tilde{\W}$, holds in $H^{2(m-1)}$,
for which, it suffices to show that \eqref{4.14} holds in $H$ by \eqref{4.13}, \eqref{2.6} and
\eqref{2.9}.

We can follow an argument similar to that in \cite[Section 7]{DGHT11} for this proof.
But the nonlinear forcing and stochastic terms need to be discussed specifically in our situation.
For sake of completeness and the reader's convenience, we present the whole procedure in details.

Since $m\geqslant1$, applying embeddings and Vitali convergence theorem again to the convergence
in \eqref{4.8} and the fact that
$$\sup_{k\in\N}\bE\(\int_0^T\|A^{m+1}\tilde{u}^{n_k}(s)\|^2\di s\)^{\frac{p}2}<\8
\hs\mb{by \eqref{3.6}},$$
one concludes that $\tilde{u}^{n_k}\ra\tilde{u}$ in $L^p(\tilde{\W};L^2(0,T;H^{2(m+1)}))$.
Then by picking a subsequence, we can further assume that
\be\label{4.15}U_k:=\tilde{u}^{n_k}-\tilde{u}\ra0\hs\mb{in }H^{2(m+1)}
\mb{ for almost surely }(t,\w)\in[0,T]\X\tilde{\W}.\ee

Fix $v\in H^{4}$.
We can infer from \eqref{4.13} and embeddings that
$$
\left|\int_0^t\langle A^2U_k,v\rangle_{2}\di s\right|
\leqslant
c\|v\|_{2}\int_0^T\|U_k\|_{2}\di s$$
\be\label{4.16}\mb{and hence}\hs\int_0^t\langle A^2\tilde{u}^{n_k},v\rangle_2\di s\ra
\int_0^t\langle A^2\tilde{u},v\rangle_{2}\di s,
\ee
for almost every $(t,\w)\in[0,T]\X\tilde\W$, where $\langle\cdot,\cdot\rangle_2$
is the dual product between $H^{-2}$ and $H^2$.

Then we address the nonlinear term $\de_N(\|\tilde{u}\|_{2})f(\tilde{u})$.
Actually,
\begin{align}&
  \left|\int_0^t\langle\de_N(\|\tilde{u}^{n_k}\|_{2})f_{n_k}(\tilde{u}^{n_k})
  -\de_N(\|\tilde{u}\|_{2})f(\tilde{u}),v\rangle_2\di s\right|\notag\\
\leqslant&
  \left|\int_0^t\de_N(\|\tilde{u}^{n_k}\|_{2})\langle f_{n_k}(\tilde{u}^{n_k})-
  f(\tilde{u}),v\rangle_2\di s\right|\notag\\
&  +\left|\int_0^t(\de_N(\|\tilde{u}^{n_k}\|_{2})-\de_N(\|\tilde{u}\|_{2}))
  \langle f(\tilde{u}),v\rangle_2\di s\right|\notag\\
\leqslant&
  \left|\int_0^t\de_N(\|\tilde{u}^{n_k}\|_{2})\langle f(\tilde{u}^{n_k})-
  f(\tilde{u}),P_{n_k}v\rangle_2\di s\right|\notag\\
&  +\left|\int_0^t\de_N(\|\tilde{u}^{n_k}\|_{2})\langle f(\tilde{u}),(1-P_{n_k})v\rangle_2
  \di s\right|\notag\\
&
  +\left|\int_0^t(\de_N(\|\tilde{u}^{n_k}\|_{2})-\de_N(\|\tilde{u}\|_{2}))
  \langle f(\tilde{u}),v\rangle_2\di s\right|:=
  K_1^{k}+K_2^k+K_3^k.\label{4.17}
\end{align}
For the estimation of $K_1^k$, first we see that
\begin{align}|&\langle f(\tilde u^{n_k})-f(\tilde{u}),P_{n_k}v\rangle_2|\notag\\
\leqslant&
  |(a+3)\langle U_k,P_{n_k}v\rangle_2|+4|\langle AU_k,P_{n_k}v\rangle_2|
  +|b\langle|\na\tilde{u}^{n_k}|^2-|\na\tilde{u}|^2,P_{n_k}v\rangle_2|\notag\\
&
  +|\langle(\tilde{u}^{n_k})^3-(\tilde{u})^3,P_{n_k}v\rangle_2|\notag\\
\leqslant&
  c\|U_k\|_2\|v\|_2+c|\langle|\na U_k|^2,P_{n_k}v\rangle_2|
  +c|\langle\na U_k\cdot\na\tilde{u}^{n_k},P_{n_k}v\rangle_2|\notag\\
&
  +c|\langle U_k^3,P_{n_k}v\rangle_2|+c|\langle\tilde{u}^{n_k}U_k^2,P_{n_k}v\rangle_2|
  +c|\langle(\tilde{u}^{n_k})^2U_k,P_{n_k}v\rangle_2|.
\label{4.18}\end{align}
Note that all the components in dual products in \eqref{4.18} are actually in $H$ by embeddings.
Hence by using Lemma \ref{le3.1}, it follows with no hard from \eqref{4.18} and embeddings that
$$|\langle f(\tilde u^{n_k})-f(\tilde{u}),P_{n_k}v\rangle_2|\leqslant
  c\|U_k\|_2\|v\|_2\(1+\|\tilde{u}^{n_k}\|^2_2+\|U_k\|^2_2\),
$$
and thus by the property of $\de_N$ and \eqref{4.15},
\be K_1^k\leqslant
  c\|v\|_2\int_0^t\|U_k\|_2\(1+\|U_k\|^2_2\)\di s\ra0,\hs\mb{as }k\ra\8.
\label{4.19}\ee
For $K_2^k$, by embeddings and \eqref{4.13}, we have as $k\ra\8$,
\begin{align}
&K_2^k\leqslant c\|(1-P_{n_k})v\|_2\int_0^t\|f(\tilde{u})\|\di s\notag\\
\leqslant&
c\|(1-P_{n_k})v\|_2\int_0^t(\|\tilde{u}\|_2+\|\tilde{u}\|_2^2+\|\tilde{u}\|_2^3)\di s
\ra0.\label{4.20}
\end{align}
As to $K_3^k$, we need to notice that $H^{2(m+1)}\subset H^2$, which with \eqref{4.15}
implies that $\tilde{u}^{n_k}-\tilde{u}\ra0$ in $H^2$ for a.s. $(t,\w)\in[0,T]\X\tilde{\W}$
and hence
\be\label{4.21}\tilde{u}^{n_k}\ra\tilde{u}\hs\mb{in }H^2\mb{ for }t\in[0,T]
\mb{ in measure and a.s. }\w\in\tilde{\W}.\ee
Note also by \eqref{4.15} that for a.s. $\w\in\tilde\W$,
\be\label{4.22}|\langle f(\tilde{u}(t)),v\rangle_2\leqslant c\|v\|_2|
\sup_{t\in[0,T]}(\|\tilde{u}(t)\|_2+\|\tilde{u}(t)\|_2^2+\|\tilde{u}(t)\|_2^3)<\8.
\ee
The facts \eqref{4.21} and \eqref{4.22} guarantee that for each $\ve>0$,
there exists a subset $E_\ve\subset[0,T]$ with its Lebesgue measure so small that
\be\label{4.23}\int_{E_\ve}|\de_N(\|\tilde{u}^{n_k}\|_{2})-\de_N(\|\tilde{u}\|_{2})|
  |\langle f(\tilde{u}),v\rangle_2|\di s<\frac{\ve}{2},
\ee
and for a.s. $\w\in\tilde{\W}$,
\be\label{4.24}\tilde{u}^{n_k}(t)\ra\tilde{u}(t),\hs\mb{in }H^2
\mb{ uniformly in }t\in[0,T]\sm E_{\ve}\mb{ as }k\ra\8.\ee
Moreover, $\de_N$ is obviously uniformly continuous on $\R^+$.
Finally, we can infer from \eqref{4.23} and \eqref{4.24} that for all $\ve>0$,
\be K_3^k\leqslant
  \(\int_{E_\ve}+\int_{[0,T]\sm E_\ve}\)|\de_N(\|\tilde{u}^{n_k}\|_{2})
  -\de_N(\|\tilde{u}\|_{2})||\langle f(\tilde{u}),v\rangle_2|\di s<\ve.\label{4.25}
\ee
when $k$ is sufficiently large.
Now the conclusions \eqref{4.17}, \eqref{4.19}, \eqref{4.20} and \eqref{4.25} imply
that as $k\ra\8$,
\be\label{4.26}
\int_0^t\langle\de_N(\|\tilde{u}^{n_k}\|_{2})f_{n_k}(\tilde{u}^{n_k}),v\rangle_2\di s
\ra\int_0^t\langle\de_N(\|\tilde{u}\|_{2})f(\tilde{u}),v\rangle_2\di s.
\ee
For the initial datum $u_0$, it is easy to see that as $k\ra\8$
\be\label{4.27}\tilde{u}^{n_k}(0)=u_0^{n_k}\ra u_0\hs\mb{in }H^{2m}.\ee

Next we cope with the convergence for the stochastic term.
First we have
\begin{align*}&\|\de_N(\|\tilde{u}^{n_k}\|_{2})P_{n_k}\phi(\tilde{u}^{n_k})
  -\de_N(\|\tilde{u}\|_{2})\phi(\tilde{u})\|_{L_2(\fU,H)}\\
\leqslant&\|\de_N(\|\tilde{u}^{n_k}\|_{2})\phi(\tilde{u}^{n_k})
-\de_N(\|\tilde{u}\|_{2})\phi(\tilde{u}))\|_{L_2(\fU,H)}\\
&  +\|\de_N(\|\tilde{u}\|_{2})(1-P_{n_k})\phi(\tilde{u})\|_{L_2(\fU,H)}.
\end{align*}
By the continuity of $\de_N$ and $\phi$ (by \eqref{2.6}), the embedding $H^2\subset H$
and the inequality,
$$\|\de_N(\|\tilde{u}\|_{2})\phi(\tilde{u})\|_{L_2(\fU,H)}\leqslant
\de_N(\|\tilde{u}\|_{2})\kappa(|\tilde{u}|_{\8})(1+\|\tilde{u}\|),$$
it yields from \eqref{4.15} that
$$\|\de_N(\|\tilde{u}^{n_k}\|_{2})P_{n_k}\phi(\tilde{u}^{n_k})
  -\de_N(\|\tilde{u}\|_{2})\phi(\tilde{u})\|_{L_2(\fU,H)}\ra0,\hs\mb{as }k\ra\8,$$
for almost every $(t,\w)\in[0,T]\X\tilde{\W}$.
Applying \eqref{2.6} again, we observe from \eqref{3.6} that
\begin{align*}
&\sup_{k\in\N}\bE\int_0^T
\|\de_N(\|\tilde{u}^{n_k}\|_{2})P_{n_k}\phi(\tilde{u}^{n_k})\|_{L_2(\fU,H)}^2\di s\\
\leqslant&c\bE\int_0^T\de_N(\|\tilde{u}^{n_k}\|_{2})\kappa(|\tilde{u}^{n_k}|_{\8})
(1+\|\tilde{u}^{n_k}\|^2)\di s\leqslant C.
\end{align*}
By the Dominated Convergence Theorem, we know
$$\de_N(\|\tilde{u}^{n_k}\|_{2})P_{n_k}\phi(\tilde{u}^{n_k})
  \ra\de_N(\|\tilde{u}\|_{2})\phi(\tilde{u})\hs\mb{in }L^2(\tilde\W;L^2(0,T;L_2(\fU,H))).
$$
This also deduces the following convergence
$$\de_N(\|\tilde{u}^{n_k}\|_{2})P_{n_k}\phi(\tilde{u}^{n_k})
  \ra\de_N(\|\tilde{u}\|_{2})\phi(\tilde{u})
$$
in probability in $L^2(0,T;L_2(\fU,H))$. Then using the conclusion (i) and Lemma \ref{le2.2},
we obtain
$$
\int_0^t\de_N(\|\tilde{u}^{n_k}\|_{2})P_{n_k}\phi(\tilde{u}^{n_k})\di\tilde{W}^{n_k}
\ra\int_0^t\de_N(\|\tilde{u}\|_{2})\phi(\tilde{u})\di\tilde{W}
$$
in probability in $L^2(0,T;H)$, and hence we can assume
\be\label{4.28}
\int_0^t\langle\de_N(\|\tilde{u}^{n_k}\|_{2})P_{n_k}\phi(\tilde{u}^{n_k})\di\tilde{W}^{n_k},
v\rangle_2
\ra\int_0^t\langle\de_N(\|\tilde{u}\|_{2})\phi(\tilde{u})\di\tilde{W},v\rangle_2
\ee
for almost every $(t,\w)\in[0,T]\X\tilde{\W}$, by picking a subsequence if necessary.
\vs

As a final result, by \eqref{4.15}, \eqref{4.16}, \eqref{4.26}, \eqref{4.27}, \eqref{4.28}
and the conclusion (iii) that $(\tilde{u}^{n_k},\tilde{W}^{n_k})$ satisfies \eqref{3.2}
with only $u^n$ and $W$ replaced by $\tilde{u}^{n_k}$ and $\tilde W^{n_k}$ therein,
we can infer that for all $v\in H^2$ and almost every $(t,\w)\in[0,T]\X\tilde\W$,
\begin{align}&\left\langle\tilde{u}(t)+\int_0^{t}[A^2\tilde{u}+\de_N(\|\tilde{u}\|_{2})
  f(\tilde{u})]\di s,v\right\rangle_2\notag\\
=&\lim_{k\ra\8}\langle\tilde{u}^{n_k}(t),v\rangle_2
  +\lim_{k\ra\8}\left\langle\int_0^{t}A^2\tilde{u}^{n_k}\di s,v\right\rangle_2\notag\\
&
  +\lim_{k\ra\8}\left\langle\int_0^{t}\de_N(\|\tilde{u}^{n_k}\|_{2})f(\tilde{u}^{n_k})\di s,
  v\right\rangle_2\notag\\
=&
  \lim_{k\ra\8}\left\langle\tilde{u}^{n_k}(t)+\int_0^{t}A^2\tilde{u}^{n_k}\di s
  +\int_0^{t}\de_N(\|\tilde{u}^{n_k}\|_{2})f(\tilde{u}^{n_k})\di s,v\right\rangle_2\notag\\
=&
  \lim_{k\ra\8}\int_0^t\langle\de_N(\|\tilde{u}^{n_k}\|_{2})
  P_{n_k}\phi(\tilde{u}^{n_k})\di\tilde{W}^{n_k},v\rangle_2\notag\\
=&
 \int_0^t\langle\de_N(\|\tilde{u}\|_{2})\phi(\tilde{u})\di\tilde{W},v\rangle_2,\label{4.29}
\end{align}
via a subsequence of $n_k$ (still denoted by $n_k$).
The density of $H^2$ in $H$ ensures that \eqref{4.29} also holds for all $v\in H$ and
with $\langle\cdot,\cdot\rangle_2$ replaced by $(\cdot,\cdot)$.
\Vs

\textbf{Step 3.} \textit{Improvement of regularity in time}

Now to finalize the existence of martingale solution, the last task is the proof of
the continuity of $\tilde{u}$ in $t$.
First we introduce an extra random variable $z$ to treat the stochastic term, such that
\be\label{4.30}\di z+A^2z\di t=\de_N(\|\tilde{u}\|_{2})\phi(\tilde{u})\di\tilde{W},
\hs z(0)=u_0.\ee
Since $\de_N(\|\tilde{u}\|_{2})\phi(\tilde{u})\in L^2(\tilde{\W};L^2(0,T;L_2(\fU,H^{2m})))$
by \eqref{2.6} and \eqref{4.13}, we have (by \cite[Chapter 8]{DaPZ14})
\be\label{4.31}
z\in L^2(\tilde{\W};\cC([0,T];H^{2m}))\cap L^2(\tilde{\W};L^2(0,T;H^{2(m+1)})).
\ee
Then define $U=\tilde{u}-z$.
By \eqref{3.1} and \eqref{4.30}, we have
\be\label{4.32}\frac{\di U}{\di t}=-A^2U-\de_N(\|U+z\|_{2})f(U+z),\hs U(0)=0.\ee
According to \eqref{4.31} and \eqref{4.13}, we obtain that
$U\in L^2(\tilde{\W},L^2(0,T;H^{2(m+1)}))$.
Each term on the right hand of \eqref{4.32} belongs to
$L^2(\tilde{\W};L^2(0,T;H^{2(m-1)}))$.
Actually, we have proved that
$$A^mU\in L^2(\tilde{\W};L^2(0,T,H^2))\hs\mb{and}\hs
\frac{\di}{\di t}A^mU\in L^2(\tilde{\W};L^2(0,T;H^{-2}))$$
Utilizing \cite[Chapter 3, Lemma 1.2]{Tem01}, we know for almost every $\w\in\tilde{\W}$,
$U\in\cC([0,T];H^{2m})$.
Then combining this with \eqref{4.13} and \eqref{4.31}, we infer that
$\tilde{u}\in L^2(\tilde{\W};\cC([0,T];H^{2m}))$.
\Vs

Now we have shown that $(\tilde{\cF},\tilde{u})$ is a global martingale solution of \eqref{3.1}.
By the definition \eqref{4.7}, we know that $\rho_N$ is a stopping time.
One sees that
$$\int_0^{t\wedge\rho_N}\de_N(\|\tilde{u}\|_{2})f(\tilde{u})\di s
=\int_0^{t\wedge\rho_N}f(\tilde{u})\di s$$
$$\mb{and}\hs
\int_0^{t\wedge\rho_N}\de_N(\|\tilde{u}\|_{2})\phi(\tilde{u})\di\tilde{W}
=\int_0^{t\wedge\rho_N}\phi(\tilde{u})\di\tilde{W}.
$$
Taking \eqref{4.29} into consideration,
we infer that $(\tilde{\cS},\tilde{u},\rho_N)$ is a local martingale solution of \eqref{2.13}.
The proof is finally finished.
\qed

\subsection{Global existence of martingale solutions}

With a further assumption on $\phi$, $b$ and $u_0$, we actually can obtain the existence of
global martingale solutions for \eqref{2.13}.
For the proof to this consequence, after an estimate on the Galerkin scheme for \eqref{2.13},
we will repeat the procedure given above in Subsection \ref{ss4.1} to obtain the existence of
local martingale solutions first.

To discuss the existence of solutions of this equation, we still need to introduce
the cut-off function.
Based on the settings in Section \ref{s3}, we consider the related equations,
$$
\di u+[A^2u+\de_N(\|u\|_2)f(u)]\di t=\phi(u)\di W
\hs\mb{and}\hs u(0)=u_0,\hs\mb{and}
$$
\be\label{4.33}
\di u^n+[A^2u^n+\de_N(\|u^n\|_2)f_n(u^n)]\di t=P_n\phi(u^n)\di W
\hs\mb{and}\hs u^n(0)=u^n_0.
\ee
The problem \eqref{4.33} has a unique global solution $u^n$ in $H_n$
for every $u_0$ satisfying \eqref{2.9}.
Now we give renewed estimates for the solution of \eqref{4.33}.

\bl\label{le4.4} Let $\cS=(\W,\cF,\{\cF_t\}_{t\geqslant0},\bP,W)$ be a stochastic basis
and $T\geqslant0$.
Suppose that $\phi$ satisfy \eqref{2.3}.
Given every $n\geqslant0$ and each initial datum $u_0$ satisfying \eqref{2.12},
there exists a positive constant $C_2$
independent of $n$ and $N$,
such that the solution $u^n$ of \eqref{3.2} has the following uniform estimates,
\be
\bE\(\sup_{s\in[0,T]}\|A^mu^n(s)\|^p+\int_{0}^T\|A^mu^n\|^{p-2}
\|A^{m+1}u^n\|^2\di s\)\leqslant C_2,
\label{4.34}\ee
\be\label{4.35}
\bE\(\int_0^T\|A^{m+1}u^n(s)\|^2\di s\)^{\frac{p}2}\leqslant C_2,
\ee
\be
\bE\left\|A^{m-1}u^n(t)-\int_0^t\de_N(\|u^n\|_{2})P_nA^{m-1}
\phi(u^n(s))\di s\right\|^2_{W^{1,2}(0,T;H)}\leqslant C_2,
\label{4.36}\ee
and given $\al\in(0,1/2)$ additionally, it holds that
\be
\bE\left\|\int_0^t\de_N(\|u^n\|_{2})P_nA^m\phi(u^n)
\di W\right\|_{W^{\al,2}(0,T;H)}^p\leqslant C_2.
\label{4.37}\ee
\el
\bo For the estimate discussed later, we first estimate
$\disp\bE\sup_{s\in[0,T]}\|u^n(s)\|^{q_0}$ with $q_0:=(2m+3)p$.
Applying It\^o's Formula to $\|u^n\|^{q_0}$, for all stopping times $\tau',\tau''$ with
$0\leqslant\tau'\leqslant s\leqslant\tau''\leqslant T$, we can deduce that
\begin{align}
&
  \frac1{q_0}\|u^n(s)\|^{q_0}+\int_{\tau'}^{s}\|u^n\|^{q_0-2}\(\|Au^n\|^2
  +\de_N(\|u^n\|_2)|u^n|_4^4\)\di s'\notag\\
=&
  \frac1{q_0}\|u^n(\tau')\|^{q_0}\!\!-\!\!\int_{\tau'}^{s}\de_N(\|u^n\|_{2})\|u^n\|^{q_0-2}
  \(u^n,(a+3)u^n-4Au^n+b|\na u^n|^2\)\di s'\notag\\
&
  +\!\frac12\!\!\int_{\tau'}^{s}\!\!\|u^n\|^{q_0-4}\!
  \left[(q_0-2)(u^n,P_n\phi(u^n))_{L_2(\fU,\R)}^2\!+\!\|u^n\|^2\|P_n\phi(u^n)\|_{L_2(\fU,H)}^2\right]
  \di s'  \notag\\
&
  +\int_{\tau'}^{s}\|u^n\|^{q_0-2}
  \(u^n,P_n\phi(u^n)\di W\).\label{4.38}
\end{align}
Observe by Lemma \ref{le3.1}, integration by parts, H\"older inequality
and Young's inequality that
\begin{align}&\left|\(u^n,(a+3)u^n+4Au^n+b|\na u^n|^2\)\right|\notag\\
\leqslant&
  c\|u^n\|^2+4|(u^n,Au^n)|+\frac{|b|}{2}\left|\((u^n)^2,\De u^n\)\right|\notag\\
\leqslant&
  c\|u^n\|^2+4|(u^n,Au^n)|+\frac{|b|}{4}\(|u^n|_4^4+(\De u^n,\De u^n)\)\notag\\
=&
  \frac{|b|}{4}\(|u^n|_4^4+\|Au^n\|^2\)+c(\|u^n\|^2+\|u^n\|\|Au^n\|)\notag\\
\leqslant&
  \frac{|b|}{4}\(|u^n|_4^4+\|Au^n\|^2\)+c\(\|u^n\|^2+\|Au^n\|^{\frac32}
  +|u^n|_4^{\frac52}\).\label{4.39}
\end{align}
Combining \eqref{4.38}, \eqref{4.39} and \eqref{3.17} with $\de_N(\|u^n\|_2)$,
$\kappa(|u^n|_\8)$ and $m$ therein replaced by $1$, $\kappa$ and $0$, correspondingly,
we have by the property of polynomials and $|b|<4$ that
\begin{align}&
  \frac1{q_0}\|u^n(s)\|^{q_0}
  +\frac12\(1-\frac{|b|}{4}\)\int_{\tau'}^s\|u^n\|^{q_0-2}
  \(\|Au^n\|^2+\de_N(\|u^n\|_2)|u^n|_4^4\)\di s'
  \notag\\
\leqslant&
  \frac1{q_0}\|u^n(\tau')\|^{q_0}+\int_{\tau'}^s\|u^n\|^{q_0-2}
  \left[c\|Au^n\|^{\frac32}-\frac12\(1-\frac{|b|}{4}\)
  \|Au^n\|^2\right]\di s'\notag\\
&
  +\int_{\tau'}^s\de_N(\|u^n\|_2)\|u^n\|^{q_0-2}\left[c\(|u^n|_4^2+|u^n|_4^{\frac52}\)
  -\frac12\(1-\frac{|b|}{4}\)|u^n|_4^4\right]\di s'\notag\\
&
  +c_{m,p}\int_{\tau'}^s\(1+\|u^n\|^{q_0}\)\di s'
  +\int_{\tau'}^{s}\|u^n\|^{q_0-2}\(u^n,P_n\phi(u^n)\di W\)\notag\\
\leqslant&
  \frac1{q_0}\|u^n(\tau')\|^{q_0}+c_{m,p}\int_{\tau'}^s\(1+\|u^n\|^{q_0}\)\di s'\notag\\
&
  +c(s-\tau')+\int_{\tau'}^{s}\|u^n\|^{q_0-2}\(u^n,P_n\phi(u^n)\di W\).
  \label{4.40}
\end{align}
Taking the expectation of the supremum over $s\in[\tau',\tau'']$ for \eqref{4.40}
and by \eqref{3.18} with $\de_N(\|u^n\|_2)$, $\kappa(|u^n|_\8)$ and $m$ therein
replaced by $1$, $\kappa$ and $0$, correspondingly, we obtain
\begin{align}&
  \bE\sup_{s\in[\tau',\tau'']}\|u^n(s)\|^{q_0}
  +\(1-\frac{|b|}{4}\)q_0\bE\int_{\tau'}^{\tau''}\|u^n\|^{q_0-2}\|Au^n\|^2\di s
  \notag\\
\leqslant&
  2\bE\|u^n(\tau')\|^{q_0}+c_{m,p}\bE\int_{\tau'}^{\tau''}\(1+\|u^n\|^{q_0}\)\di s
  +c(\tau''-\tau').\label{4.41}
\end{align}
Note here that \eqref{4.41} has no business with $N$ except $u^n$ itself.
Now similar to the argument for \eqref{3.3}, we can infer by using stopping times and
the stochastic Gronwall's inequality that
\be\label{4.42}\bE\sup_{s\in[0,T]}\|u^n(s)\|^{q_0}
  +\bE\int_0^T\|u^n\|^{q_0-2}\|Au^n\|^2\di s\leqslant C'_7,\ee
where $C'_7=C'_7(a,b,m,p,T,\kappa,u_0,\T^2)$.
\vs

Next we estimate $\disp\bE\sup_{s\in[0,T]}\|A^mu^n\|^p$.
Following the similar estimation in the proof of Lemma \ref{le3.2}
and replacing $\de_N(\|u^n\|_2)$ and $\kappa(|u^n|_\8)$ related to the stochastic term
by $1$ and $\kappa$, correspondingly,
we have an estimate similar to \eqref{3.19},
\begin{align}
&
 \bE\(\sup_{s\in[\tau',\tau'']}\|A^mu^n(s)\|^p
 +\int_{\tau'}^{\tau''}\|A^mu^n\|^{p-2}
 \|A^{m+1}u^n\|^2\di s'\)\notag\\
\leqslant&
 2\bE\|A^mu^n(\tau')\|^p+c_{m,p}\bE\int_{\tau'}^{\tau''}
 \(1+\|A^mu^n\|^p\)\di s'\notag\\
&
 +c_{m,p}\bE\int_{\tau'}^{\tau''}\|u^n\|^{(2m+3)p}\di s'.\label{4.43}
\end{align}
Also \eqref{4.43} is independent of all $n$ and $N$.
With a similar argument, \eqref{4.42} and the fact that $u_0$ satisfies \eqref{2.12},
we easily conclude \eqref{4.34}.
The inequalities \eqref{4.35}, \eqref{4.36} and \eqref{4.37} can be similarly obtained,
We omit the detailed argument.\eo

Next we need the local uniqueness of the solution in the following sense.
\bl\label{le4.5} Suppose the conditions of Lemma \ref{le4.4} holds.
Let $N_1,\,N_2\in\N^+$ with $N_2>N_1$ and $u_i^n$ be the unique solution of \eqref{4.33}
with $N$ replaced by $N_i$, $i=1,\,2$.
Define a stopping time $\zeta$ so that
$$\zeta=\inf\{t\geqslant0:\|u_1^n(t)\|_2\vee\|u_2^n(t)\|_2>N_1\}.$$
Then
\be\label{4.44}
\bP\(u_1^n(t)=u_2^n(t);\mb{ for all }t\in[0,\zeta]\)=1.\ee
\el
\bo Let $\ol v=u_2^n-u_1^n$. We know $\ol v(0)=0$ and that $\ol v$ satisfies the following equation
\begin{align}&\di\ol v
+\left[A^2\ol v+\(\de_{N_2}(\|u_2^n\|_2)f_n(u_2^n)-\de_{N_1}(\|u_1^n\|_2)f_n(u_1^n)\)\right]\di t\notag\\
=&P_n(\phi(u_2^n)-\phi(u_1^n))\di W.\label{4.45}
\end{align}
Define a stopping time for each $R>0$,
$$\varrho_R=\inf\{t\geqslant0:\|u_1^n(t)\|_{2m}+\|u^n_2(t)\|_{2m}>R\}\wedge\zeta.$$
Apply $A^m$ to \eqref{4.45} and use It\^o's Formula to $\|A^m\ol v\|^p$.
It yields for all stopping times $\tau'$, $\tau''$ with
$0\leqslant\tau'\leqslant s\leqslant\tau''\leqslant\varrho_R$ that
\begin{align}&
  \frac1p\|A^m\ol v(s)\|^p+\int_{\tau'}^{s}\|A^m\ol v\|^{p-2}\|A^{m+1}\ol v\|^2\di s'\notag\\
\leqslant&
  \frac1p\|A^m\ol v(\tau')\|^p+\frac{p-1}2\int_{\tau'}^{s}\|A^m\ol v\|^{p-2}
  \|A^m(\phi(u^n_2)-\phi(u^n_1))\|^2_{L_2(\fU,H)}\di s'\notag\\
&
  -\int_{\tau'}^{s}\|A^m\ol v\|^{p-2}
  \(A^m\ol v,\de_{N_2}(\|u^n_2\|_2)A^mf(u^n_2)
  -\de_{N_1}(\|u^n_1\|_2)A^mf(u^n_1)\)\di s'\notag\\
&
  +\left|\int_{\tau'}^{s}\|A^m\ol v\|^{p-2}
  \(A^m\ol v,P_nA^m\(\phi(u^n_2)-\phi(u^n_1)\)\di W\)\right|\notag\\
:=&
  \frac1p\|A^m\ol v(\tau')\|^p+M_1+M_2+M_3.\label{4.46}
\end{align}
Note that when $s\leqslant\varrho_R$, $\de_{N_i}(u_i^n(s))=1$, $i=1,\,2$.
Hence by
\begin{align}|M_1|=&\int_{\tau'}^{s}\|A^m\ol v\|^{p-2}|\(A^m\ol v,A^m(f(u^n_2)
  -f(u^n_1))\)|\di s'\notag\\
\leqslant&
  c\int_{\tau'}^{s}\|A^m\ol v\|^{p}\di s'
  +\int_{\tau'}^s\|A^m\ol v\|^{p-2}(A^m\ol v,A^{m+1}\ol v)\di s'\notag\\
&
  +\int_{\tau'}^s\|A^m\ol v\|^{p-2}\left|\(A^{m+1}\ol v,
  A^{m-1}\(b(|\na u^n_2|^2-|\na u^n_1|^2)\)\)\right|\di s'\notag\\
&
  +\int_{\tau'}^s\|A^m\ol v\|^{p-2}\left|\(A^{m+1}\ol v,
  A^{m-1}\((u^n_2)^3-(u^n_1)^3\)\)\right|\di s'\notag\\
\leqslant&
  \frac{2p-1}{2p}\int_{\tau'}^s\|A^m\ol v\|^{p-2}\|A^{m+1}\ol v\|^2\di s'
  +c\int_{\tau'}^{s}\|A^m\ol v\|^{p}\di s'\notag\\
&
  +c_p\int_{\tau'}^s\|A^m\ol v\|^{p-2}\left\|A^{m-1}
  \(\na\ol v\cdot\na(u^n_2+u^n_1)\)\right\|^2\di s'\notag\\
&
  +c_p\int_{\tau'}^s\|A^m\ol v\|^{p-2}\left\|A^{m-1}
  \(\ol v\((u^n_2)^2+u_2^nu_1^n+(u_1^n)^2\)\)\right\|^2\di s'.\label{4.47}
\end{align}
Since by the discussions to \eqref{3.10} and \eqref{3.13} and the embeddings,
\be\left\|A^{m-1}\(\na\ol v\cdot\na(u^n_2+u^n_1)\)\right\|^2\leqslant
  c_m\|A^m\ol v\|^2\(\|u^n_1\|_{2m}^2+\|u^n_2\|_{2m}^2\),\label{4.48}\ee
\begin{align}&\left\|A^{m-1}\(\ol v\((u^n_2)^2+u_2^nu_1^n+(u_1^n)^2\)\)\right\|^2\notag\\
\leqslant&
  c_m\|A^m\ol v\|^2\(\|u^n_1\|_{2m}^4+\|u^n_2\|_{2m}^4\),\label{4.49}
\end{align}
we have
\begin{align}|M_1|\leqslant&\frac{2p-1}{2p}\int_{\tau'}^s\|A^m\ol v\|^{p-2}
\|A^{m+1}\ol v\|^2\di s'\notag\\
&+c_{m,p}\int_{\tau'}^s\|A^m\ol v\|^p\(\|u^n_1\|_{2m}^4+\|u^n_2\|_{2m}^4+1\)\di s'.
\label{4.50}\end{align}
By \eqref{2.3}, we see
\be\label{4.51}|M_2|\leqslant c_p\int_{\tau'}^s\|A^m\ol v\|^p\di s'.\ee
For $M_3$, we use BDG inequality and obtain that
\begin{align}\bE\sup_{s\in[\tau',\tau'']}|M_3|\leqslant&
  c\bE\(\int_{\tau'}^{\tau''}\|A^m\ol v\|^{2p}\di s\)^{\frac12}\notag\\
\leqslant&
  c\bE\(\sup_{s\in[\tau',\tau'']}\|A^m\ol v(s)\|^p\)^{\frac12}
  \(\int_{\tau'}^{\tau''}\|A^m\ol v\|^p\di s\)^{\frac12}\notag\\
\leqslant&
  \frac{1}{2p}\bE\sup_{s\in[\tau',\tau'']}\|A^m\ol v(s)\|^p
  +c_p\bE\int_{\tau'}^{\tau''}\|A^m\ol v\|^p\di s.\label{4.52}
\end{align}
Consider the expectation of the supremum of \eqref{4.46} over $s\in[\tau',\tau'']$ and combine
\eqref{4.50}, \eqref{4.51} and \eqref{4.52}, we get
\begin{align}&\bE\sup_{s\in[\tau',\tau'']}\|A^m\ol v(s)\|^p
  +\bE\int_{\tau'}^{\tau''}\|A^m\ol v\|^{p-2}\|A^{m+1}\ol v\|^2\di s\notag\\
\leqslant&
  2\bE\|A^m\ol v(\tau')\|^p+c_{m,p}\bE\int_{\tau'}^{\tau''}\|A^m\ol v\|^p
  \(1+\|u^n_1\|_{2m}^4+\|u^n_2\|_{2m}^4\)\di s'.\label{4.53}\end{align}
Then apply the stochastic Gronwall's inequality to \eqref{4.53} and we conclude that
\be\label{4.54}\bE\(\sup_{s\in[0,\varrho_R]}\|A^m\ol v(s)\|^p
+\int_{0}^{\varrho_R}\|A^m\ol v\|^{p-2}\|A^{m+1}\ol v\|^2\di s\)
=C'_8\bE\|A^m\ol v(0)\|^p=0,\ee
where $C'_8$ only depends on $N_1$, $N_2$ and $m$.
Observe that $\varrho_R\ra\zeta$ as $R\ra+\8$.
We see that \eqref{4.54} indicates \eqref{4.44}.
The proof is complete.
\eo

Now we are ready to show the global existence of martingale solutions for \eqref{2.13}.
Let $n,\,N\in\N^+$, and $U^{n,N}$ be the solution of \eqref{4.33}.
We repeat the argument presented in Subsection \ref{ss4.1}
and follow the settings given there.
Let
\be\hat\mu_u^{n,N}(\cdot)=\bP(U^{n,N}\in\cdot)\in{\rm Pr}(\cX_u),
\hs\hat\mu_W^{n,N}(\cdot)=\mu_W(\cdot)\label{4.55}\ee
and $\hat\mu^{n,N}=\hat\mu_u^{n,N}\X\hat\mu_W^{n,N}$.
For each fixed $N\in\N^+$, we have the same conclusion as Lemma \ref{le4.1}, saying, for each $T>0$,
as long as $\mu_0$ satisfies \eqref{2.9},
the sequence $\{\hat\mu^{n,N}\}_{n\geqslant1}$ is tight in the phase space $\cX$.
As a final result, we infer the following consequence.

\bt\label{th4.6} Suppose that $\mu_0$ is a probability measure on $H^{2m}$ satisfying \eqref{2.10},
$\phi$ satisfies \eqref{2.3} and $|b|<4$.
Then there exists a global martingale solution to \eqref{2.13}.\et

\bo Following the same procedure to that given in the proof of Theorem \ref{th4.2},
we know that for each $T>0$ and $N\in\N^+$, there exists a subsequence $n_k$
and a probability space $(\hat\W,\hat\cF,\hat\bP)$, on which there is a sequence of $\cX$-valued
random variables $(\hat U^{n_k,N},\hat W^{n_k,N})$ such that
\benu\item[(i)] $(\hat U^{n_k,N},\hat{W}^{n_k,N})$ converges almost surely, in the topology of $\cX$,
    to an element $(\hat U^N,\hat W^N)\in\cX$,
    the law of $(\hat U^{n_k,N},\hat W^{n_k,N})$ is $\hat\mu^{n_k,N}$ for each $k$ and
    $\hat\mu^{n_k,N}$ weakly converges to the law $\hat\mu^N$ of $(\hat U^N,\hat W^N)$.
\item[(ii)] $\hat{W}^{n_k,N}$ is a cylindrical Wiener process,
    relative to the filtration $\hat\cF_t^{n_k,N}$, given by the completion of
    $$\sig((\hat U^{n_k,N}(s),\hat W^{n_k,N}(s));s\leqslant t).$$
\item[(iii)] every pair $(\hat U^{n_k,N},\tilde{W}^{n_k,N})$ satisfies \eqref{4.33}
    with only $u^n$ and $W$ replaced by $\hat U^{n_k,N}$ and $\hat W^{n_k,N}$ therein.
\eenu
For each solution $\hat U^{n,N}$ and $R\in\R^+$, define the stopping time $\zeta_R^{n,N}$ such that
\be\label{4.56}\zeta_R^{n,N}=\inf\{t\geqslant0:\|\hat U^{n,N}(t)\|_2>R\}.\ee
Note here in \eqref{4.56} $n$ denotes the dimension of $H_n$ and $N$ is
implied in the cut-off function $\de_N$.
\vs

Next we adopt the diagonal method to extend the local martingale solution to a global one.
Firstly, when $N=1$, we have a subsequence $n_k^1$ of $n$ such that
$$(\hat U^{n_k^1,1},\hat W^{n_k^1,1})\mb{ converges a.s. to }(\hat U^1,\hat W^1),$$
$$\mb{and }\hat\mu^{n_k^1,1}\mb{ weakly converges to }\hat\mu^1.$$
Then by Theorem \ref{th4.2}, we know $(\hat\cS^1,\hat U^1,\hat\rho_1)$ is a local martingale
solution of \eqref{2.13} with
$$\hat\rho_1=\inf\{t\geqslant0:\|\hat U^1\|_2>1\},$$
$\hat\cS^1=(\hat\W,\hat\cF^1,\{\hat\cF_t^1\}_{t\geqslant0},\hat\bP,\hat W)$ and
$\hat\cF_t^1$ the completion of $\sig((\hat U^1(s),\hat W^1(s));s\leqslant t)$.
For $N=2$, we further have a subsequence $n_k^2$ of $n_k^1$ such that
$$(\hat U^{n_k^2,2},\hat W^{n_k^2,2})\mb{ converges a.s. to }(\hat U^2,\hat W^2),$$
$$\mb{and }\hat\mu^{n_k^2,2}\mb{ weakly converges to }\hat\mu^2.$$
The triple $(\hat\cS^2,\hat U^2,\hat\rho_2)$ is thus also a local martingale of \eqref{2.13}
with $\hat\rho_2=\inf\{t\geqslant0:\|\hat U^2\|_2>2\}$, $\hat\cS^2=(\hat\W,\hat\cF^2,
\{\hat\cF_t^2\}_{t\geqslant0},\hat\bP,\hat W)$ and
$\hat\cF_t^2$ being the completion of $\sig((\hat U^2(s),\hat W^2(s));s\leqslant t)$.
In this sense, for each $N\in\N^+$, we have a subsequence $n_k^{N+1}$ of $n_k^N$ such that
$$(\hat U^{n_k^{N+1},N+1},\hat W^{n_k^{N+1},N+1})
\mb{ converges a.s. to }(\hat U^{N+1},\hat W^{N+1}),$$
$$\mb{and }
\hat\mu^{n_k^{N+1},N+1}\mb{ weakly converges to }\hat\mu^{N+1}.$$
And we find for each $N\in\N^+$, $(\hat\cS^N,\hat U^N,\hat\rho_N)$ is a local martingale
of \eqref{2.13} with $\hat\rho_N=\inf\{t\geqslant0:\|\hat U^N\|_2>N\}$,
$\hat\cS^N=(\hat\W,\hat\cF^N,\{\hat\cF_t^N\}_{t\geqslant0},\hat\bP,\hat W)$ and
$\hat\cF_t^N$ being the completion of $\sig((\hat U^N(s),\hat W^N(s));s\leqslant t)$.

Now we pick the sequence $\{(\hat U^{n_k^k,k},\hat W^{n_k^k,k})\}_{k\in\N^+}$.
Given each $N\in\N^+$, let
\be\label{4.57}\zeta_N=\liminf_{k\ra\8}\zeta^{n_k^k}_N\hs\mb{and}\hs
\zeta_N^{n_k^k}=\bigwedge_{j=N}^{\8}\zeta_N^{n_k^k,j}.\ee
Since obviously $\zeta_N^{n_k^k,j}\leqslant\zeta_{N+1}^{n_k^k,j}$, we have
$\zeta_N^{n_k^k}\leqslant\zeta_{N+1}^{n_k^k}$ and hence $\zeta_N\leqslant\zeta_{N+1}$.
We will show that for each $N\in\N^+$, as $k\ra\8$,
the sequence $(\hat U^{n_k^k,k},\hat W^{n_k^k,k})$ converges almost surely to the pair
$(\hat U^N,\hat W^N)$ on $[0,\zeta_N\wedge T]$,
and correspondingly, its law sequence $\mu^{n_k^k,k}$ weakly converges to $\hat\mu^N$.

By Lemma \ref{le4.5}, Skorohod's Theorem and \eqref{4.55},
we know whenever $1\leqslant N<N'\leqslant k$,
for almost surely $t\in[0,\zeta^{n_k^k}_N\wedge T]$, it holds that
\be\label{4.58}
(\hat U^{n_k^k,N'}(t),\hat W^{n_k^k,N'}(t))=
(\hat U^{n_k^k,N}(t),\hat W^{n_k^k,N}(t)).
\ee
Since the right side of \eqref{4.58} is indeed a subsequence of
$(\hat U^{n_k^N,N}(t),\hat W^{n_k^N,N}(t))$,
by the uniqueness of almost surely limit and picking all $k\in\N^+$, we conclude that
$$(\hat U^{N'}(t),\hat W^{N'}(t))=(\hat U^N(t),\hat W^N(t))\hs\mb{for all }N'>N$$
for almost every $t\in[0,\zeta_N\wedge T]$ and $\w\in\hat\W$ and due to
the almost surely equality we also have that $\hat\mu^{N'}=\hat\mu^{N}$ for $N'>N$.
These results have already implied the two convergences we desire.

We now claim that
\be\label{4.59}\zeta_N\ra\8\Hs\mb{ almost surely as }N\ra+\8.\ee
By assuming the claim, one can easily extend the pair $(\hat U^N,\hat W^N)$ on $[0,\zeta_N\wedge T]$
to $(\hat U,\hat W)$ defined on $[0,T]$ with $(\hat U(t),\hat W(t))=(\hat U^N(t),\hat W^N(t))$ on
$[0,\zeta_N\wedge T]$ for each $N\in\N^+$ and $T>0$.
Correspondingly, for each $T>0$, the law $\hat\mu$ of $(\hat U,\hat W)$ is equal to $\hat\mu^N$
for each $N\in\N^+$, and the stochastic basis
$\hat\cS=(\hat\W,\hat\cF,\{\hat\cF\}_{t\geqslant0},\hat\bP,\hat W)$ can be similarly obtained
by extending each $\hat\cF^N$ to $\hat\cF$.
Since each $(\hat\cS^N,\hat U^N,\hat\rho_N)$ is a local martingale solution of \eqref{2.13}
and $\hat\rho_N\leqslant\zeta_N$, then $(\hat\cS,\hat U,\zeta_N)=(\hat\cS^N,\hat U^N,\zeta_N)$ is
also a local martingale solution.
Due to \eqref{4.59}, we see that $(\hat\cS,\hat U)$ is a global martingale solution of \eqref{2.13}.

Eventually, it is sufficient to show the claim \eqref{4.59}.
We show it by contradiction and assume there is $\hat\W_0\subset\hat\W$ with $\hat\bP(\hat\W_0)>0$
such that $\zeta_N(\w)$ is bounded for all $\w\in\hat\W_0$ and $N\in\N^+$.
Observe that
$$\hat\W_0=\bigcup_{k\in\N^+}\hat\W_{k}\hs\mb{with }
\hat\W_k:=\{\w\in\hat\W_0:\sup_{N\in\N^+}\zeta_N(\w)\leqslant k\}.$$
Surely, there exists $k_0\in\N^+$ so that
$\ve_0:=\hat\bP(\hat\W_{k_0})>0$.
Fix each $N\in\N^+$.
By the first definition in \eqref{4.57}, we can pick a subsequence $\{k_i\}_{\geqslant1}$ of
$\{k\}_{k\geqslant1}$ so that
$\zeta^{n_{k_i}^{k_i}}_N<k_0+1$ almost surely in $\hat\W_{k_0}$.
Then by the second definition in \eqref{4.57}, we deduce that for some $N_{k_i}\geqslant N$,
$$\sup_{s\in[0,k_0+1]}\|\hat U^{n_{k_i}^{k_i},N_{k_i}}(s)\|_2\geqslant N
\hs\mb{occurs almost surely on }\hat\W_{k_0}.$$
As a result, it follows from \eqref{4.34} and the embeddings that for $i=\N^+$
and $p\geqslant2$,
\begin{align}N^p\ve_0\leqslant&
  \int_{\hat\W_{k_0}}\sup_{s\in[0,k_0+1]}\|\hat U^{n_{k_i}^{k_i},N_{k_i}}(s)\|_2^p\hat\bP(\w)
\leqslant
  \bE\sup_{s\in[0,k_0+1]}\|\hat U^{n_{k_i}^{k_i},N_{k_i}}(s)\|_2^p\notag\\
\leqslant&
  c_{m,p}\bE\sup_{s\in[0,k_0+1]}\|A^m\hat U^{n_{k_i}^{k_i},N_{k_i}}(s)\|^p
\leqslant
  C'_9,\label{4.60}
\end{align}
where $C'_9>0$ is independent of all $n_{k_i}^{k_i}$ and $N_{k_i}$.
However, the $N$ on the left hand of \eqref{4.60} can be chosen arbitrarily,
which contradicts the boundedness of the right hand.
The proof is finished.\eo

\br The procedure we adopt in this subsection can be applied to much more examples of stochastic
partial differential equations with nonlinear forcing terms and
globally Lipschitz diffusion coefficients, to establish the global existence of martingale solutions.
The equation is required to be dissipative to some extent.
More importantly, we need not to be aware of the existence of global pathwise solutions here.
\er
\section{Local and Global Existence of Pathwise Solutions}\label{s5}

We are now ready to study the existence of pathwise solutions of \eqref{2.13}.
Recalling Proposition \ref{pro2.1}, we will show that the sequence of solutions $u^n$ of \eqref{3.1}
converges almost surely in $L^2(0,T;H^{2(m+1)})\cap\cC([0,T];H^{2m})$
relative to initial stochastic basis.
We first deal with the condition \eqref{2.17}, which can be interpreted by pathwise uniqueness.

\subsection{Local existence of pathwise solutions}

We will adopt the Yamada-Wannabe theorem (see e.g. \cite{GHV14}) to give the local existence of
pathwise solutions.
Accordingly, The first result to be proved is the local pathwise uniqueness for an arbitrary pair
of solutions of the cut-off system \eqref{3.1}.

\bl\label{le5.1}Assume that $\phi$ satisfies \eqref{2.7}.
Suppose $(\cS,u^{(1)})$ and $(\cS,u^{(2)})$ are two global solutions of \eqref{3.1} given by
Theorem \ref{th4.2}, relative to the same stochastic basis
$\cS:=\{\W,\cF,\{\cF_t\}_{t\geqslant0},\bP,W\}$.
Define
\be\label{5.1}\W_0=\{u^{(1)}(0)=u^{(2)}(0)\}.\ee
Then $u^{(1)}$ and $u^{(2)}$ are indistinguishable on $\W_0$, i.e.,
\be\label{5.2}\bP\((u^{(1)}(t)-u^{(2)}(t))\chi_{\W_0}=0;\,\mb{for all }t\geqslant 0\)=1.\ee
\el
\bo Define $V=u^{(1)}-u^{(2)}$ and $\ol{V}=\chi_{\W_0}V$.
By Definition \ref{de2.1}, we know that
$$\ol{V}\in\cC([0,\8);H^{2m})\cap L_{\rm loc}^2(0,\8;H^{2(m+1)})\mb{ almost surely.}$$
Define a stopping time $\ol\varrho_R$ such that
$$\ol\varrho_R=\inf\{t\geqslant0:\|u^{(1)}(t)\|_{2m}+\|u^{(2)}(t)\|_{2m}>R\}.$$
Trivially, it can be seen that $\ol\varrho_R\ra+\8$, as $R\ra+\8$.
And the uniqueness follows immediately from the conclusion
\be\label{5.2a}\bE\sup_{s\in[0,T\wedge\ol\varrho_R]}\|\ol V(s)\|_{2m}^p=0.\ee

Now we are devoted to proving \eqref{5.2a}.
Subtracting the equations \eqref{3.1} for $u^{(2)}$ from that for $u^{(1)}$,
we get the following equation for $V$,
\be\di V+[A^2V+\Phi_1(u^{(1)},u^{(2)})]\di t=\Phi_2(u^{(1)},u^{(2)})\di W,\label{5.3}\ee
where
$$\Phi_1(u^{(1)},u^{(2)}):=\de_N(\|u^{(1)}\|_2)f(u^{(1)})-\de_N(\|u^{(2)}\|_2)f(u^{(2)}),$$
$$\Phi_2(u^{(1)},u^{(2)}):=\de_N(\|u^{(1)}\|_2)\phi(u^{(1)})-\de_N(\|u^{(2)}\|_2)\phi(u^{(2)}).$$
Apply $A^m$ to \eqref{5.3}.
It yields from It\^o's Formula for $\|A^mV\|^p$ that
for each fixed $T>0$ and $\tau'$, $\tau''$ with
$0\leqslant\tau'<s\leqslant\tau''\leqslant T\wedge\ol\varrho_R$,
we have
\begin{align}&
  \frac1p\|A^mV(s)\|^p+\int_{\tau'}^{s}\|A^mV\|^{p-2}\|A^{m+1}V\|^2\di s'\notag\\
\leqslant&
  \frac1p\|A^mV(\tau')\|^p-\int_{\tau'}^{s}\|A^mV\|^{p-2}
  \(A^mV,A^m\Phi_1(u^{(1)},u^{(2)})\)\di s'\notag\\
&
  +\frac{p-1}2\int_{\tau'}^{s}\|A^mV\|^{p-2}\|A^m\Phi_2(u^{(1)},u^{(2)})\|^2_{L_2(\fU,H)}\di s'\notag\\
&
  +\int_{\tau'}^{s}\|A^mV\|^{p-2}\(A^mV,A^m\Phi_2(u^{(1)},u^{(2)})\di W\)\notag\\
:=&
  \frac1p\|A^mV(\tau')\|^p+M'_1+M'_2+M'_3.\label{5.4}
\end{align}

Noting that $\de_N$ is globally Lipschitz (depending on $N$) and similar to \eqref{4.47},
we deduce  that
\begin{align*}&|M'_1|\\
\leqslant&
  \int_{\tau'}^{s}\|A^mV\|^{p-2}\left|\(A^{m+1}V,
  (\de_N(\|u^{(1)}\|_2)-\de_N(\|u^{(2)}\|_2))A^{m-1}f(u^{(1)})\)\right|\di s'\\
&
  +\int_{\tau'}^{s}\|A^mV\|^{p-2}\left|\(A^{m+1}V,A^{m-1}
  (f(u^{(1)})-f(u^{(2)}))\)\right|\di s'\\
\leqslant&
  \frac{2p-1}{2p}\int_{\tau'}^{s}\|A^mV\|^{p-2}\|A^{m+1}V\|^2\di s'
  +c_{m,p}\int_{\tau'}^{s}\|A^mV\|^{p}
  \|A^{m-1}f(u^{(1)})\|^2\di s'\\
&
  +c\int_{\tau'}^{s}\|A^mV\|^{p}\di s'
  +c_p\!\int_{\tau'}^{s}\|A^mV\|^{p-2}
  \left\|A^{m-1}\!\(b\na V\cdot\na(u^{(1)}+u^{(2)})\)\right\|^2\di s'\\
&
  +c_p\int_{\tau'}^{s}\|A^mV\|^{p-2}\left\|A^{m-1}
  \(V\((u^{(1)})^2+u^{(1)}u^{(2)}+(u^{(2)})^2\)\)\right\|^2\di s'.
\end{align*}
For the norm $\|A^{m-1}f(u^{(1)})\|$, we recall the discussions to \eqref{3.10}
and \eqref{3.13} and consider the embeddings into $H^{2m}$.
Then we know that
$$\|A^{m-1}f(u^{(1)})\|\leqslant
  c\(\|u^{(1)}\|_{2m}+\|u^{(1)}\|_{2m}^2+\|u^{(1)}\|_{2m}^3\).
$$
Using similar estimates to \eqref{4.48} and \eqref{4.49}, we obtain that
\begin{align}|M'_1|\leqslant&
  \frac{2p-1}{2p}\int_{\tau'}^{s}\|A^mV\|^{p-2}\|A^{m+1}V\|^2\di s'\notag\\
&  +c_{m,p}\int_{\tau'}^{s}\|A^mV\|^{p}\(\|u^{(1)}\|_{2m}^4+\|u^{(2)}\|_{2m}^4+1\)\di s'.
\label{5.5}\end{align}
For $M'_2$ and $M'_3$, we necessarily notice the following estimate,
\begin{align*}&
  \|\de_N(\|u^{(1)}\|_2)A^m\phi(u^{(1)})
  -\de_N(\|u^{(2)}\|_2)A^m\phi(u^{(2)})\|^2_{L_2(\fU,H)}\\
\leqslant&
  c\left\|\(\de_N(\|u^{(1)}\|_2)-\de_N(\|u^{(2)}\|_2)\)A^m\phi(u^{(1)})\right\|_{L_2(\fU,H)}^2\notag\\
&+c\|A^m\phi(u^{(1)})-A^m\phi(u^{(2)})\|_{L_2(\fU,H)}^2\\
\leqslant&
  c\|A^mV\|^2\kappa^2(|u^{(1)}|_\8)\(1+\|u^{(1)}\|_{2m}^2\)
  +c\kappa^2(|u^{(1)}|_\8+|u^{(2)}|_\8)\|A^mV\|^2\\
\leqslant&
  c\|A^mV\|^2\(\kappa^2(|u^{(1)}|_\8)\(1+\|u^{(1)}\|_{2m}^2\)
  +\kappa^2(|u^{(1)}|_\8+|u^{(2)}|_\8)\)\\
:=&
  c\|A^mV\|^2\Psi(u^{(1)},u^{(2)}).
\end{align*}
Then
\be\label{5.8}
|M'_2|\leqslant c_{p}\int_{\tau'}^s\|A^mV\|^p\Psi(u^{(1)},u^{(2)})\di s'.
\ee
For $M'_3$, we use BDG inequality as before and obtain that
\begin{align}&\bE\(\chi_{\W_0}\sup_{s\in[\tau',\tau'']}|M'_3(s)|\)
\leqslant
  c\bE\(\int_{\tau'}^{\tau''}\|A^m\ol V\|^{2p}\Psi(u^{(1)},u^{(2)})\di s\)^{\frac12}\notag\\
\leqslant&
  c\bE\(\sup_{s\in[\tau',\tau'']}\|A^m\ol V(s)\|^{p}\)^{\frac12}
  \(\int_{\tau'}^{\tau''}\|A^m\ol V\|^{p}\Psi(u^{(1)},u^{(2)})\di s\)^{\frac12}\notag\\
\leqslant&
  \frac1{2p}\bE\sup_{s\in[\tau',\tau'']}\|A^m\ol V(s)\|^{p}
  +c_{p}\bE\int_{\tau'}^{\tau''}\|A^m\ol V\|^{p}\Psi(u^{(1)},u^{(2)})\di s.\label{5.9}
\end{align}
Taking the expected value of the product of $\chi_{\W_0}$ and the supremum of \eqref{5.4}
over $s\in[\tau',\tau'']$, and combining \eqref{5.5}, \eqref{5.8} and \eqref{5.9}, we have
\begin{align}&\bE\sup_{s\in[\tau',\tau'']}\|A^m\ol V(s)\|^p
  +\bE\int_{\tau'}^{\tau''}\|A^m\ol V\|^{p-2}\|A^{m+1}\ol V\|^2\di s\notag\\
\leqslant&
  2\bE\|A^m\ol V(\tau')\|^p+c_{m,p}\bE\int_{\tau'}^{\tau''}\|A^m\ol V\|^p
  \(1+\|u^{(1)}\|_{2m}^4+\|u^{(2)}\|_{2m}^4\)\di s\notag\\
&
  +c_{m,p}\bE\int_{\tau'}^{\tau''}\|A^m\ol V\|^p\Psi(u^{(1)},u^{(2)})\di s.\label{5.10}
\end{align}
We apply the stochastic Gronwall's inequality to \eqref{5.10}, so that
\begin{align*}&\bE\(\sup_{s\in[0,T\wedge\ol\varrho_R]}\|A^m\ol V(s)\|^p
+\int_{0}^{T\wedge\ol\varrho_R}\|A^m\ol V\|^{p-2}\|A^{m+1}\ol V\|^2\di s\)\\
\leqslant&C'_{10}\bE\|A^m\ol V(0)\|^p=0,\end{align*}
where the positive $C'_{10}$ only depends on $m$, $p$ and $N$ via $\de_N$, $\kappa$ and the embeddings.
This implies \eqref{5.2a} and completes the proof.
\eo

Indeed, Lemma \ref{le5.1} tells us the uniqueness of local pathwise solutions $(u,\rho_N)$
(defined in Definition \ref{de2.2}), where $\rho_N$ is defined in \eqref{4.7}
with $\tilde u$ therein replaced by $u$.

Next we continue following Proposition \ref{pro2.1},
for which we recall the Galerkin solutions $u^n$
of \eqref{3.2} with the given stochastic basis $\cS$.
Let $\mu_u^{n,m}$ be the collection of joint distributions given by $(u^n,u^m)$ and
define the extended phase space
$$\cX_u^J=\cX_u\X\cX_u,\Hs\cX^J=\cX_u^J\X\cX_W$$
with $\cX_u$ and $\cX_W$ defined in \eqref{4.1}.
Corresponding to $\cX_u^J$ and $\cX^J$ and following the definitions \eqref{4.2}
and \eqref{4.2A},
we ulteriorly define
$$\mu_u^{n,m}=\mu_u^n\X\mu_u^m,\Hs\mu^{n,m}=\mu_u^{n,m}\X\mu_W.$$
Then we similarly have the following tightness for $\{\mu^{n,m}\}_{n,m\geqslant0}$.
\bl Suppose that $\phi$ satisfies \eqref{2.7} and $u_0$ satisfies \eqref{2.11}.
Then the collection $\{\mu^{n,m}\}_{n,m\geqslant0}$ is tight on $\cX^J$.
\el
\bo The proof is a slight modification of that for Lemma \ref{le4.1}.
After determining the sets $B_R^1$ and $B_R^2$, we can take compact $\cA_\ve$ and $\cB_\ve$ in
$\cX_u$ and $\cX_W$ respectively such that $\mu_u^n(\cA_\ve)\geqslant1-\frac{\ve}{4}$
and $\mu_W^n(\cB_\ve)\geqslant1-\frac{\ve}{2}$ for each $n\in\N$.
Then the set $\cA_{\ve}\X\cA_{\ve}\X\cB_\ve$ is compact in $\cX^J$ and
$$\mu^{n,m}(\cA_{\ve}\X\cA_{\ve}\X\cB_\ve)\geqslant
(1-\frac{\ve}4)^2(1-\frac{\ve}2)\geqslant 1-\ve$$
for all $\ve>0$ and $m,n\in\N$. This is the tightness.
\eo

With the tightness result, we are prepared to show the existence of local pathwise solution
of \eqref{2.13}.

\bt\label{th5.3} Suppose that $\phi$ satisfies \eqref{2.7}.
Then the problem \eqref{2.13} possesses a unique local pathwise solution.
Moreover, each local pathwise solution can be extended to be a maximal one
$(u,\tau_0,\{\tau_n\}_{n\geqslant1})$,
and the sequence $\{\tau_n\}_{n\geqslant1}$ announces any finite time blowup.
\et
\bo
The tightness of the sequence $\{\mu^{n,m}\}_{n,m\geqslant0}$ enables us to choose a subsequence
so that $\mu^{n_k,m_k}$ converges weakly to an element $\mu$ by Prohorov's Theorem.
Similarly as above, we can use Skorohod's Theorem to attain the existence of a probability space
$(\tilde{\W},\tilde{\cF},\tilde{\bP})$, on which there defines a sequence of random elements
$(\tilde{u}^{n_k},\tilde{u}^{m_k},\tilde{W}^k)$ converging almost surely to an element
$(\tilde u,\tilde{u}^*,\tilde{W})$ with their laws satisfying
$$\tilde{\bP}(\tilde{u}^{n_k},\tilde{u}^{m_k},\tilde{W}^k)=\mu^{n_k,m_k},\Hs
\tilde{\bP}(\tilde{u},\tilde{u}^*,\tilde W)=\mu.$$

For the almost surely convergences $(\tilde{u}^{n_k},\tilde{W}^{k})\ra (\tilde{u},\tilde{W})$
and $(\tilde{u}^{m_k},\tilde{W}^{k})\ra (\tilde{u}^*,\tilde{W})$, by similar argument in the proof of
Theorem \ref{th4.2}, we infer that both $\tilde{u}$ and $\tilde{u}^*$ are global martingale solutions of
\eqref{3.1} relative to the same stochastic basis
$\tilde\cS=\{\tilde\W,\tilde\cF,\{\tilde\cF_t\}_{t\geqslant0},\tilde\bP,\tilde{W}\}$
with $\tilde{\cF}_t$ the completion of $\sigma$-algebra generated by
$\{(\tilde{u}(s),\tilde{u}^*(s),\tilde{W}(s)):s\leqslant t\}$.

Note that $\mu^{n_k,m_k}_u$ converges weakly to the measure $\mu_u^J$ on $\cX^J$ with
$\mu_u^J(\cdot)=\tilde{\bP}((\tilde{u},\tilde{u}^*)\in\cdot)$.
Since the choice of $u^n$ and $u^m$ has ensured the equality $\tilde{u}(0)=\tilde{u}^*(0)=u_0$
by the convergence,
by Lemma \ref{le5.1}, we deduce that $\tilde{u}=\tilde{u}^*$ in $\cX_u$ almost surely, saying,
$$\mu^J(\{(x,y)\in\cX_u^J:x=y\})=\tilde{\bP}(\tilde{u}=\tilde{u}^*\mb{ in }\cX_u)=1.$$
This indeed indicates that the original sequence $u^n$ given in Section \ref{s3} defined on the initial
probability space $\W,\cF,\bP$ converges to an element $u$, in the topology of $\cX_u$.
Then $u$ is a global pathwise solution of \eqref{3.1}.
\vs

Based on this result,
we choose a appropriate strictly positive stopping time $\tau$ to construct the local pathwise
solution of \eqref{2.13}.
Note that $u_0\in H^{2m}$ almost surely.
First we assume $\|u_0\|_{2m}\leqslant\cM$ for some deterministic $\cM>0$.
Then by embeddings, $\|u_0\|_2\leqslant c\cM$ for some $c>0$.
Let $u$ be the global pathwise solution of \eqref{3.1}, which is in $L^2(\W;\cC([0,\8);H^{2}))$
by embeddings.
To deduce the maximal pathwise solution, we recall the stopping time $\rho_N$ defined in \eqref{4.7}.
Since $u\in\cC([0,+\8);H^{2})$, if we take $N>c\cM$, we can easily see that $\rho_N>0$.
Thus let $\tau=\rho_N$ and the pair $(u,\tau)$ is just the local pathwise solution of \eqref{2.13}.

For the general case when $u_0$ is an $H^{2m}$-valued random variable and $\cF_0$-measurable,
we proceed as what appeared in \cite[Section 4.2]{GHZ09}.
Define, for each $k\geqslant0$, $u_{0}^{\langle k\rangle}=u_0\chi_{\{k\leqslant\|u_0\|_{2m}<k+1\}}$.
By the discussion above, we get a corresponding local pathwise solution
$(u^{\langle k\rangle},\tau^{\langle k\rangle})$.
Thus we obtain a local pathwise solution $(u,\tau)$ for the initial datum $u_0$ by defining
$$u=\sum_{k\geqslant0}u^{\langle k\rangle}\chi_{\{k\leqslant\|u_0\|_{2m}<k+1\}},\hs
\tau=\sum_{k\geqslant0}\tau^{\langle k\rangle}\chi_{\{k\leqslant\|u_0\|_{2m}<k+1\}}.$$
\vs

Next, for every fixed $u_0\in H^{2m}$, we go on to extend the solution $(u,\tau)$ to a maximal time
(referring the method in \cite{GHZ09,GHV14}).
Let $\cT$ be the set of all stopping times $\tau$ for a local pathwise solution of \eqref{2.13}.
Take $\tau_0=\sup\cT$.
By the uniqueness of pathwise solutions we can obtain a process $u$ defined on $[0,\tau_0)$
such that $(u,\tau)$ is a local pathwise solution for each stopping time $\tau\in(0,\tau_0)$.

In the following, we show the existence of the strictly positive increasing stopping time sequence
$\{\tau_n\}_{n\geqslant1}$ such that the triple $(u,\tau_0,\{\tau_n\}_{n\geqslant1})$ is
a maximal pathwise solution of \eqref{2.13}.
For each $n\in\N^+$, take
\be\label{5.11}\tilde{\rho}_n=\inf\{t\geqslant0:\|u(t)\|_2>n\}\wedge\tau_0.\ee
The continuity of $u$ on $H^2$ ensures that $\tilde{\rho}_n$ is a well-defined stopping time.
Moreover, by uniqueness we know that $(u,\tilde{\rho}_n)$ is a local pathwise solution
on $\{\tilde{\rho}_n>0\}$ for each $\w\in\W$ and $n>0$ (Note that each $\w$ has sufficiently large
$N$ so that $\tilde{\rho}_n(\w)>0$).
Suppose by contradiction that, for some $n,\,T>0$, we have $\bP(\tau_0=\tilde{\rho}_n\wedge T)>0$.
Since $\tau_0>0$, we see that $(u,\tilde{\rho}_n\wedge T)$ is a local pathwise solution.
By the unique existence of local pathwise solutions established above, we have another stopping time
$\tau'>\tilde{\rho}_n\wedge T$ and a process $u'$ such that $(u',\tau')$ is a local pathwise solution
corresponding to $u(\tilde{\rho}_n\wedge T)$ and then also to $u_0$.
This denies the maximality of $\tau_0$.

We have shown that $\bP(\tau_0=\tilde{\rho}_n\wedge T)=0$, for all $n,\,T>0$.
Actually, on the set $\{\tau_0<\8\}$, it can be seen that $\tilde\rho_n<\tau_0$
for all $n>0$ almost surely by appropriate choices of $T$.
By \eqref{5.11}, we infer that $\sup_{s\in[0,\tilde{\rho}_n]}\|u(s)\|_2=n$ and so as $n\ra\8$,
$$\sup_{s\in[0,\tilde{\rho}_n]}\|u(s)\|_{2}^2+\int_0^{\tilde{\rho}_n}\|A^2u(s)\|^2\di s\ra+\8.$$
Now let $\tilde\tau_n$ be a strictly positive increasing stopping time with $\tilde{\tau}_n<\tau_0$ and
$\tilde{\tau}_n\ra\8$ as $n\ra\8$.
We have known $(u,\tilde{\tau}_n)$ is a local pathwise solution for each $n$.
Then the sequence $\tau_n:=\tilde{\rho}_n\vee\tilde{\tau}_n$ is the sequence we require.
We have finally accomplished the proof now.
\eo
\subsection{Global existence of pathwise solutions}

In this subsection we establish the global existence of pathwise solution of \eqref{2.13}.
To this aim, we need to assure that for the maximal pathwise solution $(u,\tau_0,\{\tau_n\}_{n\geqslant1})$,
the probability $\bP(\tau_0=\8)=1$.
This is stated in the following theorem.

\bt\label{th5.4} Let $\cS=(\W,\cF,\{\cF_t\}_{t\geqslant0},\bP,W)$ be a stochastic basis.
Suppose that $u_0$ is an $H^{2m}$-valued ($m\geqslant1$) random variable (relative to $\cS$)
and $\cF_0$-measurable, $\phi$ satisfies \eqref{2.3} and $|b|<4$.
Let $(u,\tau_0,\{\tau_n\}_{n\geqslant1})$ be the maximal pathwise solution of \eqref{2.13}.
Then $\bP(\tau_0=\8)=1$ and $(u,\tau_0)$ is naturally a unique global pathwise solution of \eqref{2.13}
corresponding to the initial datum $u_0$.
\et

To show Theorem \ref{th5.4}, we need a small result as follows.
\bl\label{le5.5} Let $(\W,\mu)$ be a positive measure space and $\W_0\subset\W$
be a measurable subset of $\W$ with $0<\mu(\W_0)<\8$.
Suppose that $\{f_k\}_{k\geqslant1}$ is a sequence of nonnegative $\mu$-integrable functions
defined on $\W$ such that $f_k(\w)\ra+\8$ as $k\ra\8$ for almost every $\w\in\W_0$.
Then
\be\label{5.11A}\int_{\W}f_k(\w)\mu(\di\w)\ra+\8,\Hs\mb{as }k\ra\8.\ee
\el
\bo It is sufficient to show that for each $R>0$, there exists $K>0$ so that when $k\geqslant K$,
$\disp\int_{\W_0}f_k(\w)\mu(\di\w)>R$.
Indeed, if we set $g_k(\w)=\arctan f_k(\w)$ for each $w\in\W_0$, we know that $g_k(\w)\ra\frac{\pi}{2}$
from below for a.e. $\w\in\W_0$.
By Egorov's theorem, we know that there is a measurable subset $\W'_0$ of $\W$ with
$\mu(\W'_0)=\mu(\W_0)/2$, such that $\{g_k\}_{k\geqslant1}$ uniformly converges to $\frac{\pi}2$
when restricted to $\W'_0$.
This means that we can find $K>0$ such that whenever $k\geqslant K$,
$$g_k(\w)>\arctan\frac{2R}{\mu(\W_0)},\Hs\mb{for every }\w\in\W'_0,$$
and hence $f_k(\w)>\frac{2R}{\mu(\W_0)}$.
As a result, we obtain that
$$\int_{\W_0}f_k(\w)\mu(\di\w)\geqslant\int_{\W'_0}f_k(\w)\mu(\di\w)
>\frac{2R}{\mu(\W_0)}\cdot\frac{\mu(\W_0)}{2}=R.$$
The proof is finished.\eo

\br Lemma \ref{le5.5} is a simple result in measure theory and should have appeared somewhere,
but we have not found any one in other references.
For the completeness of this paper, we give a concise proof above.
\er

Now we go on to verify Theorem \ref{th5.4}.

\noindent\textit{Proof of Theorem \ref{th5.4}.} The uniqueness of local pathwise solution has ensured
the uniqueness of the global pathwise solution. We hence only need to show $\bP(\tau_0=\8)=1$.

Since $(u,\tau_n)$ is a local pathwise solution \eqref{2.13} for each $n\in\N^+$,
we can apply the estimation method in the proof of Lemma \ref{le4.4} so that for each $T>0$, $n\in\N^+$
and $0\leqslant\tau'\leqslant\tau''\leqslant\tau_n\wedge T$,
\begin{align*}
  &\bE\sup_{s\in[\tau',\tau'']}\|u(s)\|^{2(2m+3)}\\
  \leqslant&
  2\bE\|u(\tau')\|^{2(2m+3)}+c_m\bE\int_{\tau'}^{\tau''}\(1+\|u\|^{2(2m+3)}\)\di s
  +c(\tau''-\tau'),
\end{align*}
similar to \eqref{4.41} by observing that essentially $\de_N$ does not involve the estimation
in that proof except the maximum of $\de_N$ being $1$.
Define a subspace $D_r$ of $\W$ for each $r>0$ such that
$$D_r:=\{\w\in\W:\|u_0\|_2^2+\|u_0\|^{2(2m+3)}\leqslant r\}.$$
Using stochastic Gronwall's inequality, we similarly obtain a positive constant $C'_{11}=C'_{11}(r)$
independent of $n$, such that
\be\label{5.12}
\bE\sup_{s\in[0,\tau_n\wedge T]}\|u(s)\|^{2(2m+3)}\chi_{D_r}\leqslant C'_{11}\hs\mb{and}\ee
\be\label{5.12A}
\bE\(\sup_{s\in[0,\tau_n\wedge T]}\|u(s)\|_2^2+\int_0^{\tau_n\wedge T}
\|A^2u\|^2\di s\)\chi_{D_r}\leqslant C'_{11}.
\ee

Now we show the conclusion $\bP(\tau_0=\8)=1$ by contradiction and suppose the contrary.
Then there is a subset $\W_0\subset\W$ with $\bP(\W_0)>0$ such that for all $\w\in\W_0$,
$\tau_0(\w)<\8$.
Note that
$$\W_0=\bigcup_{k\in\N^+}\W_k\hs\mb{with }\W_k=\{\w\in\W_0:\tau_0(\w)\leqslant k\}.$$
We immediately have $r_0\in\R^+$ and $k_0\in\N^+$ so that $\ve_0:=\bP(\W_{k_0}\cap D_{r_0})>0$.
Apparently, $\tau_n(\w)\leqslant k_0$ for all $\w\in\W_{k_0}\cap D_{r_0}$ and $n\in\N^+$.
According to \eqref{2.16} in the definition of maximal pathwise solution, we can pick $T=k_0$ and
obtain that for almost every $\w\in\W_{k_0}$, \eqref{2.16} holds with $\tau_n$ therein replaced by
$\tau_n\wedge T$.
Then taking $n\ra\8$, we see by Lemma \ref{le5.5} the left side of the second inequality tends
to the infinity, but $C'_{11}$ is independent of $n$, which causes a contradiction and finishes the proof.
\qed

\section{Existence of Ergodic Invariant Measures}\label{s6}

In this section, we consider the existence of ergodic invariant measures
for the 2-dimensional stochastic modified Swift-Hohenberg periodic problem \eqref{1.1} -- \eqref{1.3}.
In what follows in this section, let $\cS=(\W,\cF,\{\cF_t\}_{t\geqslant0},\bP,W)$ be a stochastic basis.
We assume that $|b|<4$ and $\phi$ satisfies \eqref{2.3} and \eqref{2.4}.
This is sufficient to guarantee the unique existence of global pathwise solutions by Theorem \ref{th5.4}
for each $H^{2m}$-valued random variable $u_0$ measured by $\cF_0$.

Let $u(t;u_0)$ be the unique global pathwise solution in $H^{2m}$ of the problem \eqref{2.13}.
This defines a stochastic process $u(t;\cdot)$ on the space $H^{2m}$.
For notational convenience, we often use $u(t)$ to denote
stochastic process $u(t;\cdot)$.
Let $\cP_t$ be the Markovian transition semigroup associated with the stochastic process
$u(t)$ on $\cB(H^{2m})$.
In order to prove the existence of invariant measures for the semigroup $\cP_t$,
we first present some auxiliary consequences of
the estimates of moment bounds and some related probabilities for the solution $u(t;u_0)$.

\subsection{Estimates on moment bounds of solutions}
In the following,
we let the initial data $u_0$ be an $H^{2m}$-valued random variable.
For each $u_0$ given above and $T>0$, recall the deduction of \eqref{5.12} and \eqref{5.12A}.
In that $u(t;u_0)$ is a given global pathwise solution of \eqref{2.13}, we can repeat the procedure
from \eqref{4.38} to \eqref{4.43} by replacing $u^n$ and $\de_N$ therein by $u$ and $1$.
Thus we can obtain similar estimates to \eqref{4.42} and \eqref{4.34} as follows,
for each $t\in[0,T]$ and $u_0$ satisfying \eqref{2.12} with $p=2$,
\begin{align}
&\bE\(\sup_{s\in[0,t]}\|u(s)\|^{2(2m+3)}+\int_{0}^t\|u\|^{2(2m+2)}
\|Au\|^2\di s\)\notag\\
\leqslant&C_4\(\bE\|u_0\|^{2(2m+3)}+t\),\label{6.2}\end{align}
\begin{align}
&\bE\(\sup_{s\in[0,t]}\|A^mu(s)\|^2+\int_{0}^t
\|A^{m+1}u\|^2\di s\)\notag\\
\leqslant&
C_4\(\bE\|u_0\|_{2m}^2+\bE\|u_0\|^{2(2m+3)}+t\),
\label{6.3}\end{align}
where $C_4=C_4(T)$ is a positive constant.

Define a stopping time $\xi_r(u_0)$ for each $r>0$ and $H^{2m}$-valued stochastic variable $u_0$ so that
\be\label{6.4}\xi_r(u_{0}):=\inf\{t\geqslant0:\|A^mu(t;u_{0})\|^2>r\}.\ee
In order to obtain a description of the instant parabolic regularization for the solutions of \eqref{2.13},
we define another stopping time $\eta_r(u_0)$ for each $r>0$ and
$H^{2m}$-valued stochastic variable $u_0$ such that
$$\eta_r(u_0)=\inf\{t\geqslant 0:t\|A^{m+1}u(t;u_0)\|^2>r\}.$$
To estimate the probability concerning $\eta_r$, we first give a suitable local moment bound
on the $H^{2(m+1)}$-norm of the solutions.

\bl\label{le6.1}Let $\phi$ satisfy \eqref{2.3} and \eqref{2.4}.
For each $T>0$ and $u_0$ satisfying \eqref{2.12} with $p=2$ and $q'\geqslant2(2m+5)$,
it holds that
\be\label{6.5}\bE\sup_{s\in[0,t]}\(s\|A^{m+1} u(s)\|^2\)
\leqslant C_5\(\bE\|u_0\|_{2m}^2+\bE\|u_0\|^{2(2m+5)}+t\),\ee
for all $t\in[0,T]$, where $C_5:=C_5(T)$ is a positive constant.
\el

\bo Apply $A^{m+1}$ to \eqref{2.13} and It\^o's Formula to $s\|A^{m+1}u(s)\|^2$,
we deduce that for stopping times $0\leqslant\tau'\leqslant \vsig\leqslant s\leqslant\tau''
\leqslant t\leqslant T$,
\begin{align} &s\|A^{m+1}u(s)\|^2+2\int_{\vsig}^ss'\|A^{m+2}u\|^2\di s'\notag\\
=&
  \vsig\|A^{m+1}u(\vsig)\|^2-\int_{\vsig}^s
  \left[2s'\(A^{m+1}u,A^{m+1}f(u)\)-\|A^{m+1}u\|^2\right]\di s'\notag\\
&
  +\int_{\vsig}^ss'\|A^{m+1}\phi(u)\|_{L_2(\fU,H)}^2\di s'
  +2\int_{\vsig}^ss'(A^{m+1}u,A^{m+1}\phi(u)\di W).\label{6.6}
\end{align}
Similar to the treatment in the proof of Lemma \ref{le3.2}, we have by \eqref{2.4} that
\begin{align}\left|\(A^{m+1}u,A^{m+1}f(u)\)\right|\leqslant&
\|A^{m+2}u\|^2+c\(\|A^{m+1}u\|^2+\|u\|^{2(2m+5)}+1\)\\
\label{6.8}\|A^{m+1}\phi(u)\|_{L_2(\fU,H)}^2\leqslant&
c\(\|A^{m+1}u\|^2+1\).\end{align}
The estimates from \eqref{6.6} to \eqref{6.8} indicate that
\begin{align}&s\|A^{m+1}u(s)\|^2\notag\\
\leqslant&
  \vsig\|A^{m+1}u(\vsig)\|^2+c\int_{\vsig}^s
  s'\(\|A^{m+1}u\|^2+\|u\|^{2(2m+5)}+1\)\di s'\notag\\
&
  +\int_{\vsig}^s\|A^{m+1}u\|^2\di s'
  +2\left|\int_{\vsig}^ss'(A^{m+1}u,A^{m+1}\phi(u)\di W)\right|.
\label{6.9}\end{align}
Integrating \eqref{6.9} with respect to $\vsig$ over $[\tau',s]$ and dividing it by $s-\tau'$,
we have
\begin{align}&s\|A^{m+1}u(s)\|^2\notag\\
\leqslant&
  \frac{1}{s-\tau'}\int_{\tau'}^s\vsig\|A^{m+1}u(\vsig)\|^2\di\vsig+c\int_{\tau'}^s
  s'\(\|A^{m+1}u\|^2+\|u\|^{2(2m+5)}+1\)\di s'\notag\\
&
  +\int_{\tau'}^s\|A^{m+1}u\|^2\di s'
  +2\sup_{\vsig\in[\tau',s]}\left|\int_{\vsig}^ss'(A^{m+1}u,A^{m+1}\phi(u)\di W)\right|.
\label{6.10}\end{align}
Note also by BDG inequality that
\begin{align*}&
  \bE\sup_{s\in[\tau',\tau'']}\sup_{\vsig\in[\tau',s]}\left|\int_{\vsig}^s
  s'(A^{m+1}u,A^{m+1}\phi(u)\di W)\right|\notag\\
=&
  \bE\sup_{s\in[\tau',\tau'']}\sup_{\vsig\in[\tau',\tau'']}
  \(\left|\(\int_{\tau'}^s-\int_{\tau'}^{\vsig}\)
  s'(A^{m+1}u,A^{m+1}\phi(u)\di W)\right|\)\notag\\
\leqslant&
  2\bE\sup_{s\in[\tau',\tau'']}\left|\int_{\tau'}^s
  s'(A^{m+1}u,A^{m+1}\phi(u)\di W)\right|\notag\\
\leqslant&
  c\bE\(\int_{\tau'}^{\tau''}(s')^2\|A^{m+1}u\|^2
  (1+\|A^{m+1}u\|^2)\di s'\)^{\frac12}\notag\\
\leqslant&
  c\bE\(\sup_{s'\in[\tau',\tau'']}\(s'\|A^{m+1}u(s')\|^2\)\)^{\frac12}
  \(\int_{\tau'}^{\tau''}s'\(1+\|A^{m+1}u(s')\|^2\)\di s'\)^{\frac12}\\
\leqslant&
  \frac14\bE\sup_{s\in[\tau',\tau'']}\(s\|A^{m+1}u(s)\|^2\)
  +c\bE\int_{\tau'}^{\tau''}s'\(\|A^{m+1}u(s')\|^2+1\)\di s'.
\end{align*}
Thus, by considering the supremum of \eqref{6.10} over $s\in[\tau',\tau'']$,
and then taking the expected value, we get
\begin{align}
&\bE\sup_{s\in[\tau',\tau'']}\(s\|A^{m+1}u(s)\|^2\)\notag\\
\leqslant&
  c\bE\sup_{s\in[\tau',\tau'']}\frac{1}{s-\tau'}
  \int_{\tau'}^s\vsig\|A^{m+1}u(\vsig)\|^2\di\vsig
  +c\bE\int_{\tau'}^{\tau''}\|A^{m+1}u\|^2\di s'\notag\\
&+c\bE\int_{\tau'}^{\tau''}
  s'\(\|A^{m+1}u\|^2+\|u\|^{2(2m+5)}+1\)\di s'.\label{6.11}
\end{align}
By \eqref{6.2} and \eqref{6.3}, we know the right side of \eqref{6.11} makes sense
even when $\tau'=0$ and $\tau''=t$ and hence we immediately have
\begin{align*}
&\bE\sup_{s\in[0,t]}\(s\|A^{m+1}u(s)\|^2\)\\
\leqslant&
  c\bE\sup_{s\in[0,t]}\frac1s\int_{0}^s\vsig\|A^{m+1}u(\vsig)\|^2\di\vsig
  +c\bE\int_{0}^{t}\|A^{m+1}u\|^2\di s'\\
&+c\bE\int_{0}^{t}
  s'\(\|A^{m+1}u\|^2+\|u\|^{2(2m+5)}+1\)\di s'\\
\leqslant&c\bE\int_{0}^{t}\|A^{m+1}u\|^2\di s'
  +c\bE\int_{0}^{t}s'\(\|A^{m+1}u\|^2+\|u\|^{2(2m+5)}+1\)\di s'.
\end{align*}
which implies \eqref{6.5} by \eqref{6.2} and \eqref{6.3}.
The proof is over now.\eo

\br Note that $s\|A^{m+1}u(s)\|^2$ is not defined for $s=0$.
This prompts us to integrate \eqref{6.9} for $\vsig$ over $[\tau',s]$ and then we successfully
avoid the abuse of $\left.\(s\|A^{m+1}u(s)\|^2\)\right|_{s=0}$.
\er

Now we give the estimate of the probabilities of $\xi_r$ and $\eta_r$ in the following lemma.
\bl\label{le6.3} For each $t\in(0,T]$, $r>0$ and $u_0$ satisfying \eqref{2.12} with $p=2$
and $q'\geqslant2(2m+5)$,
\be\label{6.12}
\bP(\xi_r(u_0)<t)\leqslant\frac{C_4}{r}\(\bE\|u_0\|_{2m}^2+\bE\|u_0\|^{2(2m+3)}+t\)\ee
\be\label{6.13}\mb{and}\hs\bP(\eta_r(u_0)<t)\leqslant\frac{C_5}{r}
\(\bE\|u_0\|_{2m}^2+\bE\|u_0\|^{2(2m+5)}+t\),
\ee
where $C_4$ and $C_5$ are given by \eqref{6.3} and \eqref{6.5}, respectively.
\el
\bo By applying Markovian inequality, we see that for each $t\in(0,T]$,
\begin{align*}
&\bP(\xi_r(u_0)<t)
=\bP\(\sup_{s\in[0,t]}\|A^mu(s)\|^2>r\)\\
\leqslant&\frac{1}{r}\bE\sup_{s\in[0,t]}\|A^mu(s)\|^2\leqslant
\frac{C_4}{r}\(\bE\|u_0\|_{2m}^2+\bE\|u_0\|^{2(2m+3)}+t\),
\end{align*}
which is exactly \eqref{6.12}.
The inequality \eqref{6.13} follows in the same way.
\eo

\subsection{Feller property}

In this subsection we apply ourselves to confirming the Feller property of $\cP_t$,
for which we first extend the stopping times $\xi_r$ and $\eta_r$ to
relying on two initial data.
For arbitrary deterministic initial data $u_{10}$ and $u_{20}$, denote
$$\xi_r(u_{10},u_{20}):=\xi_r(u_{10})\wedge\xi_r(u_{20})\hs\mb{and}\hs
\eta_r(u_{10},u_{20}):=\eta_r(u_{10})\wedge\eta_r(u_{20}).$$
\bl\label{le6.4} For each fixed deterministic $T>0$, $r>0$ and $u_{10},\,u_{20}\in H^{2m}$,
it holds that
\be\label{6.14}\bE\sup_{s\in[0,t\wedge\xi_r(u_{10},u_{20})]}\|u(s,u_{10})-u(s,u_{20})\|^2_{2m}
\leqslant C_6\|u_{10}-u_{20}\|^2_{2m},\ee
for all $t\in[0,T]$, where $C_6=C_6(T,r)$ is a positive constant.
\el
\bo
Let $t\in[0,T]$, $u_i(t):=u(t,u_{i0})$, $i=1,2$ and $Z(t)=u_1(t)-u_2(t)$, for $t\in[0,T]$.
Then by \eqref{2.13} we have
\be\label{6.15}\di Z+[A^2Z+\tilde{\Phi}_1(u_1,u_2)]\di t=\tilde{\Phi}_2(u_1,u_2)\di W,\ee
where
$$\tilde{\Phi}_1(u_1,u_2):=f(u_1)-f(u_2)\mb{ and }\tilde{\Phi}_2(u_1,u_2):=\phi(u_1)-\phi(u_2).$$
Applying $A^m$ to \eqref{6.15} and It\^o's Formula to $\|A^mZ(t)\|^2$,
and similar to the argument in the proof of Lemma \ref{le4.5}, we have for all stopping times
$0\leqslant\tau'\leqslant s\leqslant\tau''\leqslant t\wedge\xi_r(u_{10},u_{20})$,
\begin{align*}
&\frac12\|A^mZ(s)\|^2+\int_{\tau'}^s\|A^{m+1}Z\|^2\di s'\\
\leqslant&
  \frac12\|A^mZ(\tau')\|^2+\left|\int_{\tau'}^{s}
  \(A^mZ,A^m\tilde{\Phi}_1(u_1,u_2)\)
  \di s'\right|\\
& +\frac12\int_{\tau'}^{s}\|A^m
  \tilde{\Phi}_2(u_1,u_2)\|^2_{L_2(\fU,H)}\di s'+\left|\int_{\tau'}^{t}\(A^mZ,A^m
  \tilde{\Phi}_2(u_1,u_2)\di W\)\right|\\
\leqslant&
  \frac12\|A^mZ(\tau')\|^2+\int_{\tau'}^{s}\|A^{m+1}Z\|^2\di s'
  +\left|\int_{\tau'}^{t}\(A^mZ,A^m\tilde{\Phi}_2(u_1,u_2)\di W\)\right|\\
&  +c\int_{\tau'}^s\|A^mZ\|^2\(\|u_1\|_{2m}^4+\|\|u_2\|_{2m}^4+1\)\di s'.
\end{align*}
Combining the estimate
\begin{align*}&\bE\sup_{s\in[\tau',\tau'']}\left|\int_{\tau'}^{s}
\(A^mZ,A^m\tilde{\Phi}_2(u_1,u_2)\di W\)\right|\\
\leqslant&\frac12\bE\sup_{s\in[\tau',\tau'']}\|A^mZ(s)\|^2
+c\bE\int_{\tau'}^{\tau''}\|A^mZ\|^2\di s',
\end{align*}
we similarly obtain
\begin{align*}&\bE\sup_{s\in[\tau',\tau'']}\|A^mZ(s)\|^2\\
\leqslant&
  2\bE\|A^mZ(\tau')\|^2+c\bE\int_{\tau'}^{\tau''}
\|A^mZ\|^2\(\|u_1\|_{2m}^4+\|u_2\|_{2m}^4+1\)\di s'.
\end{align*}
The stochastic Gronwall's inequality and embeddings apply and we obtain \eqref{6.14}.
The proof is finished.\eo

Now we can show the Feller property of the semigroup $\cP_t$.

\bt\label{th6.5} For each $\vp\in \cC_{\rm b}(H^{2m})$ and $t\geqslant0$,
the mapping $u\mapsto\cP_t\vp(u)$ is continuous.
\et
\bo Pick $u_{10},\,u_{20}\in H^{2m}$ with $u_{10}$ fixed.
Let $T>t$ be fixed and $r_1,\,r_2>0$.
Due to the embedding of $H^{2m}$ into $H$, we can fix $r_3>0$ such that when
$\|u_{10}-u_{20}\|_{2m}<r_3$,
$$\|u_{20}\|_{2m}^2\leqslant\|u_{10}\|_{2m}^2+1
\hs\mb{and}\hs
\|u_{20}\|^{2(2m+3)}\leqslant\|u_{10}\|^{2(2m+5)}+1.$$

Now we divide $|\cP_t\vp(u_{10})-\cP_t\vp(u_{20})|$ into three parts as follows,
\begin{align}&|\cP_t\vp(u_{10})-\cP_t\vp(u_{20})|=|\bE(\vp(u(t,u_{10}))-\vp(u(t,u_{20})))|\notag\\
\leqslant&
  |\bE(\vp(u(t,u_{10}))-\vp(u(t,u_{20})))\chi_{\{\xi_{r_1}(u_{10},u_{20})<t\}}|\notag\\
&  +|\bE(\vp(u(t,u_{10}))-\vp(u(t,u_{20})))\chi_{\{\eta_{r_2}(u_{10},u_{20})<t\}}|\notag\\
&
  +|\bE(\vp(u(t,u_{10}))-\vp(u(t,u_{20})))\chi_{\{\xi_{r_1}(u_{10},u_{20})\geqslant t\}}
  \chi_{\{\eta_{r_2}(u_{10},u_{20})\geqslant t\}}|\notag\\
:=&E_1+E_2+E_3.\label{6.16}
\end{align}
Let $\|\vp\|_\8:=\sup_{u\in H^{2m}}|\vp(u)|$ be the sup-norm of $\cC_{\rm b}(H^{2m})$.
We first consider $E_1$.
By Lemma \ref{le6.3}, we have
\begin{align}E_1\leqslant&2\|\vp\|_{\8}\bP(\xi_{r_1}(u_{10},u_{20})<t)
\leqslant2\|\vp\|_{\8}\(\bP(\xi_{r_1}(u_{10})<t)+\bP(\xi_{r_1}(u_{20})<t)\)\notag\\
\leqslant&
  \frac{2C_4}{r_1}\|\vp\|_\8\(\|u_{10}\|_{2m}^2+\|u_{10}\|^{2(2m+3)}
  +\|u_{20}\|_{2m}^2+\|u_{20}\|^{2(2m+3)}+t\)\notag\\
\leqslant&
  \frac{4C_4}{r_1}\|\vp\|_\8\(\|u_{10}\|_{2m}^2+\|u_{10}\|^{2(2m+3)}
  +T+1\).\label{6.17}
\end{align}
For $E_2$, by Lemma \ref{le6.3}, we have
\begin{align}E_2\leqslant&2\|\vp\|_{\8}\bP(\eta_{r_2}(u_{10},u_{20})<t)\leqslant
2\|\vp\|_{\8}\(\bP(\eta_{r_2}(u_{10})<t)+\bP(\eta_{r_2}(u_{20})<t)\)\notag\\
\leqslant&\frac{4C_5}{r_2}\|\vp\|_\8\(\|u_{10}\|_{2m}^2+\|u_{10}\|^{2(2m+5)}
  +T+1\).\label{6.18}
\end{align}
For each given $\ve>0$, by taking $r_1$, $r_2$ sufficiently large and fixed
in \eqref{6.17} and \eqref{6.18}, we can ensure that
\be E_1+E_2<\frac\ve2.\ee

Next we address $E_3$.
For this, we approximate $\vp$ by a Lipschitz function $\tilde\vp$ (to be chosen below).
Observe that on the set $\{\eta_{r_2}(u_{10},u_{20})\geqslant t\}$,
$\|A^{\frac{m+1}2}u(t,u_{i0})\|^2\leqslant r_2/t$.
Then it yields that
\begin{align}E_3\leqslant&
  2\sup_{u\in B}|\vp(u)-\tilde\vp(u)|
  +|\bE(\tilde\vp(u(t,u_{10}))-\tilde\vp(u(t,u_{20})))\chi_{\{\xi_{r_1}(u_{10},u_{20})
  \geqslant t\}}|\notag\\
\leqslant&
  2\sup_{u\in B}|\vp(u)-\tilde\vp(u)|\notag\\
&  +L_{\tilde\vp}\bE\|\(u(t,u_{10})-u(t,u_{20})\)\chi_{\{\xi_{r_1}(u_{10},u_{20})
  \geqslant t\}}\|_{2m},
\end{align}
where $B:=\ol{\bf B}_{2(m+1)}\((r_2/t)^{\frac12}\)$, $\ol{\bf B}_{\al}(r)$ is the closed ball in $H^{\al}$
centered at $0$ with radius $r$, and $L_{\tilde\vp}$ is the Lipschitz constant of $\tilde\vp$.
Since $H^{2(m+1)}$ is compactly embedded into $H^{2m}$, we know that $B$ is compact in $H^{2m}$.
By the density of Lipschitz functions in $\cC_{\rm b}(B)$,
we can find a Lipschitz function
$\tilde\vp\in \cC_{\rm b}(B)$
such that
\be\sup_{u\in B}|\vp(u)-\tilde\vp(u)|<\frac\ve8.\ee
The choice of $\tilde\vp$ determines $L_{\tilde\vp}$.
By Jensen's inequality and Lemma \ref{le6.4}, we have
\begin{align}&
  L_{\tilde\vp}\bE\|\(u(t,u_{10})-u(t,u_{20})\)
  \chi_{\{\xi_{r_1}(u_{10},u_{20})\geqslant t\}}\|_{2m}\notag\\
\leqslant&
  L_{\tilde\vp}\(\bE\sup_{s\in[0,t\wedge\xi_{r_1}(u_{10},u_{20})]}
  \|u(s,u_{10})-u(s,u_{20})\|^2_{2m}\)^{\frac12}
  \notag\\
\leqslant&C_6L_{\tilde\vp}\|u_{10}-u_{20}\|_{2m}.
\end{align}
Hence, for the above $\ve$, we have $r_4>0$ such that when $\|u_{10}-u_{20}\|_{2m}<r_4$,
\be\label{6.23} C_6L_{\tilde\vp}\|u_{10}-u_{20}\|_{2m}\leqslant\frac\ve4.\ee
Combining the estimations from \eqref{6.16} to \eqref{6.23}, we see that for each $\ve>0$,
whenever $\|u_{10}-u_{20}\|_{2m}<\min\{r_3,r_4\}$,
$$|\cP_t\vp(u_{10})-\cP_t\vp(u_{20})|<\ve,$$
which proves the continuity.
\eo

\subsection{Existence of ergodic invariant measures}

In this subsection, we first apply the classical procedure -- Krylov-Bogoliubov existence theorem
(\cite[Corollary 11.8]{DaPZ14}) and Prokhorov's theorem (\cite[Theorem 2.3]{DaPZ14})
to verify the the existence of the invariant measures for the semigroup $u(t)$ on $H^{2m}$.
And then we give the existence of ergodic invariant measures for $\cP_t$.
We present the classical procedure in our situation here.

\bl\label{le6.6} Assume that $\cP_t$ is Feller on $\cC_{\rm b}(H^{2m})$.
For each $u\in L^2(\W,H^{2m})$, define a family of Borel probability measures on $H^{2m}$,
\be\label{6.24}\nu_T(\cdot):=\frac1T\int_0^T\cP_t(u,\cdot)\di t,\hs T>0.\ee
Suppose that the family $\{\nu_T\}_{T>0}$ is tight.
Then each sequence $\{\nu_{T_n}\}_{n\in\N}$ with $T_n\ra\8$ as $n\ra\8$ has
a weakly convergent subsequence and the corresponding weak limit is an invariant measure
for $\cP_t$.\el

Now we show the existence of the invariant measures.

\bt\label{th6.7} Suppose that $|b|<4$ and $\phi$ satisfies \eqref{2.3} and \eqref{2.4}.
Then there exists an invariant measure on $H^{2m}$ for the semigroup $\cP_t$.\et

\bo Following Theorem \ref{th6.5} and Lemma \ref{le6.6},
we only need to show the tightness of the $\nu_T$ defined as \eqref{6.24}
for some suitable initial datum $u_0$, saying,
given some $H^{2m}$-valued $\cF_0$-measurable random variable $u_0$,
for arbitrary $\ve>0$, there exists a compact set $\cK_\ve\subset H^{2m}$ such that
\be\label{6.25}\nu_T(\cK_\ve)\geqslant 1-\ve,\hs\mb{for all }T>0,\ee

Indeed, if we set $u_0=0$, then due to the compact embedding of $H^{2(m+1)}$ into $H^{2m}$
and the fact $u(t;0)\in H^{2(m+1)}$ for all $t>0$, by Chebyshev's inequality and
\eqref{6.3} and choosing $\cK_\ve=\ol{\bf B}_{2(m+1)}(\cR)$, we have
\begin{align*}\nu_{T}(H^{2m}\sm\ol{\bf B}_{2(m+1)}(\cR))=&
  \nu_{T}(H^{2m}\sm\ol{\bf B}_{2(m+1)}(\cR))\\
=&\frac1{T}\int_0^{T}
  \bP(\|A^{m+1} u(t,0)\|>\cR)\di t\\
\leqslant&
  \frac1{T\cR^2}\bE\int_0^T\|A^{m+1} u(t,0)\|^2\di t
  \leqslant\frac{C_4}{\cR^2}<\ve,
\end{align*}
as long as $\cR$ is sufficiently large.
This asserts the existence of $\cK_\ve$ for \eqref{6.25}.
Therefore, Lemma \ref{le6.6} implies the existence of invariant Borel probability measures.
\eo

The second task now is to investigate the existence of ergodic invariant measures
for the Markovian semigroup $\cP_t$.
For this goal, it is necessary for us to show the tightness of $\cI$,
the set of all the invariant measures for $u(t)$ on $H^{2m}$.
We first prove the uniform boundedness of some integrations over $H^{2m}$
with respect to each invariant measure.

\bl\label{le6.8} Let $\nu\in\cI$ and $q\geqslant 2$.
Then we have $\tilde\rho_0^q>0$ such that
\be\label{6.26}\int_{H^{2m}}\|u\|^{q}\nu(\di u)\leqslant\tilde\rho_0^q.\ee\el
\bo We utilize the property of invariant measures to give the proof.
We again start with the equation \eqref{2.13}.
Apply It\^o's Formula to $\|u\|^q$ and then we obtain
\begin{align}&
  \di\|u(\vsig)\|^q+\|u(\vsig)\|^q\di\vsig
  +q\|u\|^{q-2}\|Au\|^2\di\vsig\notag\\
=&
  \frac{q}2\|u\|^{q-4}\left[(q-2)\(u,\phi(u)\)_{L_2(\fU,\R)}
  +\|u\|^2\|\phi(u)\|_{L_2(\fU,H)}^2\right]\di\vsig\notag\\
&
  +\|u(s)\|^q\di\vsig-q\|u\|^{q-2}\(u,f(u)\)\di\vsig
  +q\|u\|^{q-2}\(u,\phi(u)\di W\).\label{6.27}
\end{align}
Let $\xi^{0}_r(u_0)$ be a stopping time such that
$\xi^{0}_r(u_0)=\inf\{t\geqslant0:\|u(t)\|^2>r\}$.
Then multiply \eqref{6.27} by $\me^\vsig$ and integrate the obtained equality
over $\vsig\in[0,s]$ with $s\in(0,t\wedge\xi^{0}_r(u_0)]$.
Similar to the treatment of \eqref{4.38} and by noting that it is not necessary to take
the abstract value for the stochastic integral, we deduce that
\begin{align}&\|u(s)\|^{q}+\frac{q}2\(1-\frac{|b|}{4}\)
  \int_0^{s}\me^{\vsig-s}\|u\|^{q-2}\(\|Au\|^2+|u|_4^4\)\di\vsig\notag\\
\leqslant&
  \|u_0\|^q\me^{-s}+q\int_0^{s}\me^{\vsig-s}\|u\|^{q-2}\left[c\|Au\|^{\frac32}
  -\frac12\(1-\frac{|b|}{4}\)\|Au\|^2\right]\di\vsig\notag\\
&
  +q\int_0^{s}\me^{\vsig-s}\|u\|^{q-2}
  \left[c_q\(1+|u|_4^2+|u|_4^{\frac52}\)-\frac12\(1-\frac{|b|}4\)|u|_4^4\right]
  \di\vsig\notag\\
&  +q\int_0^{s}\me^{\vsig-s}\|u\|^{q-2}(u,\phi(u)\di W)\notag\\
\leqslant&
  \|u_0\|^q\me^{-s}+\tilde\rho'_0+q\int_0^{s}
  \me^{\vsig-s}\|u\|^{q-2}(u,\phi(u)\di W),
  \label{6.28}
\end{align}
where $\tilde\rho'_0=\tilde\rho'_0(q)$ is a positive constant independent of $u$ and $s$.
Let $R>0$ and consider the set $D_R=\{\w\in\W:\|u_0\|^q\leqslant R\}$.
Taking $s=t\wedge\xi^{0}_r(u_0)$, multiplying \eqref{6.28} by $\chi_{D_R}$
and taking the expected value, we have
$$\bE\|u(t\wedge\xi^{0}_r(u_0))\|^q\chi_{D_R}\leqslant
\bE\|u_0\|^q\me^{-t\wedge\xi^{0}_r(u_0)}\chi_{D_R}+\tilde{\rho}'_0\leqslant\bE\|u_0\|^q\chi_{D_R}
+\tilde{\rho}'_0,$$
by the property of stochastic integrals (\cite{DaPZ14}), which annihilates the stochastic integral.
Fix $R>0$. Applying the Dominated Convergence Theorem, we let $r$ tend to the infinity, and obtain
$$\bE\|u(t)\|^q\chi_{D_R}\leqslant\me^{-t}\bE\|u_0\|^q\chi_{D_R}+\tilde{\rho}'_0.$$
For each fixed $R>0$, we can find $t_R>0$ such that
$$\bE\|u(t_R)\|^q\chi_{D_R}\leqslant2\tilde{\rho}'_0.$$
By the property of invariant measures, we know that the distribution of $\|u(t_R)\|^{q}$ is $\nu$.
Hence
$$
\bE\|u(t_R)\|^q\chi_{D_R}=\int_{H^{2m}}\|u\|^q\chi_{u(t_R;D_R)}\nu(\di u)\ra\int_{H^{2m}}\|u\|^q\nu(\di u).
$$
as $R\ra+\8$.
Note that $\tilde\rho'_0$ is independent of $R$.
By taking $\tilde{\rho}_0^q=2\tilde{\rho}'_0$ and setting $R\ra+\8$,
we obtain \eqref{6.26} and complete the proof.
\eo

\bl\label{le6.9} Let $\nu\in\cI$.
Then for the $m$ given in this section, we have $\tilde\rho_m>0$ such that
\be\label{6.29}\int_{H^{2m}}\|u\|_{2m}^2\nu(\di u)\leqslant\tilde\rho_m.\ee
\el

\bo First we follow \eqref{6.28}, and deduce a more concise form
\begin{align}&\|u(s)\|^q+\int_0^s\me^{\vsig-s}\|u\|^{q-2}\|Au\|^2\di\vsig\notag\\
\leqslant&\|u_0\|^q\me^{-s}+\tilde\rho'_{0,q}
+c_q\int_0^{s}\me^{\vsig-s}\|u\|^{q-2}(u,\phi(u)\di W),\label{6.30}\end{align}
for $s\in(0,t\wedge\xi_r^0(u_0)]$, where $\tilde\rho'_{0,q}$ is independent of $u$ and $s$.
For $r>0$ and each $0\leqslant k\leqslant m$, we define
$$\tilde\xi_r^k(u_0):
=\inf\left\{t\geqslant0:\inf_{0\leqslant k'\leqslant k}
\|A^{k'}u(t)\|^2>r\right\}.$$
After \eqref{6.30} we claim that for each $s\in(0,t\wedge\tilde\xi_r^k(u_0)]$ with
$t>0$,
\begin{align}&\|A^ku(s)\|^2+\int_0^s\me^{\vsig-s}\|A^{k+1}u\|^2\di\vsig\notag\\
\leqslant& c_k\(\|A^ku_0\|^2+Q_k(u_0)s\)\me^{-s}+\tilde\rho'_{k}+\Theta_k(s),
\label{6.31}\end{align}
for each $0\leqslant k\leqslant m$ and some deterministic $\tilde\rho'_k$
independent of $u$ and $s$, where $Q_0(s,u_0)\equiv0$ and for $k\geqslant1$,
\be\label{6.31AA}Q_k(u_0):=\sum_{l=1}^k\|u_0\|^{2(2l+3)},\ee
$\Theta_k(s)$ is the sum of the terms for
the stochastic integral with respect to $W$ and hence
$\bE\Theta_k(t\wedge\tau)=0$, for each stopping time $\tau\geqslant0$.

Now we show the claim by the induction method.
Firstly, the inequality \eqref{6.30} has implied the claim
for the case when $k=0$.
Secondly, we show that when the claim holds for $k=l-1$ with $l\in\{1,2,\cdots,m\}$,
then it also holds for $k=l$.

Still apply $A^l$ to \eqref{2.13}, use It\^o's Formula to $\|A^lu\|^2$ and
add $\|A^lu(\vsig)\|^2\di\vsig$ to both sides of the obtained equality.
We have
\begin{align}
  &\di\|A^lu(\vsig)\|^2+\|A^lu(\vsig)\|^2\di\vsig+2\|A^{l+1}u\|^2\di\vsig\notag\\
=&[\|A^l\phi(u)\|_{L_2(\fU,H)}^2-2(A^lu,A^lf(u))+\|A^lu\|^2]\di\vsig\notag\\
&  +2\(A^lu,A^l\phi(u)\di W\).\label{6.32}
\end{align}
According to the treatment of \eqref{3.7}, we notice that
$$|\|A^l\phi(u)\|_{L_2(\fU,H)}^2-2(A^lu,A^lf(u))|\!\leqslant\!
\|A^{l+1}u\|^2+c_l\(\|A^lu\|^2+\|u\|^{2(2l+3)}+1\).$$
We multiply \eqref{6.32} with $\me^{\vsig}$, integrate \eqref{6.29} over
$\vsig\in[0,s]$ with $s\in(0,t\wedge\xi_r^{l}(u_0)]$, and obtain
\begin{align*}&\|A^lu(s)\|^2+\int_0^s\me^{\vsig-s}\|A^{l+1}u\|^2\di\vsig\\
\leqslant&
  \|A^lu_0\|^2\me^{-s}+c_l\int_0^{s}\me^{\vsig-s}
  \(\|A^lu\|^2+\|u\|^{2(2l+3)}+1\)\di\vsig\\
&
  +2\int_0^{s}\me^{\vsig-s}\(A^lu,A^l\phi(u)\di W\).
\end{align*}
By \eqref{6.31} with $k=l-1$ and \eqref{6.30}, it follows by embeddings that
\begin{align}&\|A^lu(s)\|^2+\int_0^s\me^{\vsig-s}\|A^{l+1}u\|^2\di\vsig\notag\\
\leqslant&
  \left[\|A^lu_0\|^2+c_l\(\|A^{l-1}u_0\|^2+Q_{l-1}(u_0)s\)\right]\me^{-s}
  +c_l\tilde\rho'_{l-1}+c_l\Theta_{l-1}(s)\notag\\
&
  +c_l\int_0^{s}\me^{\vsig-s}
  \(\|u_0\|^{2(2l+3)}\me^{-\vsig}+\(\tilde\rho'_{0,2(2l+3)}+1\)\)\di\vsig\notag\\
&
  +c_l\int_0^s\me^{\vsig-s}\int_0^{\vsig}\me^{\vsig'-\vsig}\|u\|^{2(2l+2)}(u,\phi(u)\di W(\vsig'))\di\vsig\notag\\
&
  +2\int_0^{s}\me^{\vsig-s}\(A^lu,A^l\phi(u)\di W(\vsig)\)\notag\\
\leqslant&
  c_l\(\|A^lu_0\|^2+Q_l(u_0)s\)\me^{-s}
  +c_l\(\tilde\rho'_{l-1}+\tilde\rho'_{0,2(2l+3)}+1\)\notag\\
&
  +c_l\Theta_{l-1}(s)+c_l\int_0^s\int_0^{\vsig}\me^{\vsig'-s}
  \|u\|^{2(2l+2)}\(u,\phi(u)\di W(\vsig')\)\di\vsig\notag\\
&
  +2\int_0^{s}\me^{\vsig-s}\(A^lu,A^l\phi(u)\di W(\vsig)\).\label{6.33}
\end{align}
By setting
$$\tilde\rho'_l=c_l\(\tilde\rho'_{l-1}+\tilde\rho'_{0,2(2l+3)}+1\)\Hs\mb{and}$$
\begin{align*}\Theta_l(s)=&c_l\Theta_{l-1}(s)+c_l\int_0^s\int_0^{\vsig}\me^{\vsig'-s}
  \|u\|^{2(2l+2)}\(u,\phi(u)\di W(\vsig')\)\di\vsig\\
&  +2\int_0^{s}\me^{\vsig-s}\(A^lu,A^l\phi(u)\di W\),\end{align*}
we see by Fubini theorem that for each stopping time $\tau\geqslant0$,
$$\bE\Theta_l(t\wedge\tau)=c_l
\bE\int_0^{t\wedge\tau}(t\wedge\tau-\vsig')\me^{\vsig'-t\wedge\tau}
\|u\|^{2(2l+2)}\(u,\phi(u)\di W(\vsig')\)=0,
$$
and \eqref{6.33} indicates the claim for $k=l$.

By the claim we have
\begin{align}&\|A^mu(t\wedge\tilde\xi_r^m(u_0))\|^2+\int_0^{t\wedge\tilde\xi_r^m(u_0)}
\me^{\vsig-s}\|A^{m+1}u\|^2\di\vsig\notag\\
\leqslant&c_m\(\|A^mu_0\|^2+Q_m(u_0)t\wedge\tilde\xi_r^m(u_0)\)
\me^{-t\wedge\tilde\xi_r^m(u_0)}\notag\\
&
 +\tilde\rho'_{m}+\Theta_m(t\wedge\tilde\xi_r^m(u_0)).\label{6.34}\end{align}
Now let $R>0$ and define
$$\cD_R:=\{\w\in\W:\|A^mu_0\|^2+Q_m(u_0)\leqslant R\}.$$
Similarly multiply \eqref{6.34} with $\chi_{\cD_R}$ and take the expected value of
the obtained inequality.
One sees that
\begin{align}&\bE\|A^mu(t\wedge\tilde\xi^m_r(u_0))\|^2\chi_{\cD_R}\notag\\
\leqslant&
  c_m\bE\(\|A^mu_0\|^2+Q_m(u_0)t\wedge\tilde\xi_r^m(u_0)\)\me^{-t\wedge\tilde\xi_r^m(u_0)}
  \chi_{\cD_{R}}+\tilde\rho'_m\notag\\
\leqslant&
  c_m\bE\(\|A^mu_0\|^2+Q_m(u_0)\)\chi_{\cD_{R}}+\tilde\rho'_m.\label{6.35}
\end{align}
Fix $R>0$, use the Dominated Convergence Theorem, and let $r$ tend to the infinity.
We obtain that
$$\bE\|A^mu(t)\|^2\chi_{\cD_R}\leqslant
c_m\bE[\|A^mu_0\|^2+Q_m(u_0)t]\chi_{\cD_{R}}\me^{-t}+\tilde\rho'_m$$
For each fixed $R>0$, we have a $t'_R>0$ such that
\be\label{6.36}\bE\|A^mu(t)\|^2\chi_{\cD_R}\leqslant2\tilde\rho'_m.\ee
Again, using the property of invariant measures, we infer that as $R\ra+\8$,
$$\bE\|A^mu(t)\|^2\chi_{\cD_R}=\int_{H^{2m}}\|A^mu\|^2\chi_{u(t'_R;\cD_R)}\nu(\di u)
\ra\int_{H^{2m}}\|A^{2m}u\|^2\nu(\di u).$$
Since $\tilde\rho'_m$ is independent of $R$, we set $\tilde\rho_m=2\tilde\rho'_m$ and $R\ra+\8$
in \eqref{6.36} and get \eqref{6.29}.
\eo

We now can present the tightness of $\cI$ as follows.
\bl\label{le4.12} The invariant measure set $\cI$ is tight.\el

\bo According to the definition, for every $\ve>0$, we look for a compact subset
$\ol{\bf B}_{2(m+1)}(\cR)$ with a sufficiently large $\cR$ such that
\be\label{6.37}\nu(H^{2m}\sm\ol{\bf B}_{2(m+1)}(\cR))<\ve\hs\mb{for all }\nu\in\cI.\ee
Indeed, for some fixed $t>0$ and by the invariance and Lemma \ref{le6.1},
we have
\begin{align}&\nu(H^{2m}\sm\ol{\bf B}_{2(m+1)}(\cR))
  =\int_{H^{2m}}\cP_t(\tilde u,H^{2m}\sm\ol{\bf B}_{2(m+1)}(\cR))\nu(\di\tilde u)\notag\\
=&
  \int_{H^{2m}}\bP(t\|A^{m+1}u(t,\tilde u)\|^2>t\cR^2)\nu(\di\tilde u)
  \leqslant\int_{H^{2m}}\bP\(\eta_{t\cR^2}(\tilde u)<t\)\nu(\di\tilde u)\notag\\
\leqslant&
  \frac{C_5}{t\cR^2}\int_{H^{2m}}\(\|\tilde u\|_{2m}^2+\|\tilde u\|^{2(2m+5)}+t\)
  \nu(\di\tilde u),
\label{6.38}
\end{align}
where $C_5=C_5(t)$ is given by \eqref{6.6}.
By Lemmas \ref{le6.8} and \ref{le6.9}, we know the integral with respect to $\nu$ over $H^{2m}$
in \eqref{6.38} is bounded only by a deterministic positive constant depending on $t$, i.e.,
$$\nu(H^{2m}\sm\ol{\bf B}_{2(m+1)}(\cR))\leqslant\frac{C'_{12}(t)}{t\cR^2}.$$
Thus, for each $\ve>0$, if we choose $\cR$ sufficiently large, we are able to ensure \eqref{6.37}
and complete the proof.
\eo

Now we are prepared to verify the existence of ergodic invariant measures.
\bt Suppose that $|b|<4$ and $\phi$ satisfies \eqref{2.3} and \eqref{2.4}.
Then there exists an ergodic invariant measure on $H^{2m}$
for the Feller Markovian semigroup $\cP_t$.
\et
\bo
We adopt a similar routine to \cite[Theorem IV.13]{WLYJ21} to address this proof.
Since an invariant measure is ergodic if and only if it is an extreme point of $\cI$
(see \cite[Proposition 11.12]{DaPZ14}),
we only need to show the existence of extreme points of $\cI$,
for which, it suffices to show that $\cI$ is nonempty, compact and convex,
by the Krein-Milman Theorem (see \cite[Theorem 7.4, Chapter V]{Conway90}).

The existence of invariant measures implies that $\cI\ne\es$.
The linearity of $\cP^*_t$ with respect to $\nu$ implies that $\cI$ is convex.
For the compactness, we first show that $\cI$ is closed.

Let $\{\nu_n\}_{n\in\N}$ be a given weakly convergent sequence in $\cI$
such that $\nu_n$ weakly converges to $\nu$, as $n\ra\8$.
By the invariance, we know that for all $\Gam\in\cB(H^{2m})$,
when $n$ tends to the infinity,
\begin{align*}
&\nu(\Gam)\la\nu_n(\Gam)=\cP^*_t\nu_n(\Gam)\\
=&\int_H\cP_t(u,\Gam)\nu_n(\di u)
\ra\int_H\cP_t(u,\Gam)\nu(\di u)=\cP^*_t\nu(\Gam).
\end{align*}
By the uniqueness of weak limit, we know that $\nu$ is also an invariant measure.
Hence the closedness.

The closedness and Lemma \ref{le4.12} guarantee the compactness of $\cI$,
which eventually proves this theorem.
\eo
\section{Infinite Regularity of Invariant Measures}\label{s7}
In this section we are to give the infinite regularity of the invariant measures obtained
in Section \ref{s6}.
Let $m\geqslant1$ and $\nu$ be the invariant measure for the Markovian semigroup $\cP_t$ on $H^{2m}$
given in Theorem \ref{th6.7}.
First we increase the regularity of $\nu$ to the smoother space $H^{2(m+1)}$.

For this aim, recall the discussion in Section \ref{s6}.
We note that if it could hold that (stated in Theorem \ref{th7.1} below)
\be\label{7.1}\int_{H^{2m}}\|u\|_{2(m+1)}^2\nu(\di u)<+\8,\ee
we can conclude \textit{a posteriori} that the support of $\nu$ is contained in $H^{2(m+1)}$
(using the ideas in \cite{BKR96,GKV14}).
Unfortunately, this can not come true in the way presented
in the proof of Lemma \ref{le6.9},
since $A^{m+1}u(t)$ makes no sense at $t=0$.
However, we can make it by reusing the method given in the proof of Lemma \ref{le6.1}.

\bt\label{th7.1} Let $|b|<4$ and $\phi$ satisfy \eqref{2.3} and \eqref{2.4}.
Then \eqref{7.1} holds true.
Moreover, the support of $\nu$ is contained in $H^{2(m+1)}$.
\et
\bo We only show \eqref{7.1}.
Similar to the estimation in the proof of Lemma \ref{le6.1},
we apply It\^o's Formula to $\|A^{m+1}u\|^2$ and deduce like
the argument in the proof of Lemma \ref{le6.8} that
for stopping times $0\leqslant\vsig\leqslant s\leqslant t\wedge\tilde\xi_r^{m+1}(u_0)$,
with $\tilde\xi_r^{m+1}(u_0)$ defined for each $r>0$ as
$$\tilde\xi_r^{m+1}(u_0):
=\inf\left\{t\geqslant0:\inf_{0\leqslant k'\leqslant m+1}
\|A^{k'}u(t)\|^2>r\right\}.$$
it holds that
\begin{align}&\|A^{m+1}u(s)\|^2\notag\\
\leqslant&
  \|A^{m+1}u(\vsig)\|^2\me^{\vsig-s}+c\int_{\vsig}^s\me^{s'-s}
  \(\|A^{m+1}u\|^2+\|u\|^{2(2m+5)}+1\)\di s'\notag\\
&
  +2\int_{\vsig}^s\me^{s'-s}(A^{m+1}u,A^{m+1}\phi(u)\di W).
\label{7.2}\end{align}
Integrating \eqref{7.2} over $\vsig\in[0,s]$ and recalling \eqref{6.30}
and the claim \eqref{6.31}, we infer that
\begin{align}&s\|A^{m+1}u(s)\|^2\notag\\
\leqslant&
  \int_0^s\|A^{m+1}u(\vsig)\|^2\me^{\vsig-s}\di\vsig
  +cs\int_0^s\me^{s'-s}\(\|A^{m+1}u\|^2+\|u\|^{2(2m+5)}+1\)\di s'\notag\\
&
  +\int_0^s\int_{\vsig}^s\me^{s'-s}(A^{m+1}u,A^{m+1}\phi(u)\di W(s'))
  \di\vsig\notag\\
\leqslant&
  (1+cs)\(c_m(\|A^mu_0\|^2+Q_m(u_0)s)\me^{-s}+\tilde\rho'_m+\Theta_m(s)\)
  \notag\\
&
  +cs\int_0^s\me^{s'-s}\(\|u_0\|^{2(2m+5)}\me^{-s'}+\tilde\rho'_{0,2(2m+5)}+1\)\di s'\notag\\
&
  +c_ms\int_0^s\int_0^{s'}\me^{\vsig'-s}\|u\|^{2(2m+4)}(u,\phi(u)\di W(\vsig'))
  \di s'\notag\\
&
  +2\int_0^ss'\me^{s'-s}(A^{m+1}u,A^{m+1}\phi(u)\di W(s'))\notag\\
\leqslant&
  c_m(1+s)\(\|A^{m}u_0\|^2+Q_{m+1}(u_0)s\)\me^{-s}+c(1+s)\tilde\rho'_{m+1}
  +\Theta_{m+1}(s),
\label{7.3}\end{align}
where $Q_{m+1}(u_0)$ is given in \eqref{6.31AA},
$$\tilde\rho'_{m+1}=\tilde\rho'_m+\tilde\rho'_{0,2(2m+5)}+1\hs\mb{and}$$
\begin{align*}\Theta_{m+1}(s)=&
  \(1+cs\)\Theta_m(s)+c_ms\int_0^s(s-\vsig')\me^{\vsig'-s}\|u\|^{2(2m+4)}(u,\phi(u)
  \di W(\vsig'))\\
&
  +2\int_0^ss'\me^{s'-s}(A^{m+1}u,A^{m+1}\phi(u)\di W(s')).
\end{align*}
Similar to the proof of Lemmas \ref{le6.8} and \ref{le6.9}, we again define
$$\cD_\cR=\{\w\in\W:\|A^mu_0\|^2+\|u_0\|^{2(2m+5)}\leqslant\cR\}.$$
Picking $s=t\wedge\tilde\xi_r^{m+1}(u_0)$ in \eqref{7.3},
multiplying the obtained inequality by $\chi_{\cD_\cR}$
and taking the expectation of what is just obtained,
we annihilate the stochastic term $\Theta_{m+1}(t\wedge\tilde\xi_r^{m+1}(u_0))$ and easily obtain
\begin{align}&\bE\(t\wedge\tilde\xi_r^{m+1}(u_0)\)\|A^{m+1}u(t\wedge\tilde\xi_r^{m+1}(u_0))\|^2
  \chi_{\cD_\cR}\notag\\
\leqslant&
  c_m\bE(\|A^{m}u_0\|^2+Q_{m+1}(u_0))+c(1+t)\tilde\rho'_{m+1},\label{7.4}
\end{align}
by the boundedness of $s^k\me^{-s}$ with $s\geqslant0$ and $k\in\N^+$.
Fix $\cR>0$. Using the Dominated Convergence Theorem to \eqref{7.4} and considering $r\ra+\8$,
we can obtain that
\begin{align*}&\bE\|A^{m+1}u(t)\|^2\chi_{\cD_\cR}\\\leqslant&
c_m\(1+\frac1t\)\me^{-t}\bE\(\|A^{m}u_0\|^2+Q_{m+1}(u_0)t\)\chi_{\cD_\cR}+c\(1+\frac1t\)\tilde\rho'_{m+1}.
\end{align*}
Given each fixed $\cR>0$, it is easy to know that there exists $t_{\cR}>0$ such that
\be\bE\|A^{m+1}u(t)\|^2\chi_{\cD_\cR}\leqslant2c\tilde\rho'_{m+1}.\label{7.5}\ee
Then as $\cR\ra+\8$, we obtain
\begin{align*}\bE\|A^{m+1}u(t)\|^2\chi_{\cD_\cR}=&
\int_{H^{2m}}\|A^{m+1}u\|^2\chi_{u(t_\cR;\cD_\cR)}\nu(\di u)\\
\ra&\int_{H^{2m}}\|A^{m+1}u\|^2\nu(\di u).\end{align*}
Owing to the independence of $\cR$ for $\tilde\rho'_{m+1}$, we immediately obtain \eqref{7.1}
by letting $\cR\ra+\8$ in \eqref{7.5}. The proof is complete.
\eo

Based on the regularity of $\nu$ obtained in Theorem \ref{th7.1}, the infinite regularity of $\nu$
is a trivial consequence by induction and embeddings.
\bt Let $|b|<4$.
Suppose that $\phi$ satisfies \eqref{2.5}.
Then each invariant measure supported on $H^2$ for the Markovian semigroup $\cP_t$
is supported on $\cC^{\8}\cap H$.
\et
\section{Summary and Remarks}\label{s8}

In this paper, we aim to obtain the existence and regularity of invariant measures for the stochastic
MSHE with multiplicative noise and periodic boundary.
To this end, we first confirm the unique existence of global pathwise solutions,
and by Yamada-Wannabe theorem, we are also required to provide the existence of martingale solutions
in advance.
With a further study, we also obtain the global existence of martingale solutions.

In our situation, the concept of global martingale solution is stronger a bit than that in a previous
work \cite{MSY16}, and we use a new method to obtain this kind of global martingale solutions.
Moreover, our new method can be applied extensively to other stochastic nonlinear partial differential
equations to get a global martingale solution, once the drift is dissipative somehow and the diffusion
coefficient is globally Lipschitz, or more weakly, the solutions coincide locally in time for Galerkin approximating
cut-off systems with different cut-off functions.
However, we are still unaware whether the global martingale solution is unique.

In the proof of existence of ergodic invariant measures in $H^{2m}$, we actually follow the estimation
methods used frequently in attractor theory to study the existence of the global attractor in spaces of
high regularity (see \cite{SZM10,WWL22,WZZ21}).
Correspondingly, to improve the regularity of the invariant measures, we also use the
regularity estimation methods for deterministic systems.
All in all, we eventually obtain the infinite regularity of invariant measures for stochastic MSHE with
multiplicative noise and periodic boundary under appropriate conditions on the diffusion coefficient $\phi$,
which is just what Glatt-Holtz expected in \cite{GKV14} for the stochastic primitive equation.

Combining the research on invariant measures and attractors
in Chen et al. \cite{CDJNZ,CDJZ} and Huang et al. \cite{HJLY1,HJLY2,HJLY3,HJLY4,HJLY5,HJLY6},
we have some ideas about the relation between invariant measures for stochastic systems
and attractors for deterministic systems as follows.

Consider the following abstract deterministic differential equation
\be\label{8.1}u_t=Au+f(u),\Hs t>0,\ee
and the corresponding stochastic differential equation with multiplicative noise
\be\label{8.2}\di u=[Au+f(u)]\di t+\phi(u)\di W,\Hs t>0,\ee
where $u:[0,+\8)\ra X$ is the unknown function, $X$ is a Hilbert (or Banach) space,
$A:\cD(A)\subset X\ra X$ is linear operator, $f$ is the external force, $\phi(u)\di W$ is
the stochastic term and $W$ is a Wiener process.

Firstly, we conjecture that, \textit{if \eqref{8.1} possesses a unique global strong solution
in $X$ for each initial datum, then \eqref{8.2} possesses a unique global pathwise solution
for a given stochastic basis under some appropriate global Lipschitz condition on $\phi$.}
Secondly, we conjecture that, \textit{if the dynamical system generated by \eqref{8.1} has
a global attractor $\sA$ in $X$, then the Markovian transition semigroup associated with
\eqref{8.2} has an invariant measure $\mu$ on $X$ under some appropriate global Lipschitz
condition on $\phi$}.
Moreover, \textit{if $\sA$ is contained in a subspace $Y$ with higher regularity of $X$,
then $\mu$ is supported on $Y$ as well, under  some appropriate global Lipschitz condition
on $\phi$}.
Thirdly, if the conjectures mentioned above hold true, perhaps the converse of
these conclusions are also hopeful to be valid under some appropriate conditions.

The proof of these conjectures may need very elaborate investigations, but at least,
there should be many examples indicating these conclusions.
What is more, we are also interested in the existence of (ergodic) invariant measures
when the drift is not dissipative, for example, when $|b|\geqslant 4$ for the stochastic
MSHE \eqref{1.1}.

\section*{Data Availability}

Data sharing is not applicable to this article as no new data were created or analyzed in this study.


\begin{thebibliography}{99}
%
%

\bibitem{B95}Bensoussan, A.:
 Stochastic Navier-Stokes equations.
 Acta Appl. Math. {\bf 38}(3), 267-304 (1995)

\bibitem{BKR96}Bogachev, V.I., Krylov, N., R\"ockner, M.:
 Regularity of invariant measures: the case of non-constant diffusion part.
 J. Funct. Anal. {\bf 138}, 223-242 (1996)

\bibitem{CDJNZ}Chen, L.F., Dong, Z., Jiang, J.F., Niu, L., Zhai, J.L.:
 Decomposition formula and stationary measures for stochastic Lotka-Volterra system
 with applications to turbulent convection.
 J. Math. Pures Appl. {\bf 125}, 43-93 (2019)

\bibitem{CDJZ}Chen, L.F., Dong, Z., Jiang, J.F., Zhai, J.L.:
 On limiting behavior of stationary measures for stochastic evolution systems
 with small noise intensity.
 Sci. China Math. {\bf 63}, 1463-1504 (2020)

\bibitem{C15}Choi, Y.:
 Dynamical bifurcation of the one dimensional modified Swift-Hohenberg equation.
 Bull. Korean. Math. Soc. {\bf 52}, 1241-1252 (2015)

\bibitem{CHHL17}Choi, Y., Ha, T., Han, J., Lee, D.S.:
 Bifurcation and final patters of a modified Swift-Hohenberg equation.
 Discrete Contin. Dyn. Syst., Ser. B {\bf 22}, 2543-2567 (2017)

\bibitem{Conway90}Conway, J.B.:
 A Course in Functional Analysis (Second Edition).
 Springer-Verlag, New York (1990)

\bibitem{DaPZ14}Da Prato, G., Zabczyk, J.:
 Stochastic Equations in Infinite Dimensions (Second Edition).
 Cambridge University Press, Cambridge (2014)

\bibitem{DGHT11}Debussche, A., Glatt-Holtz, N., Temam, R.:
 Local martingale and pathwise solutions for an abstract fluids model.
 Phys. D {\bf 240}, 1123-1144, (2011)

\bibitem{DJZ19}Dhariwal, G., J\"ungel, A., Zamponi N.,
 Global martingale solutions for a stochastic population cross-diffusion system.
 Stochastic Process. Appl. {\bf 129}, 3792-3820 (2019)

\bibitem{DZ20}Du, L.H., Zhang, T.:
 Local and global existence of pathwise solution for
 the stochastic Boussinesq equations with multiplicative noises,
 Stochastic Process. Appl. {\bf 130}, 1545-1567 (2020)


\bibitem{DG12}Duan, N., Gao, W.J.:
 Optimal control of a modified Swift-Hohenberg equation.
 Electron. J. Differ. Equations {\bf 2012}, 155 (2012)

\bibitem{EH01}Eckmann, J.P., Hairer, M.:
 Invariant measures for stochastic partial differential equations
 in unbounded domains.
 Nonlinearity {\bf 14}, 133-151 (2001)

\bibitem{EKZ17}Ekren, I., Kukavica, I., Ziane, M.:
 Existence of invariant measures for the stochastic damped Schr\"odinger equation.
 Stoch. PDE: Anal. Comp. {\bf 5}, 343-367 (2017)

\bibitem{E13}Evans, L.C.,
 An Introduction to Stochastic Differential Equations (Second Edition).
 American Mathematical Society, USA (2013)

\bibitem{FG95}Flandoli, F., Gatarek, D.:
 Martingale and stationary solutions for stochastic Navier-Stokes equations.
 Probab. Theory Related Fields {\bf 102}(3), 367-391 (1995)

\bibitem{FMRT01}Foias, C., Manley, O., Rosa, R., Temam, R.:
 Navier-Stokes Equations and Turbulence.
 Cambridge University Press, Cambrige (2001)

\bibitem{F48}Folland, G.B.:
 Real Analysis (second edition), in: Pure and Applied Mathematics.
 John Wiley \& Sons Inc, New York (1999)

\bibitem{GKV14}Glatt-Holtz, N., Kukavica, I., Vicol, V., Ziane, M.:
 Existence and regularity of invariant measures for
 the three dimensional stochastic primitive equations.
 J. Math. Phys. {\bf 55}, 277-304 (2014)

\bibitem{GHV14}Glatt-Holtz, N., Vicol, V.C.:
 Local and global existence of smooth solutions for
 the stochastic Euler equations with multiplicative noise.
 Ann. Probab. {\bf 42}(1), 80-145 (2014)

\bibitem{GHZ09}Glatt-Holtz, N., Ziane, M.:
 Strong pathwise solutions of the stochastic Navier-Stokes system.
 Adv. Differ. Equations {\bf 14}, 567-600 (2009)

\bibitem{GGL17}Guo, C.X., Guo, Y.F., Li, C.Y.:
 Dynamic behaviors of a local modified stochastic Swift-Hohenberg equation
 with multiplicative noise.
 Bound. Value Probl. {\bf 2017}, 19-1-13 (2017)

\bibitem{GK96}Gy\"ongy, I., Krylov, N.:
 Existence of strong solutions for It\^o's stochastic equations via approximations.
 Probab. Theory Related Fields, {\bf 105}(2), 143-158 (1996)

\bibitem{Hen81}Henry, D.:
 Geometric theory of semilinear parabolic equations.
 Lect. Notes in Math. 840, Springer-Verlag, Berlin, New York (1981)

\bibitem{HS92}Hohenberg, P.C., Swift, J.B.:
 Effects of additive noise at the onset of Rayleigh-B\'enard's convection.
 Phys. Rev. A {\bf 46}(8), 4773-4785 (1992)

\bibitem{HJLY1}Huang, W., Ji, M., Liu, Z.X., Yi, Y.F.:
 Integral identity and measure estimates for stationary Fokker-Planck equations.
 Ann. Probab. {\bf 43}, 1712-1730 (2015)

\bibitem{HJLY2}Huang, W., Ji, M., Liu, Z.X., Yi, Y.F.:
 Steady states of Fokker-Planck equations: I. Existence.
 J. Dynam. Differential Equations {\bf 27}, 721-742 (2015)

\bibitem{HJLY3}Huang, W., Ji, M., Liu, Z.X., Yi, Y.F.:
 Steady states of Fokker-Planck equations: II. Non-existence.
 J. Dynam. Differential Equations {\bf 27}, 743-762 (2015)

\bibitem{HJLY4}Huang, W., Ji, M., Liu, Z.X., Yi, Y.F.:
 Steady states of Fokker-Planck equations: III. Degenerate diffusion.
 J. Dynam. Differential Equations {\bf 28}, 127-141 (2016)

\bibitem{HJLY5}Huang, W., Ji, M., Liu, Z.X., Yi, Y.F.:
 Stochastic stability of measures in gradient systems.
 Phys. D {\bf 314}, 9-17 (2016)

\bibitem{HJLY6}Huang, W., Ji, M., Liu, Z.X., Yi, Y.F.:
 Concentration and limit behaviors of stationary measures.
 Phys. D {\bf 369}, 1-17 (2018)

\bibitem{LMR75}LaQuey, R.E., Mahajan, S.M., Rutherford, P.H., et al.:
 Nonlinear saturation of the trapped-ion mode.
 Phys. Rev. Lett. {\bf 34}, 391-394 (1975)

\bibitem{LMN94}Lega, J., Moloney, J.V., Newell, A.C.:
 Swift-Hohenberg equation for lasers.
 Phys. Rev. Lett. {\bf 73}, 2978-2981 (1994)

\bibitem{LLW21}Li, C.Q., Li, D.S., Wang, J.T.:
 A remark on attractor bifurcation.
 Dyn. Partial Differ. Equ. {\bf 18}(2), 157-172 (2021)

\bibitem{LW18}Li, D.S., Wang, Z.-Q.:
 Local and global dynamic bifurcations of nonlinear evolution equations.
 Indiana U. Math. J. {\bf 67}, 583-621 (2018)

\bibitem{LWZ20} Li, Y.J., Wu, H.Q., Zhao, T.G.:
 Random pullback attractor of a non-autonomous local modified stochastic
 Swift-Hohenberg equation with multiplicative noise.
 J. Math. Phys. {\bf 61}, 1-13 (2020)

\bibitem{MSY16}Misiats, O., Stanzhytskyi, O., Yip, N.K.:
 Existence and uniqueness of invariant measures for stochastic reaction-diffusion equations
 in unbounded domains.
 J. Theor. Probab {\bf 29}, 996-1026 (2016)

\bibitem{PM80}Pomeau, Y., Manneville, P.:
 Wave length seletion in cellular flows.
 Phys. Lett. A {\bf 75}, 296-298 (1980)

\bibitem{SY02}Sell, G.R., You, Y.:
 Dynamics of Evolutionary Equations.
 Springer, New York (2002)

\bibitem{S77}Sivashinsky, G.L.:
 Nonlinear analysis of hydrodynamic instability in laminar flames
 i. derivation of basic equations.
 Acta Astron. {\bf 4}, 11-12 (1977)

\bibitem{SZM10}Song, L.Y., Zhang, Y.D., Ma, T.:
 Global attractor of a modified Swift-Hohenberg equation in $H^k$ spaces.
 Nonlinear Anal. {\bf 72}, 183-191 (2010)

\bibitem{S18}Sun, B.:
 Optimal disturbed contral problem for the modified Swift-Hohenberg equations.
 Electron. J. Differ. Equations {\bf 2018}, 131 (2018)

\bibitem{SH77}Swift, J.B., Hohenberg, P.C.:
 Hydrodynamic fluctuations at the convective instability.
 Phys. Rev. A {\bf 15}, 319-328 (1977)

\bibitem{Tem01}Temam, R.:
 Navier-Stokes equations: theory and numerical analysis.
 AMS Chelsea Publishing, Providence, RI, Reprint of the 1984 edition (2001)

\bibitem{WLYJ21}Wang, J.T., Li, C.Q., Yang, L., Jia, M.:
 Upper semi-continuity of random attractors and existence of invariant measures
 for nonlocal stochastic Swift-Hohenberg equation with multiplicative noise.
 J. Math. Phys., {\bf 62}, 111507 (2021)

\bibitem{WL21}Wang, J.T., Li, D.S.:
 On relative category and Morse decompositions for infinite-dimensional dynamical systems.
 Topology Appl. {\bf 291}, 107624-1-16 (2021)

\bibitem{WLD16}Wang, J.T., Li, D.S., Duan, J.Q.:
 On the shape Conley index theory of semiflows on complete metric spaces.
 Discrete Contin. Dyn. Syst. {\bf 36}(3), 1629-1647 (2016)

\bibitem{WLD21}Wang, J.T., Li, D.S., Duan, J.Q.:
 Compactly generated shape index theory and its application to
 a retarded nonautonomous parabolic equation.
 Topol. Methods Nonlinear Anal., {\bf 59}(1), 1-33 (2022).

\bibitem{WYD20}Wang, J.T., Yang, L., Duan, J.Q.:
 Recurrent solutions of a nonautonomous modified Swift-Hohenberg equation.
 Appl. Math. Comput. {\bf 379}, 125270 (2020)

\bibitem{WZZ21}Wang, J.T., Zhang, X.Q., Zhao, C.D.:
 Statistical solutions for a nonautonomous modified Swift-Hohenberg equation.
 Math. Methods Appl. Sci. {\bf 44}, 14502-14516 (2021)

\bibitem{WZ23}Wang, J.T., Zhang, X.Q.:
  Invariant sample measures and random Liouville type theorem for a nonautonomous stochastic $p$-Laplacian equation,
  Discrete Contin. Dyn. Syst. Ser. B {\bf 28}(4), 2803-2827 (2023).

\bibitem{WZC20}Wang, J.T., Zhao, C.D., Caraballo, T.:
 Invariant measures for the 3D globally modified Navier-Stokes equations
 with unbounded variable delays.
 Commun. Nonlinear Sci. Numer. Simul. {\bf 91}, 105459 (2020)

\bibitem{WWL22}Wang, X.J., Wang, J.T., Li, C.Q.:
 Invariant measures and statistical solutions for a nonautonomous
 nonlocal Swift-Hohenberg equation.
 Dyn. Syst., {\bf 37}(1), 136-158 (2022).

\bibitem{WSD05}Wang, W., Sun, J.-H., Duan, J.Q.:
 Ergodic dynamics of the stochastic Swift-Hohenberg system.
 Nonlinear Analysis: Real World Appl. {\bf 6}, 273-295 (2005)

\bibitem{WD17}Wang, Z., Du, X.Y.:
 Pullback attractors for modified Swift-Hohenberg equation on unbounded domains with
 non-autonomous deterministic and stochastic forcing terms.
 J. Appl. Anal. Comput. {\bf 7}(1): 207-223 (2017)

\bibitem{XG10}Xiao, Q.K., Gao, H.J.:
 Bifurcation analysis of a modified Swift-Hohenberg equation.
 Nonlinear Analysis: Real World Appl. {\bf 11}, 4451-4464 (2010)

\bibitem{XM15}Xu, L., Ma, Q.Z.:
 Existence of the uniform attractors for a non-autonomous modified Swift-Hohenberg equation.
 Adv. Differ. Equ. {\bf 2015}, 1-11 (2015)


\end{thebibliography}
\end{document}